\documentclass[10pt,a4paper]{article}
\usepackage[utf8]{inputenc}
\usepackage[T1]{fontenc}
\usepackage[english]{babel}
\usepackage[bottom=3cm, top=3cm, left=2cm, right=2cm]{geometry}
\usepackage{amsmath}
\usepackage{amsfonts}
\usepackage{amssymb}
\usepackage{amsthm}
\usepackage{graphicx}
\usepackage{color}
\usepackage{mathrsfs}
\usepackage{lmodern}
\usepackage[most]{tcolorbox}
\usepackage{framed}
\usepackage{fancybox}
\usepackage{tikz}
\usepackage{bbm}
\usepackage{dsfont}
\usepackage{wrapfig}
\usepackage{hyperref}

\newtcolorbox{boxrd}[2][]{enhanced,colback=white,width={\textwidth},
attach boxed title to top left={yshift={-0.5\baselineskip},xshift=1cm}, 
title={#2},
boxrule=0.5pt,
coltitle=black,
boxed title style={
  borderline={-0.5mm}{black}
  colframe=white,
  colback=gray!50,
  colupper={black},
},
}

\newtcolorbox{boxsp}[2][]{%
  enhanced,colback=white,colframe=black,coltitle=black,
  boxrule=0.4pt,
  fonttitle=\itshape,
  attach boxed title to top left={yshift=-0.5\baselineskip-0.3pt,xshift=2mm},
  boxed title style={tile,size=minimal,left=0.5mm,right=0.5mm,
    colback=white,before upper=\strut},
  title=#2,#1
}

\newcommand{\R}{\mathbb{R}}
\newcommand{\C}{\mathbb{C}}

\newcommand{\sech}{\text{sech}}

\usepackage{fancyhdr}
\usepackage{changepage}

\title{\textbf{Asymptotic stability of solitary waves for the 1D near-cubic non-linear Schrödinger equation in the absence of internal modes}}
\author{Guillaume Rialland}
\date{\footnotesize{Université de Paris-Saclay, UVSQ, CNRS, Laboratoire de Mathématiques de Versailles, 78000 Versailles \\ \texttt{guillaume.rialland@uvsq.fr}}}

\begin{document}
\maketitle

\begin{adjustwidth}{80pt}{80pt}
\small{\textsc{Abstract.} We consider perturbations of the one-dimensional cubic Schrödinger equation, under the form $i \, \partial_t \psi + \partial_x^2 \psi + |\psi|^2 \psi - g( |\psi|^2 ) \psi = 0$. Under hypotheses on the function $g$ that can be easily verified in some cases, we show that the linearized problem around a solitary wave does not have internal mode (nor resonance) and we prove the asymptotic stability of these solitary waves, for small frequencies.}
\end{adjustwidth}

\textcolor{white}{a} \\ \\ \textcolor{white}{a}

\noindent We consider the non-linear Schrödinger equation
\begin{equation}
    i \, \partial_t \psi + \partial_x^2 \psi + |\psi|^2 \psi - g( |\psi|^2 ) \psi = 0, \ \, \, \, \, \, \, \, \, (t \, , x) \in \R \times \R,
    \label{NLS}
\end{equation}
which is a perturbation of the cubic NLS equation $i \, \partial_t \psi + \partial_x^2 \psi + | \psi |^2 \psi = 0$. Here, $g \, : \, \R_+ \to \R$ is a function so that the term $g ( |\psi|^2 ) \psi$ is small compared to $| \psi|^2 \psi$ for $| \psi |$ small. We refer to \cite{Pe} or \cite{Ki} for the physical interest of such equations. \\
\\ The corresponding Cauchy problem is globally well-posed in the energy space $H^1 ( \R )$ (see for example \cite{Ca2}) and we recall the Galilean transform, translation and phase invariances of this equation: if $\psi (t \, , x)$ is a solution then, for any $\beta,\sigma,\gamma \in \R$, $\widetilde{\psi} (t \, , x) = e^{i( \beta x - \beta^2 t + \gamma )} \psi (t \, , x-2 \beta t - \sigma )$ is also a solution to the same equation. \\
\\ Solitary waves are solutions of \eqref{NLS} which take the form $\psi (t \, , x) = e^{i \omega t} \phi_\omega (x)$ where 
\begin{equation}
    \phi_\omega '' = \omega \phi_\omega - \phi_\omega^3 + \phi_\omega g(\phi_\omega^2). 
    \label{eqphi}
\end{equation}
It will be proven in the first section below that, under minor hypotheses on $g$ and provided that $\omega$ is small enough, the equation \eqref{eqphi} has a unique solution $\phi_\omega \in H^1 ( \R )$ that is nonnegative, even and that vanishes at infinity. The invariances previously described generate a family of traveling waves given by $\psi (t \, , x) = e^{i ( \beta x - \beta^2 t + \omega t + \gamma )} \phi_\omega (x-2 \beta t - \sigma )$. To begin with, we recall the following standard orbital stability result (see \cite{Ca}, \cite{Gr}, \cite{Il}, \cite{We2}).

\begin{leftbar}
\noindent \textbf{Proposition 1.} For $\omega_0$ small enough and any $\epsilon > 0$, there exists $\delta > 0$ so that, for any $\psi_0 \in H^1 ( \R )$ satisfying $|| \psi_0 - \phi_{\omega_0} ||_{H^1 ( \R )} \leqslant \delta$, if we let $\psi$ be the solution of \eqref{NLS} with initial data $\psi (0) = \psi_0$, then
\[ \sup_{t \in \R} \inf_{( \gamma , \sigma ) \in \R^2} || \psi (t \, , \cdot + \sigma ) - e^{i \gamma} \phi_{\omega_0} ||_{H^1 ( \R )} \leqslant \epsilon.  \]
\end{leftbar}

\noindent In this paper we are interested in the asymptotic stability of solitary waves. There is a vast literature about the asymptotic stability of solitary waves for nonlinear Schrödinger equations, in different cases (various nonlinearities, with or without potential, in different dimensions), see for example \cite{Co}, \cite{Cu1}, \cite{Cu2}, \cite{Ma1} and the review \cite{Ma4}. Before stating our main results, we need to introduce a few hypotheses. First introduce $G (s) = \int_0^s g$. Let us denote $[\![ 0 \, , 5 ]\!] := [0 \, , 5] \cap \mathbb{N}$. Now let us consider the following hypotheses:
\[ \begin{array}{rl} (H_1) & \displaystyle{g \in \mathscr{C}^5 ((0 \, , + \infty )) \cap \mathscr{C}^1 ( [0 \, , + \infty )) \, , \, \, g^{(k)}(s) \, \underset{s \to 0}{=} \, o \left ( s^{1-k} \right ) \, \, \text{for all $k \in [\![ 0 \, , 5 ]\!]$} \, \, \text{and} \, \, g \not\equiv 0 \, \, \text{near $0$},} \\
\\ (H_2) & \displaystyle{\lim_{\omega \to 0} \frac{1}{\varepsilon_\omega^2 \sqrt{\omega}} \int_{\R} \left ( -3 g( \phi_\omega^2 ) + \phi_\omega^2 g'(\phi_\omega^2) + 4 \frac{G(\phi_\omega^2)}{\phi_\omega^2} \right ) \, \text{d}x = + \infty,} \end{array} \]
where $\varepsilon_\omega := \sup\limits_{0 \leqslant s \leqslant 3 \omega} | sg''(s)|$. In this definition, as we shall see in the incoming proofs, $3 \omega$ can be replaced by $2^+ \omega$ where $2^+$ is any constant strictly greater than $2$. Note that the hypothesis $(H_1)$ implies that $\varepsilon_\omega$ exists and is not zero for $\omega>0$ small enough ($\varepsilon_\omega =0$ for $\omega>0$ small would imply that $g'' \equiv 0$ near $0$, thus $g \equiv 0$ since $g(0)=g'(0)=0$). The hypothesis $(H_1)$ also implies that $\varepsilon_\omega \longrightarrow 0$ when $\omega \to 0$. \\
\\ Depending on the function $g$, the equation \eqref{NLS} may (or may not) involve what are called \textit{internal modes}. An internal mode is a solution to the system \eqref{IntM}. It generates periodic solutions to the linearized equation around the solitary wave. For example, $g(s)=s^2$ is a case without internal mode (see the particular study of this case in \cite{Ma1}) while $g(s) = -s^2$ is a case with an internal mode (see \cite{Pe}). In the case $g=0$, there is a resonance (see \cite{Ch}). These considerations justify why we ask for $g \not\equiv 0$ in hypothesis $(H_1)$. The hypothesis $(H_2)$ is a repulsion hypothesis, which involves in particular the sign of the function $g$; the previous remarks let us see that this sign is indeed important. See \cite{Pe}, \cite{Ch} and \cite{CG} for related discussions. Internal modes are potential obstacles to the asymptotic stability of solitons, and we do not address this issue here. We will show that, under the two hypotheses $(H_1)$ and $(H_2)$, there does not exist any internal mode to our problem, in the sense below. Corollary 2 will also assure that there does not exist resonance in this case either. We introduce the following operators, that appear when we linearize \eqref{NLS} around $\phi_\omega$:
\[ L_+ = - \partial_x^2 + \omega - 3 \phi_\omega^2 + g ( \phi_\omega^2) + 2 \phi_\omega^2 g'( \phi_\omega^2 ) \, \, \, \, \, \text{and} \, \, \, \, \, L_- = - \partial_x^2 + \omega - \phi_\omega^2 + g ( \phi_\omega^2 ). \]

\begin{leftbar}
\noindent \textbf{Theorem 1.} Assume that hypotheses $(H_1)$ and $(H_2)$ are satisfied. Then, for $\omega$ small enough, the only solutions $(X \, , Y \, , \lambda ) \in H^1 ( \R )^2 \times \C$ to the system
\begin{equation}
    \left \{ \begin{array}{ccl} L_- X &=& \lambda Y \\ L_+ Y &=& \lambda X \end{array} \right.
    \label{IntM}
\end{equation}
are $X=Y=0$ (and any $\lambda \in \C$) or $\lambda = 0$, $X \in \text{span} ( \phi_\omega )$ and $Y \in \text{span} ( \phi_\omega ')$. 
\end{leftbar}

\noindent Under the same assumptions, we get the following result that ensures the asymptotic stability of the solitons of equation \eqref{NLS}. 

\begin{leftbar}
\noindent \textbf{Theorem 2.} Assume that hypotheses $(H_1)$ and $(H_2)$ are satisfied. For $\omega_0$ small enough, there exists $\delta > 0$ so that, for any $\psi_0 \in H^1 ( \R )$ satisfying $|| \psi_0 - \phi_{\omega_0} ||_{H^1 ( \R )} \leqslant \delta$, if we let $\psi$ be the solution of \eqref{NLS} with initial data $\psi (0) = \psi_0$, then there exists $\beta_+ \in \R$ and $\omega_+ > 0$ such that, for any bounded interval $I \subset \R$,
\[ \lim_{t \to + \infty} \inf_{( \gamma , \sigma ) \in \R^2} \sup_{x \in I} | \psi (t \, , x+ \sigma ) - e^{i \gamma} e^{i \beta_+ x} \phi_{\omega_+} (x) | = 0.  \]
\end{leftbar}

\noindent \textit{Remarks.} A few remarks can be given about this result. Most of them are already in the paper \cite{Ma1} and shall not be recalled here. 
\begin{itemize}
    \item The result is written with an "$\inf\limits_{ \gamma , \sigma}$" formulation. It can be stated in another way, which is the actual way the proof will be led: there exists $\mathscr{C}^1$ functions $\beta , \sigma , \gamma , \gamma \, : \, [0 \, , + \infty ) \mapsto \R^4$ such that $\lim\limits_{t \to + \infty} \beta (t) = \beta_+$, $\lim\limits_{t \to + \infty} \omega (t) = \omega_+$ and
    \[ \lim_{t \to + \infty} \sup_{x \in I} \left | \psi (t \, , x + \sigma (t)) - e^{i \gamma (t)} e^{i \beta (t) x} \phi_{\omega (t)} (x) \right | = 0. \]
    \item The proof will show that $\omega (t) , \omega_+ \in \left ( \frac{\omega_0}{2} \, , \frac{3 \omega_0}{2} \right )$. In fact, we could show that, for any $\eta > 0$, $\delta$ can be chosen small enough such that $\omega (t) , \omega_+ \in \left ( \omega_0 - \eta \, , \omega_0 + \eta \right )$. 
\end{itemize}

\noindent \textcolor{white}{a} \\ The hypothesis $(H_2)$ might appear a little bit cryptic. Let us see how it can be verified in simple cases. Consider for example $g(s) = s^\sigma$ with $\sigma > 1$. We have $sg'(s) = \sigma s^\sigma$, $\frac{G(s)}{s} = \frac{s^\sigma}{\sigma +1}$ and $\varepsilon_\omega = \sigma ( \sigma -1) (3 \omega )^{\sigma -1}$. The hypothesis $(H_1)$ is clearly satisfied. To verify the hypothesis $(H_2)$, we need the following lower bound which will be proved in the first section below: $\phi_\omega (x) \geqslant c \sqrt{\omega} \, e^{- \sqrt{\omega} |x|}$ where $c>0$ does not depend on $\omega$. We see that
\[ \begin{array}{rcl} \displaystyle{\frac{1}{\varepsilon_\omega^2 \sqrt{\omega}} \int_{\R} \left ( -3 g(\phi_\omega^2) + \phi_\omega^2 g'(\phi_\omega^2) + 4 \frac{G(\phi_\omega^2)}{\phi_\omega^2}  \right ) \, \text{d}x} &=& \displaystyle{\frac{1}{\sigma^2 ( \sigma - 1) ( \sigma +1)} \frac{\omega^{- \left ( 2 \sigma - \frac{3}{2} \right )}}{3^{2 \sigma -2}} \int_{\R} \phi_\omega^{2 \sigma}} \\
\\ & \geqslant & \displaystyle{c_\sigma \omega^{-(\sigma -1)} \, \underset{\omega \to 0^+}{\longrightarrow} \, + \infty} \end{array} \]
therefore $(H_2)$ is satisfied. Hence the theorem stated above holds for $g(s) = s^\sigma$ (with $\sigma >1$). \\
\\ Consider a more general situation where $g$ verifies $(H_1)$ and $g''(s) \sim a s^p$ as $s \to 0$, with $a>0$ and $p>-1$. Denote $\sigma := p+2$. Since $\sigma > 1$, $g''(s) \sim as^{\sigma -2}$ leads to $g'(s) \sim \frac{a}{\sigma -1} s^{\sigma -1}$, $g(s) \sim \frac{a}{\sigma (\sigma -1)} s^\sigma$ and $G(s) \sim \frac{a}{( \sigma +1) \sigma ( \sigma - 1)} s^{\sigma+1}$. We get
\[ -3g(s) + sg'(s) + 4 \frac{G(s)}{s} \sim \frac{(\sigma -1)a}{\sigma (\sigma +1)} s^\sigma \, \, \, \, \, \, \, \text{where} \, \, \frac{(\sigma -1)a}{\sigma (\sigma +1)} > 0, \]
which gives $-3 g(s) + sg'(s) + 4 \frac{G(s)}{s} \geqslant \frac{(\sigma - 1)a}{2 \sigma (\sigma+1)} s^\sigma = c_{a , \sigma} s^\sigma$ for $s$ small enough, with $c_{a,\sigma} > 0$. We will see in the first section below that $||\phi_\omega ||_{\infty} \leqslant \sqrt{3 \omega}$ for $\omega$ small enough. Thus, taking $\omega$ small enough, we see that
\[ -3 g(\phi_\omega^2) + \phi_\omega^2 g'(\phi_\omega^2) + 4 \frac{G(\phi_\omega^2)}{\phi_\omega^2} \geqslant c_{a,\sigma} \phi_\omega^{2 \sigma} \geqslant c_{a,\sigma} \omega^{\sigma} e^{-2 \sigma \sqrt{\omega} |x|}. \]
On the other hand, from $g''(s) \sim as^{\sigma -2}$ we deduce that, for $s$ small enough, $|sg''(s)| \leqslant 2 as^{\sigma-1}$ and thus, for $\omega$ small enough, $\varepsilon_\omega \leqslant 2 a (3 \omega )^{\sigma -1} = c_{a,\sigma} \omega^{\sigma -1}$. Gathering these estimates and integrating, we get
\[ \frac{1}{\varepsilon_\omega^2 \sqrt{\omega}} \int_{\R} \left ( -3 g(\phi_\omega^2) + \phi_\omega^2 g'(\phi_\omega^2) + 4 \frac{G(\phi_\omega^2)}{\phi_\omega^2}  \right ) \, \text{d}x \geqslant \frac{c_{a,\sigma}}{\omega^{2 (\sigma-1)} \sqrt{\omega}} \, \omega^\sigma \omega^{-1/2} = c_{a,\sigma} \omega^{1-\sigma} \, \underset{\omega \to 0^+}{\longrightarrow} \, + \infty \]
hence $(H_2)$ is satisfied here too. This case includes functions such as $g(s) = a_1 s^{\sigma_1} + a_2 s^{\sigma_2} + \cdots$ where $1 < \sigma_1 < \sigma_2 < \cdots$, $a_1 > 0$ and $a_i$ (for $i \geqslant 2$) are real numbers whose signs do not matter. \\
\\ We will first prove Theorem 1, which ensures there is no internal mode for our problem. This will be the object of our second part. The third part of this paper is dedicated to the proof of Theorem 2 in itself. The proof extends the one of the analogous result for the case $g(s) = s^2$, which can be found in \cite{Ma1}. It relies on virial arguments, the study of a transformed problem and spectral properties of the linearized operators ($L_+$, $L_-$) and their transformed versions ($M_+$, $M_-$). \\
\\ One can find in \cite{Co} a different approach to the asymptotic stability of the solitons of equation \eqref{NLS}. The functional setting is different, with the use of weighted spaces, and a stronger conclusion about the convergence (often called \textit{full asymptotic stability}). The result of \cite{Co} relies on a natural spectral assumption, namely the non-existence of internal mode and resonance, which was another motivation for Theorem 1 and Corollary 2. Our hypotheses $(H_1)$ and $(H_2)$ and the discussion above thus give concrete situations where the result in \cite{Co} can be applied. \\
\\ The letters $u$, $v$, $w$ and $z$ will denote complex-valued functions; we will index by $1$ their real part and by $2$ their imaginary part (for example, $u=u_1+iu_2$ with $u_1,u_2 \in \R$). The Fourier transform of a function $w$ is denoted by $\widehat{w}$. For $\alpha > 0$, we will use the operator
\[ X_\alpha = (1 - \alpha \partial_x^2 )^{-1} \, \, \, \, \, \, \text{i.e.} \, \, \, \, \, \, \widehat{X_\alpha w} ( \xi ) = \frac{\widehat{w} ( \xi )}{1 + \alpha \xi^2} \, \, \text{for $w \in L^2 ( \R )$.} \]
The $L^2$ scalar product is denoted by $\langle u \, , v \rangle = \text{Re} \left ( \int_{\R} u \overline{v} \, \text{d}x \right )$ and the $L^2$ norm is denoted by $|| \cdot ||$. The $H^1$ norm will be denoted by $|| \cdot ||_{H^1 ( \R )}$. \\
\\ About the virial arguments, we fix a smooth even function $\chi \, : \, \R \to \R$ satisfying $\chi = 2$ on $[0 \, , 1]$, $\chi = 0$ on $[2 \, , + \infty )$ and $\chi ' \leqslant 0$ on $[0 \, , + \infty )$. For $K>0$ we define
\[  \begin{array}{ll} \chi_K (x) = \chi \left ( \frac{x}{K} \right ), & \eta_K (x) = \sech \left ( \frac{2 x}{K} \right ), \\
\\ \zeta_K (x) = \exp \left ( - \frac{|x|}{K} \left ( 1 - \chi ( \sqrt{\omega_0} x) \right ) \right ), & \displaystyle{\Phi_K (x) = \int_0^x \zeta_K(y)^2 \, \text{d}y.} \end{array} \]
We take $A$ and $B$ two large constants that we will fix later (and that depend on $\omega_0$); the idea is to have $A \gg B \gg \frac{1}{\sqrt{\omega_0}} \gg 1$. In everything that follows, $A$ and $B$ are constants (that depend on $\omega_0$) which are assumed to satisfy $A>B>\omega_0^{-1/2} > 1$. Such an inequality will be verified when we indeed fix $A$ and $B$ (in the proof of Proposition 4 for $B$, in the proof of Theorem 2 for $A$). We then define $\Psi_{A,B} = \chi_A^2 \Phi_B$. Most of the bounds we will use and the sketches of the proofs are drawn from \cite{Ma2}, \cite{Ma3}, \cite{Ma1}. Finally we introduce the following weight function
\[ \rho (x) = \sech \left ( \frac{\sqrt{\omega_0}}{10} \, x \right ). \]
Lastly, in this paper, the letter $C$ denotes various positive constants whose expression change from one line to another. The concerned constants do not depend on the parameters $\omega_0$, $\epsilon$, $\alpha$, $A$ and $B$, except in the last part of the proof of Proposition 4, when parameters such as $B$, $\alpha$, $A$ are already fixed. \\
\\ This paper is the result of many discussions with Yvan Martel. The motivation of this paper and its proof are based on his paper \cite{Ma1}. May he be warmly thanked for it here.

\section{Preliminaries}

\subsection{Solitary waves}
\noindent Our proof relies on estimates on the solitons $\phi_\omega$, hence we first have to gather such estimates. The task was easier in the case of the defocusing cubic-quintic NLS equation (see \cite{Ma1}), where solitons were known explicitely. Here, solitons are not know explicitely, but we can prove the following bounds. 

\begin{leftbar}
\noindent \textbf{Lemma 1.} Assume $g$ to be $\mathscr{C}^{5} ( (0 \, , + \infty ))$, $\mathscr{C}^1 ([0 \, , \infty ))$ and such that $g(0)=g'(0)=0$. There exists $\omega_0 > 0$ (depending on $g$) such that, for all $\omega \in (0 \, , \omega_0 )$, there exists a unique solution $\phi_\omega \in H^1 ( \R )$ to the equation $\phi_\omega '' - \omega \phi_\omega + \phi_\omega^3 - g( \phi_\omega^2 ) \phi_\omega = 0$ such that $\phi_\omega$ is even and nonnegative. \\ Moreover, the application $(x \, , \omega) \in \R \times (0 \, , \omega_0) \mapsto \phi_\omega (x)$ is $\mathscr{C}^6$.
\end{leftbar}

\noindent \textit{Proof.} Let us denote $f_\omega (\zeta) = - \omega \zeta + \zeta^3 - g(\zeta^2) \zeta$ and $F_\omega (\zeta) = \int_0^\zeta f_\omega$. We know from \cite{Be} that a solution $\phi_\omega$ verifying all wanted conditions exists if and only if $\zeta_\omega := \inf \{ \zeta > 0 \, \, | \, \, F_\omega ( \zeta ) = 0 \}$ exists and is not zero, and $f_\omega ( \zeta_\omega ) > 0$. In our case, since $g(0)=0$, $f_\omega (\zeta_\omega) > 0$ implies $\zeta_\omega \neq 0$. First, we check that $F_\omega ( \zeta ) = - \frac{\omega \zeta^2}{2} + \frac{\zeta^4}{4} - \frac{G(\zeta^2)}{2}$. By the change of variable $s = \zeta^2$, we have the equivalence
\[ F_\omega ( \zeta )=0 \, \, \, \, \Longleftrightarrow \, \, \, \, \frac{s}{2} - \frac{G(s)}{s} = \omega. \]
Let us denote $J(s) = \frac{s}{2} - \frac{G(s)}{s}$. We take $J(0)$ to be $0$. Indeed, since $g(0)=g'(0)=0$, we have $g(s) = o(s)$ and then $G(s) = o(s^2)$ as $s \to 0^+$. Therefore, $J(s) \sim \frac{s}{2}$. $J$ is clearly $\mathscr{C}^6$ on $(0 \, , + \infty)$ and it is $\mathscr{C}^2$ on $[0 \, , + \infty)$, verifying $J'(0)=\frac{1}{2}$, $J''(0)=0$. Since $J'(0) \neq 0$, by local inversion we know that there exists $s_0 > 0$ such that $J$ is bijective from $[0 \, , s_0]$ to $[0 \, , J( s_0 )]$. Taking $\omega_0 = J(s_0)$, it is now clear that, for every $\omega \in (0 \, , \omega_0)$, there exists a unique $s_\omega \in (0 \, , s_0)$ such that $J(s_\omega) = \omega$. The uniqueness shows that $\zeta_\omega = \sqrt{s_\omega}$ is the quantity $\inf \{ \zeta > 0 \, \, | \, \, F_\omega ( \zeta ) = 0 \}$ we look for. \\
\\ Now, $f_\omega ( \zeta_\omega ) = \zeta_\omega ( - \omega - g(s_\omega) + s_\omega)$. We aim to prove that this is positive. First, we have $J \left ( \frac{3 \omega}{2} \right ) = \frac{3 \omega}{4} - \frac{G(3 \omega /2)}{3 \omega /2}$. Since $G(s) = o(s^2)$, $J \left ( \frac{3 \omega}{2} \right ) \sim \frac{3 \omega}{4}$ as $\omega \to 0$; thus we can take a smaller $\omega_0$ to be sure that $J \left ( \frac{3 \omega}{2} \right ) < \omega$ for all $\omega \in (0 \, , \omega_0 )$. From now on we make that assumption. This proves that, for all $\omega \in (0 \, , \omega_0 )$, we have $\frac{3 \omega}{2} < s_\omega$. \\
\\ Since $g (s) = o(s)$, we can assume that $|g(s)| \leqslant \frac{s}{3}$ for all $s \in [0 \, , s_1]$. On the foregoing, we may have assumed that $s_0 \leqslant s_1$. From now on we make that assumption. Now, we can check that, for all $\omega \in (0 \, , \omega_0 )$,
\[ g(s_\omega) \leqslant \frac{s_\omega}{3} \, \, \, \, \, \text{thus} \, \, \, \, \, - \omega - g(s_\omega) + s_\omega \geqslant - \omega - \frac{s_\omega}{3} + s_\omega = \frac{2 s_\omega}{3} - \omega > 0 \]
as we have seen. This shows that $f_\omega ( \zeta_\omega ) > 0$ and completes the first part of the lemma. \\
\\ The regularity of the function $(x \, , \omega) \mapsto \phi_\omega (x)$ comes from standard arguments. We recall from \cite{Be} that the solution $\phi_\omega$ is the only solution of the Cauchy problem
\[ \left \{ \begin{array}{l} \phi_\omega '' - \omega \phi_\omega + \phi_\omega^3 - g( \phi_\omega^2 ) \phi_\omega = 0, \\ \phi_\omega (0) = \zeta_\omega , \, \, \, \phi_\omega '(0)=0. \end{array} \right. \]
We have to check that $\omega \mapsto \zeta_\omega$ is $\mathscr{C}^6$ on $(0 \, , \omega_0 )$. This is the case since $\omega \mapsto \zeta_\omega$ is nothing else than $\sqrt{J^{-1}}$ and that $J$ is $\mathscr{C}^6$ on $(0 \, , \omega_0 )$. Note that $J$ itself is possibly not $\mathscr{C}^6$ near $0$. That is not a problem, since the solutions $\phi_\omega$ take their values in $(0 \, , + \infty )$; hence $(0 \, , + \infty )$ is the arrival domain of the Cauchy-Lipschitz theorem with parameter we apply. We then get the $\mathscr{C}^6$ regularity we seek. \hfill \qedsymbol

\textcolor{white}{a} \\ \noindent The hypotheses above about $g$ will always be assumed: they are implied by hypothesis $(H_1)$. We have $G(s) = o(s^2)$ and thus $\omega = \frac{\zeta_\omega^2}{2} - \frac{G(\zeta_\omega^2)}{\zeta_\omega^2} = \frac{\zeta_\omega^2}{2} + o(\zeta_\omega^2)$. Hence, $\zeta_\omega^2 \sim 2 \omega$ i.e. $\zeta_\omega \sim \sqrt{2 \omega}$. We will suppose in the whole paper that $\omega$ is chosen small enough so that $\zeta_\omega \leqslant \sqrt{3 \omega}$. We also suppose that $\omega$ is chosen small enough so that $|g(s)| < s$ for any $s \in [0 \, , 3 \omega ]$. Moreover, we will need an equivalent of $\frac{\text{d} \zeta_\omega}{\text{d} \omega}$. Recalling that $( \omega \mapsto \zeta_\omega ) = \sqrt{J^{-1}}$, we write that
\[ \frac{\text{d} \zeta_\omega}{\text{d} \omega} = \frac{1}{2 J'(J^{-1} ( \omega )) \sqrt{J^{-1} ( \omega )}} = \frac{1}{\zeta_\omega \left ( 1 - \frac{2 g(\zeta_\omega^2)}{\zeta_\omega} + \frac{2 G ( \zeta_\omega^2)}{\zeta_\omega^3} \right )} \sim \frac{1}{\zeta_\omega} \sim \frac{1}{\sqrt{2 \omega}} \]
since $\frac{g(\zeta_\omega^2)}{\zeta_\omega} = o ( \zeta_\omega ) = o(1)$ and $\frac{G(\zeta_\omega^2)}{\zeta_\omega^3} = o(\zeta_\omega) = o(1)$. \\
\\ In what follows, we always take $\omega \in (0 \, , \omega_0)$. We drop the notation $\omega_0$ and only say that $\omega$ is \textit{"small enough"}. We might have to reduce the range to which $\omega$ belongs in what follows, which is not a problem. Let $Q_\omega$ be the solitary-wave solution of the cubic Schrödinger stationary equation $Q_\omega '' - \omega Q_\omega + Q_\omega^3 = 0$. That is to say, $Q_\omega$ corresponds to the case $g=0$. We know $Q_\omega$ explicitly: denoting $Q(x) = \frac{\sqrt{2}}{\cosh (x)}$, $Q_\omega$ is given by $Q_\omega (x) = \sqrt{\omega} \, Q ( \sqrt{\omega} \, x)$. We can guess that $\phi_\omega$ has growth properties that are similar to $Q_\omega$. This is the object of the following lemma. Besides, since $\phi_\omega$ is $\mathscr{C}^6$ with regards to $\omega$ (provided $(H_1)$ is satisfied), it makes sense to consider $\Lambda_\omega := \omega \, \frac{\partial \phi_\omega}{\partial \omega}$ and we know that $\Lambda_\omega$ is the solution on $\R$ of the following Cauchy system
\[ \left \{ \begin{array}{l} - \Lambda_\omega '' = - \omega \phi_\omega - \omega \Lambda_\omega + 3 \phi_\omega^2 \Lambda_\omega - 2 \phi_\omega^2 g'(\phi_\omega^2) \Lambda_\omega - g(\phi_\omega^2) \Lambda_\omega \\ \\ \Lambda_\omega (0) = \omega \frac{\text{d} \zeta_\omega}{\text{d} \omega} \sim \sqrt{\frac{\omega}{2}} , \, \, \, \Lambda_\omega '(0) = 0, \end{array} \right. \]
where we recognise the first line to be $L_+ \Lambda_\omega = - \omega \phi_\omega$. Controlling $\Lambda_\omega$ and its derivative will be the object of Lemma 5. 

\begin{leftbar}
\noindent \textbf{Lemma 2.} Assume $g$ to be $\mathscr{C}^{5} ( (0 \, , + \infty ))$, $\mathscr{C}^1 ([0 \, , \infty ))$ and such that $g(0)=g'(0)=0$. For any $k \in [\![ 0 \, , 6 ]\!]$, there exists $C_k > 0$ such that, for any $\omega > 0$ small enough and any $x \in \R$,
\[ | \phi_\omega^{(k)} (x) | \leqslant C_k \omega^{\frac{1+k}{2}} e^{- \sqrt{\omega} |x|}. \]
Moreover, for every $\varepsilon > 0$, for any $\omega > 0$ small enough, 
\[ | \phi_\omega (x) - Q_\omega (x) | \leqslant \varepsilon \sqrt{\omega} \, e^{- \sqrt{\omega} |x|}.  \]
Lastly, there exists $c>0$ such that $\phi_\omega (x) \geqslant c \sqrt{\omega} \, e^{- \sqrt{\omega} |x|}$. 
\end{leftbar}

\noindent \textit{Proof.} This proof will require several steps and is based on standard ordinary differential equations arguments that can be found in \cite{Be}. We will denote $P_\omega = \phi_\omega - Q_\omega$. Let $\varepsilon > 0$. Let us take $x_0 > 0$ such that $Q(x) < \varepsilon$ for $x \geqslant x_0$ ($x_0$ does not depend on $\omega$). Now, for $x \geqslant x_0 / \sqrt{\omega}$, $Q_\omega (x) < \varepsilon \sqrt{\omega}$. Considering the equations satisfied by $\phi_\omega$ and $Q_\omega$, we get
\[ P_\omega '' - \omega P_\omega = - P_\omega ( Q_\omega^2 + \phi_\omega Q_\omega + \phi_\omega^2 ) + g ( \phi_\omega^2 ) \phi_\omega. \]
It is clear that $0 \leqslant Q_\omega (x) \leqslant \sqrt{2 \omega}$ for all $x \in \R$. Now, since $\phi_\omega$ is nonincreasing on $\R_+$ and even, $0 \leqslant \phi_\omega (x) \leqslant \phi_\omega (0) = \zeta_\omega \leqslant C \sqrt{\omega}$. Thus, we get
\[  | P_\omega '' | \leqslant \omega | P_\omega | + 2 \left ( |Q_\omega|^2 + | \phi_\omega |^2 \right ) |P_\omega| + \varepsilon | \phi_\omega |^3 \leqslant C \omega | P_\omega | + C \varepsilon \omega^{3/2}. \]
Considering the vectorial function $\overrightarrow{\textbf{P}}_\omega (x) = (P_\omega (x) \, , P_\omega '(x) / \sqrt{\omega} )$, we have $|| \overrightarrow{\textbf{P}}_\omega ' (x) ||_{1} \leqslant C \sqrt{\omega} || \overrightarrow{\textbf{P}}_\omega (x) ||_1 + C \varepsilon \omega$ where $|| (p_1 \, , p_2)^\top ||_1 := |p_1| + |p_2|$. We then use Grönwall's lemma and the fact that $P_\omega '(0)=0$ to see that
\[ |P_\omega (x)| \leqslant || \overrightarrow{\textbf{P}}_\omega ||_1 \leqslant - C \varepsilon \sqrt{\omega} + (P_\omega (0) + \varepsilon \sqrt{\omega} ) e^{C \sqrt{\omega} x}. \]
As $|P_\omega (0)| = | \zeta_\omega - \sqrt{2 \omega} | = o ( \sqrt{\omega} )$, we get that, $|P_\omega (x_0 / \sqrt{\omega})| \leqslant C \varepsilon \sqrt{\omega}$ and thus $| \phi_\omega (x_0 / \sqrt{\omega} ) | \leqslant C \varepsilon \sqrt{\omega}$. Now, $\phi_\omega$ being nonincreasing, we get, for any $x \geqslant x_0 / \sqrt{\omega}$, that $|\phi_\omega (x)| \leqslant C \varepsilon \sqrt{\omega}$. Now, let us use standard arguments from \cite{Be}. Setting $v_\omega (x) = \phi_\omega (x)^2$, we have, for any $x \geqslant x_0 / \sqrt{\omega}$,
\[ v_\omega ''(x) = 2 \phi_\omega '(x)^2 + 2 \left ( \omega - \phi_\omega (x)^2 + g ( \phi_\omega (x)^2) \right ) v_\omega (x) \geqslant 2 \left ( \omega - 4 C \varepsilon^2 \omega - 4 C \varepsilon^2 \omega \right ) v_\omega (x) \geqslant \omega v_\omega (x) \]
providing we take $\varepsilon$ small enough so that $1 - 8 C \varepsilon^2 > \frac{1}{2}$. Now taking $z_\omega (x) = e^{- \sqrt{\omega}x} (v_\omega '(x) + \sqrt{\omega} \, v_\omega (x))$, we have $z_\omega '(x) = e^{- \sqrt{\omega}x} (v_\omega ''(x) - \omega v_\omega (x)) \geqslant 0$ for $x \geqslant x_0 / \sqrt{\omega}$. Therefore $z$ is nondecreasing on $\left [ \frac{x_0}{\sqrt{\omega}} \, , + \infty \right )$. Suppose that $z_\omega (y) > 0$ for some $y > x_0 / \sqrt{\omega}$. Then, for all $x \geqslant y$, $z_\omega (x) \geqslant z_\omega (y) > 0$ thus $v_\omega '(x) + \sqrt{\omega} \, v_\omega (x) \geqslant z_\omega (y) e^{\sqrt{\omega}x}$, showing that $v_\omega + \sqrt{\omega} \, v_\omega \not\in L^1 ([y \, , + \infty ))$. However we know that $\phi_\omega \in H^1 ( \R )$, hence $\phi_\omega , \phi_\omega ' \in L^2 ( \R )$ and $v_\omega = \phi_\omega^2 \in L^1 ( \R )$ and $v_\omega ' = 2 \phi_\omega \phi_\omega ' \in L^1 ( \R )$ too. Finally, this is absurd: we conclude that $z_\omega$ remains nonpositive for all $x \geqslant x_0 / \sqrt{\omega}$. This shows that $x \mapsto e^{\sqrt{\omega} x} v_\omega (x)$ is nonincreasing on $\left [ \frac{x_0}{\sqrt{\omega}} \, , + \infty \right )$ and then 
\[  \forall x \geqslant \frac{x_0}{\sqrt{\omega}} , \, \, \, \, 0 \leqslant v_\omega (x) \leqslant e^{x_0} v_\omega \left ( \frac{x_0}{\sqrt{\omega}} \right ) e^{- \sqrt{\omega} x}. \]
Since $v_\omega \left ( \frac{x_0}{\sqrt{\omega}} \right ) \leqslant 4 \varepsilon^2 \omega$, we finally get that $v_\omega (x) \leqslant C \varepsilon^2 \omega e^{- \sqrt{\omega} x}$ and thus $\phi_\omega (x) \leqslant C \varepsilon \sqrt{\omega} \, e^{- \frac{\sqrt{\omega}}{2} x}$ for any $x \geqslant x_0 / \sqrt{\omega}$. \\
\\ Now we see that, by the variation of the constants, using the initial conditions $\phi_\omega (0) = \zeta_\omega$ and $\phi_\omega '(0) = 0$,
\[ \phi_\omega (x) = \frac{\zeta_\omega}{2} \, e^{\sqrt{\omega} x} + \frac{\zeta_\omega}{2} \, e^{- \sqrt{\omega} x} + \frac{e^{\sqrt{\omega} x}}{2 \sqrt{\omega}} \int_0^x \ell_\omega (y) e^{- \sqrt{\omega} y} \, \text{d}y - \frac{e^{ - \sqrt{\omega} x}}{2 \sqrt{\omega}} \int_0^x \ell_\omega (y) e^{\sqrt{\omega} y} \, \text{d}y \]
where $\ell_\omega (y) = - \phi_\omega (y)^3 + g(\phi_\omega (y)^2) \phi_\omega (y)$. We introduce $I_\omega^{\pm} = \int_0^{\infty} \ell_\omega (y) e^{ \pm \sqrt{\omega} y} \, \text{d}y$. Both of these integrals indeed converge, as $|\ell_\omega (y)| \leqslant C \omega^{3/2} e^{- 3 \sqrt{\omega} y/2}$ when $y \to \infty$. We then write $\phi_\omega (x)$ as
\[ \phi_\omega (x) = \left ( \frac{\zeta_\omega}{2} + \frac{I_\omega^-}{2 \sqrt{\omega}} \right ) e^{\sqrt{\omega} x} + \left ( \frac{\zeta_\omega}{2} - \frac{I_\omega^+}{2 \sqrt{\omega}} \right ) e^{- \sqrt{\omega} x} - \frac{e^{\sqrt{\omega} x}}{2 \sqrt{\omega}} \int_x^{\infty} \ell_\omega (y) e^{- \sqrt{\omega} y} \, \text{d}y + \frac{e^{- \sqrt{\omega} x}}{2 \sqrt{\omega}} \int_x^{\infty} \ell_\omega (y) e^{\sqrt{\omega} y} \, \text{d}y . \]
Since $\phi_\omega (x) \longrightarrow 0$ as $x \to + \infty$, $\frac{\zeta_\omega}{2} + \frac{I_\omega^-}{2 \sqrt{\omega}} = 0$ and we get the following expression:
\[ \phi_\omega (x) = \left ( \frac{\zeta_\omega}{2} - \frac{I_\omega^+}{2 \sqrt{\omega}} \right ) e^{- \sqrt{\omega} x} - \frac{e^{\sqrt{\omega} x}}{2 \sqrt{\omega}} \int_x^{\infty} \ell_\omega (y) e^{- \sqrt{\omega} y} \, \text{d}y + \frac{e^{- \sqrt{\omega} x}}{2 \sqrt{\omega}} \int_x^{\infty} \ell_\omega (y) e^{\sqrt{\omega} y} \, \text{d}y . \]
Separating the integral $I_\omega$ at $x_0 / \sqrt{\omega}$ and using respectively the control $\phi_\omega (y) \leqslant C \sqrt{\omega} e^{- \sqrt{\omega} y/2}$ if $y \geqslant x_0 / \sqrt{\omega}$, and the control $\phi_\omega (y) \leqslant C \sqrt{\omega}$ if $0 \leqslant y < x_0 / \sqrt{\omega}$, we get that $|I_\omega| \leqslant C \omega$. Hence $\left | \frac{\zeta_\omega}{2} - \frac{I_\omega^+}{2 \sqrt{\omega}} \right | \leqslant C \sqrt{\omega}$. \\
\\ About the integral terms, we shall separate the integral at the point $x_0 / \sqrt{\omega}$ too. If $x \geqslant x_0 / \sqrt{\omega}$, there is no need to separate: the upper bound $\phi_\omega (y) \leqslant C \sqrt{\omega} \, e^{- \sqrt{\omega} y/2}$ directly gives $\left | \int_x^\infty \ell_\omega (y) e^{- \sqrt{\omega} y} \, \text{d}y \right | \leqslant C \omega e^{-5 \sqrt{\omega} x/2}$. If $0 \leqslant x < x_0 / \sqrt{\omega}$, we separate the integral and use the same upper bounds as for $I_\omega$; we get
\[ \left | \int_x^{\infty} \ell_\omega (y) e^{- \sqrt{\omega} y} \, \text{d}y \right | \leqslant C \omega \left ( e^{- \sqrt{\omega} x} - e^{-x_0} \right ) + C \omega e^{-5 x_0/2} \leqslant C \omega . \]
We then get
\[ \left | \frac{e^{\sqrt{\omega} x}}{2 \sqrt{\omega}} \int_x^\infty \ell_\omega (y) e^{- \sqrt{\omega} y} \, \text{d}y \right | \leqslant C \sqrt{\omega} \, e^{\sqrt{\omega} x} \leqslant C \sqrt{\omega} \, e^{x_0} = C \sqrt{\omega} \leqslant C \sqrt{\omega} \, e^{- \sqrt{\omega} x} \]
thanks to the lower bound $e^{- \sqrt{\omega} x} \geqslant e^{-x_0}$. In the lines above, the important fact is that the constant $C$ (which changes from one expression to another) does not depend on $\omega$. We thus have proved that, for any $x \geqslant 0$,
\[ \left | \frac{e^{\sqrt{\omega} x}}{2 \sqrt{\omega}} \int_x^\infty \ell_\omega (y) e^{- \sqrt{\omega} y} \, \text{d}y \right | \leqslant C \sqrt{\omega} \, e^{- \sqrt{\omega} x}. \]
The reasoning is exactly the same for the second integral: a direct exponential control when $x \geqslant x_0 / \sqrt{\omega}$ thanks to the previous upper bound, and a bounded control when $x < x_0 / \sqrt{\omega}$, which is sufficient for our purpose. Finally, we get that, for any $x \in \R$, $|\phi_\omega (x)| \leqslant C \sqrt{\omega} \, e^{- \sqrt{\omega} |x|}$ where $C$ does not depend on $\omega$. \\
\\ The estimates on the derivatives of $\phi_\omega$ follow from the expression obtained previously. We indeed have
\[ \phi_\omega '(x) = \left ( - \sqrt{\omega} \zeta_\omega - \frac{I_\omega}{2} \right ) e^{- \sqrt{\omega} x} - \frac{e^{\sqrt{\omega}x}}{2} \int_x^\infty \ell_\omega (y) e^{- \sqrt{\omega} y} \, \text{d}y - \frac{e^{- \sqrt{\omega} x}}{2} \int_x^\infty \ell_\omega (y) e^{\sqrt{\omega} y} \, \text{d}y. \]
With the bounds shown above about the integral terms, we get $|\phi_\omega '(x)| \leqslant C \omega e^{- \sqrt{\omega} x}$ with the same proof. To control $\phi_\omega ''$ and further derivatives, we use the equation satisfied by $\phi_\omega$; the conclusion follows. \\
\\ Now let us prove the bound on $P_\omega = \phi_\omega - Q_\omega$. To start, let us prove that $||P_\omega||_{\infty} = o ( \sqrt{\omega} )$ as $\omega \to 0$. Let $\varepsilon >0$. We know, from the exponential decays of $\phi_\omega$ and $Q_\omega$, that $|P_\omega (x)| \leqslant C \sqrt{\omega} \, e^{- \sqrt{\omega} x}$ for all $x \in \R$. Let us take $\omega$ sufficiently small so that $\zeta_\omega \leqslant \sqrt{3 \omega}$, $|g(s)| \leqslant \delta_1 s$ for all $s \in [0 \, , 3 \omega ]$, and finally $| \zeta_\omega - \sqrt{2 \omega} | \leqslant \delta_2 \sqrt{\omega}$; where we have denoted $\delta_1 = \frac{\varepsilon}{4} \, e^{-12 \ln (C / \varepsilon )}$ and $\delta_2 = \frac{\varepsilon}{2} \, e^{-12 \ln (C / \varepsilon )}$. The previous lines imply that $\phi_\omega \leqslant \sqrt{3 \omega}$, $g ( \phi_\omega^2 ) \leqslant \delta_1 \phi_\omega^2 \leqslant 3 \delta_1 \omega$ and $|P_\omega (0)| \leqslant \delta_2 \sqrt{\omega}$. We then have, thanks to the equation $P_\omega '' - \omega P_\omega = - P_\omega (Q_\omega^2 + \phi_\omega Q_\omega + \phi_\omega^2 ) + g ( \phi_\omega^2 ) \phi_\omega$ verified by $P_\omega$,
\[ | P_\omega '' | \leqslant \omega | P_\omega | + (2 \omega + \sqrt{6} \omega + 3 \sqrt{\omega} ) |P_\omega | + 3 \sqrt{3} \delta_1 \omega^{3/2} \leqslant 12 \omega | P_\omega | + 6 \delta_1 \omega^{3/2}. \]
Let $x_\omega = \frac{\ln (C / \varepsilon)}{\sqrt{\omega}}$, such that $|P_\omega (x)| \leqslant C \sqrt{\omega} \, e^{- \sqrt{\omega} x_\omega} = \varepsilon \sqrt{\omega}$ for all $x \geqslant x_\omega$. Following the same computation as in the first part of the proof, we get by Gronwäll's lemma that, for all $x \in [0 \, , x_\omega]$,
\[ |P_\omega (x)| \leqslant \sqrt{\omega} \left [ \frac{\delta_1}{2} + e^{12 \sqrt{\omega} x} \left ( \delta_2 + \frac{\delta_1}{2} \right ) \right ] \leqslant \sqrt{\omega} \left [ \frac{\delta_1}{2} + e^{12 \sqrt{\omega} x_\omega} \left ( \delta_2 + \frac{\delta_1}{2} \right ) \right ] = \varepsilon \]
thanks to the judicious choices of $\delta_1$ and $\delta_2$. Therefore, we have proved that $|P_\omega (x)| \leqslant \varepsilon \sqrt{\omega}$ for all $x \geqslant 0$ and all $\omega >0$ small enough. \\
\\ Now, the variation of the constants and the fact that $P_\omega '(0)=0$ give the expression
\[ P_\omega (x) = A \, e^{\sqrt{\omega} x} + \left ( \frac{P_\omega (0)}{2} - \frac{J_\omega}{2 \sqrt{\omega}} \right ) e^{- \sqrt{\omega} x} - \frac{e^{\sqrt{\omega} x}}{2 \sqrt{\omega}} \int_x^{+ \infty} S_\omega (y) e^{- \sqrt{\omega} y} \, \text{d}y + \frac{e^{-\sqrt{\omega} x}}{2 \sqrt{\omega}} \int_x^{+ \infty} S_\omega (y) e^{\sqrt{\omega} y} \, \text{d}y \]
where $S_\omega = Q_\omega^3 - \phi_\omega^3 + g ( \phi_\omega^2 ) \phi_\omega = - P_\omega (Q_\omega^2 + \phi_\omega Q_\omega + \phi_\omega^2 ) + g ( \phi_\omega^2 ) \phi_\omega$ and $J_\omega = \int_0^{+ \infty} S_\omega (y) e^{\sqrt{\omega} y} \, \text{d}y$. Taking $\omega$ sufficiently small so that $g(s) \leqslant \varepsilon s$ for all $s \in [0 \, , 3 \omega ]$ and $|\zeta_\omega| \leqslant \sqrt{3 \omega}$, and also using the inequality $||P_\omega||_\infty \leqslant \varepsilon \sqrt{\omega}$ we have just proved, we find that $| S_\omega (y) | \leqslant C \varepsilon \omega^{3/2} e^{-2 \sqrt{\omega}y}$ for all $y \geqslant 0$. This gives 
\[ |J_\omega| \leqslant C \varepsilon \omega , \, \, \, \, \, \, \, \left | \int_x^{+ \infty} S_\omega (y) e^{- \sqrt{\omega} y} \, \text{d}y \right | \leqslant C \varepsilon \omega e^{-3 \sqrt{\omega} x} \, \, \, \, \, \, \, \text{and} \, \, \, \, \, \, \, \left | \int_x^{+ \infty} S_\omega (y) e^{\sqrt{\omega} y} \, \text{d}y \right | \leqslant C \varepsilon \omega e^{- \sqrt{\omega} x}. \]
Since $P_\omega$ vanishes at infinity, this shows that $A=0$, and gathering all the upper bounds we get that $| P_\omega (x) | \leqslant C \varepsilon \sqrt{\omega} e^{- \sqrt{\omega} x}$. \\
\\ For the last bound, we know from the explicit expression of $Q_\omega$ that $Q_\omega (x) \geqslant c \sqrt{\omega} \, e^{- \sqrt{\omega} |x|}$. Taking $\varepsilon$ small enough in $| \phi_\omega (x) - Q_\omega (x) | \leqslant \varepsilon \sqrt{\omega} \, e^{- \sqrt{\omega} |x|}$, we obtain the desired lower bound: $\phi_\omega (x) \geqslant c \sqrt{\omega} \, e^{- \sqrt{\omega} |x|}$. \hfill \qedsymbol

\textcolor{white}{a} \\ \noindent We recall that the linearization of \eqref{NLS} around $\phi_\omega$ involves the operators
\[ L_+ = - \partial_x^2 + \omega - 3 \phi_\omega^2 + g ( \phi_\omega^2) + 2 \phi_\omega^2 g'( \phi_\omega^2 ) \, \, \, \, \, \text{and} \, \, \, \, \, L_- = - \partial_x^2 + \omega - \phi_\omega^2 + g ( \phi_\omega^2 ) \]
as we can see in \cite{We2} for instance. Some spectral properties are known about $L_+$ and $L_-$ (see \cite{We1}). Both operators are self-adjoint in $L^2$. In $L^2$, the operator $L_+$ has exactly one negative eigenvalue and its kernel is generated by $\phi_\omega '$. On the other hand, still in $L^2$, the kernel of $L_-$ is generated by $\phi_\omega$. \\
\\ Let us discuss some aspects about the invertibility of $L_+$, which will intervene in the last part of our proof. The invertibility of $M_-$ also intervenes at the same point; it is discussed in section 1.3. We denote by $A_\omega$ the even solution of $L_+ A_\omega = 0$ such that $\phi_\omega '' A_\omega - \phi_\omega ' A_\omega ' = 1$ on $\R$. The variation of the constants shows that, if $A_\omega$ was bounded, then we would have $A_\omega (x) , A_\omega '(x) \, \underset{x \to + \infty}{\longrightarrow} \, 0$, which clearly contradicts the relation $\phi_\omega '' A_\omega - \phi_\omega ' A_\omega ' = 1$; thus $A_\omega$ is not bounded on $\R$. We will need the following estimate on $A_\omega$. 

\begin{leftbar}
\noindent \textbf{Lemma 3.} Assume that hypothesis $(H_1)$ holds. For $\omega > 0$ small enough and for any $k \in [\![ 0 \, , 6 ]\!]$, there exists some constants $C_k > 0$ such that $|A_\omega^{(k)} (x)| \leqslant C_k \omega^{\frac{k-3}{2}} e^{\sqrt{\omega} |x|}$ for all $x \in \R$.
\end{leftbar}

\noindent \textit{Proof.} Starting with the wronskian relation, we have $A_\omega ' - \frac{\phi_\omega ''}{\phi_\omega '} \, A_\omega = - \frac{1}{\phi_\omega '}$ on $(0 \, , + \infty )$ and thus we get
\[ A_\omega (x) = \phi_\omega '(x) \left [ \alpha_\omega + \int_x^{1/ \sqrt{\omega}} \frac{\text{d}y}{\phi_\omega '(y)^2} \right ] \]
where $\alpha_\omega$ is an unknown constant (that depends on $\omega$). Now, let us define $\text{res}_\omega (x) := \frac{1}{\phi_\omega '(x)^2} - \frac{1}{\phi_\omega ''(0)^2 x^2}$. Using $\phi_\omega '(x) = \phi_\omega ''(0) x + \mathcal{O} (x^3)$ as $x \to 0$, we see that $\text{res}_\omega (x) = \mathcal{O} (1)$. Differentiating the expression of $A_\omega$ above, we find that
\[ A_\omega ' (x) \, \underset{x \to 0}{=} \, - \frac{\sqrt{\omega}}{\phi_\omega ''(0)} + \left ( \alpha_\omega + \int_0^{1/ \sqrt{\omega}} \text{res}_\omega \right ) \phi_\omega ''(0) + o(1). \]
Since $A_\omega$ is even, $A_\omega '(0) = 0$ thus $\alpha_\omega = \frac{\sqrt{\omega}}{\phi_\omega ''(0)^2} - \int_0^{1/\sqrt{\omega}} \text{res}_\omega$. \\
\\ Now let us take $\varepsilon >0$ and introduce $D_\omega := P_\omega ' = \phi_\omega ' - Q_\omega '$ where we recall that $P_\omega = \phi_\omega - Q_\omega$. We see that $D_\omega ' = \omega P_\omega - P_\omega (Q_\omega^2 + \phi_\omega Q_\omega + \phi_\omega^2) + g(\phi_\omega^2) \phi_\omega$. Using the estimates of Lemma 2 we obtain, for $\omega >0$ small enough, $|D_\omega'(x)| \leqslant C \varepsilon \omega^{3/2} e^{- \sqrt{\omega} x}$ for all $x>0$. For $x > \omega^{-1/2}$, we get
\[ |D_\omega(x)| \leqslant \int_x^{+ \infty} C \varepsilon \omega^{3/2} e^{- \sqrt{\omega} y} \, \text{d}y \leqslant C \varepsilon \omega e^{- \sqrt{\omega} x} \leqslant C \varepsilon \omega^{3/2} x e^{- \sqrt{\omega} x}, \]
and for $0 < x < \omega^{-1/2}$ we get
\[ |D_\omega(x)| \leqslant \int_0^x C \varepsilon \omega^{3/2} e^{- \sqrt{\omega} y} \, \text{d}y \leqslant C \varepsilon \omega (e^{- \sqrt{\omega} x} -1) \leqslant C \varepsilon \omega \sqrt{\omega} x \leqslant C \varepsilon \omega^{3/2} x e^{- \sqrt{\omega} x}. \]
Thus, $|D_\omega(x)| \leqslant C \varepsilon \omega^{3/2} x e^{- \sqrt{\omega} x}$ for all $x>0$. Note that we have used the fact that $D_\omega(0)=0$. Also note that it is also true that $|D_\omega(x)| \leqslant C \varepsilon \omega e^{- \sqrt{\omega} x}$ for all $x>0$. Now, using the explicit expression $Q_\omega ' (x) = - \sqrt{2} \, \omega \, \frac{\text{sinh} ( \sqrt{\omega} \, x)}{\text{cosh}^2 ( \sqrt{\omega} \, x)}$, we see that $| Q_\omega '(x) | \geqslant C \omega^{3/2} x$ for $x \in (0 \, , \omega^{-1/2} )$ and that $| Q_\omega '(x)| \geqslant C \omega e^{- \sqrt{\omega} |x|}$ for all $x \in \R$. This shows that, for all $x \in \R$, $| \phi_\omega '(x) | \geqslant C (1 - C' \varepsilon ) | Q_\omega '(x)|$. Choosing $\varepsilon >0$ correctly, we obtain $| \phi_\omega '(x) | \geqslant C | Q_\omega '(x) | \geqslant C \omega e^{- \sqrt{\omega} |x|}$ for all $x \in \R$. For $x \in (0 \, , \omega^{-1/2})$, this leads to $\phi_\omega '(x)^2 \geqslant C Q_\omega'(x)^2 \geqslant C \omega^3 x^2$. \\
\\ On the other hand, differentiating four times the quantity $(\phi_\omega ')^2$ thanks to \eqref{eqphi} and using Lemma 2, we easily see that, for $\omega >0$ small enough and all $x>0$,
\[ \left | \frac{\text{d}^4}{\text{d}x^4} \left ( \phi_\omega ''(0)^2 x^2 - \phi_\omega '(x)^2 \right ) \right | \leqslant C \omega^4. \]
Since the function $x \mapsto  \phi_\omega ''(0)^2 x^2 - \phi_\omega'(x)^2$ and its first three derivatives vanish at $x=0$, we obtain $| \phi_\omega ''(0)^2 x^2 - \phi_\omega '(x)^2 | \leqslant C \omega^4 x^4$ for all $x>0$. \\
\\ Finally, using \eqref{eqphi}, we see that $\phi_\omega ''(0) \sim - \sqrt{2} \, \omega^{3/2}$ as $\omega \to 0$. Thus, for $\omega > 0$ small enough, $\phi_\omega ''(0)^2 \geqslant C \omega^3$. Putting these estimates together, we find that, for $x \in (0 \, , \omega^{-1/2})$,
\[ | \text{res}_\omega (x) | = \frac{| \phi_\omega ''(0)^2 x^2 - \phi_\omega'(x)^2 |}{\phi_\omega '(x)^2 \phi_\omega ''(0)^2 x^2} \leqslant \frac{C \omega^4 x^4}{C \omega^3 x^2 \cdot C \omega^3 x^2} \leqslant C \omega^{-2}. \]
Integrating on $(0 \, , \omega^{-1/2})$ and recalling that $\phi_\omega ''(0)^{-2} \leqslant C \omega^{-3}$, we obtain $| \alpha_\omega | \leqslant C \omega^{-5/2}$. Thus $| \alpha_\omega \phi_\omega '(x)| \leqslant C \omega^{-3/2}$. The conclusion now follows easily. For $0 < x \leqslant \omega^{-1/2}$, using the previous upper bounds and the explicit expression of $Q_\omega '$ we see that $| \phi_\omega '(x)| \leqslant C \varepsilon \omega^{3/2} x + | Q_\omega '(x)| \leqslant C \omega^{3/2} x$ and thus
\[ \left | \phi_\omega '(x) \int_x^{1/\sqrt{\omega}} \frac{\text{d}y}{\phi_\omega '(y)^2} \right | \leqslant C \omega^{3/2} x \int_x^{1/\sqrt{\omega}} \frac{C \omega^{-3} \, \text{d}y}{y^2} \leqslant C \omega^{3/2} x \cdot \frac{C \omega^{-3}}{x} \leqslant C \omega^{-3/2} \leqslant C \omega^{-3/2} e^{\sqrt{\omega} x}. \]
On the other hand, for $x> \omega^{-1/2}$, we have
\[ \int_{1/\sqrt{\omega}}^x \frac{\text{d}y}{\phi_\omega '(y)^2} \leqslant \int_{1/\sqrt{\omega}}^x \frac{C \omega^{-2} \, \text{d}y}{e^{-2 \sqrt{\omega} y}} \leqslant C \omega^{-5/2} e^{2 \sqrt{\omega} x}. \]
Using Lemma 2, we obtain
\[ \left | \phi_\omega '(x) \int_{1/ \sqrt{\omega}}^x \frac{\text{d}y}{\phi_\omega '(y)^2} \right | \leqslant C \omega e^{-\sqrt{\omega} x} \cdot C \omega^{-5/2} e^{2 \sqrt{\omega} x} \leqslant C \omega^{-3/2} e^{\sqrt{\omega} x}. \]
Hence the bound for $A_\omega$ is proved for all $x>0$. The bounds for its derivatives are similar and do not show additional difficulties, now that $\alpha_\omega$ is estimated. \hfill \qedsymbol

\noindent \textcolor{white}{a} \\ \noindent For any bounded continuous function $W$, define
\[ I_+ [W] (x) := \left | \begin{array}{ll} \displaystyle{- \phi_\omega '(x) \int_0^x A_\omega W - A_\omega (x) \int_x^{+ \infty} \phi_\omega ' W} & \text{if $x \geqslant 0$} \\
\\ \displaystyle{\phi_\omega '(x) \int_x^0 A_\omega W + A_\omega (x) \int_{- \infty}^x \phi_\omega ' W} & \text{if $x<0$.} \end{array} \right. \]
Note that if $\langle W \, , \phi_\omega ' \rangle = 0$ then $- \int_x^{+ \infty} \phi_\omega ' W = \int_{- \infty}^x \phi_\omega ' W$ and therefore the two expressions above coincide at $x=0$ and provide a solution to the equation $L_+U = W$. We will now provide estimates on $\Lambda_\omega$. In what follows, let us denote $\Lambda_\omega^Q := \omega \frac{\partial Q_\omega}{\partial \omega}$. First, we shall prove the following result, only here to be used in the next proof.

\begin{leftbar}
\noindent \textbf{Lemma 4.} For $\omega > 0$ small enough (as in the previous lemmas), $\Lambda_\omega$ is bounded on $\R$. 
\end{leftbar}

\noindent \textit{Proof.} Our proof relies on spectral arguments. To this end, let $L_+^Q := - \partial_x^2 + \omega - 3 Q_\omega^2$ and $L_+^{Q0} := - \partial_x^2 + 1 - 3 Q^2$. We know from \cite{Ch} that $L_+^{Q0}$ has only one negative eigenvalue which is $-3$, associated to the eigenfunction $Q^2$. The kernel of $L_+^{Q0}$ is generated by $Q'$. We know the following spectral coercivity property from \cite{Tao}: for any $u \in H^1 ( \R )$,
\[ \langle L_+^{Q0} u \, , u \rangle \geqslant c_1 ||u||_{H^1}^2 - c_2 | \langle u \, , Q^2 \rangle |^2 - c_3 | \langle u \, , Q' \rangle |^2  \]
with $c_1,c_2,c_3$ positive constants. Let $\text{Ev}_\omega u (x) = u \left ( \frac{x}{\sqrt{\omega}} \right )$. We see that $\text{Ev}_\omega^{-1} u(x) = u ( \sqrt{\omega} \, x)$, $\text{Ev}_\omega^\star = \sqrt{\omega} \, \text{Ev}_\omega^{-1}$ and $L_+^Q = \omega \, \text{Ev}_\omega^{-1} L_+^{Q0} \text{Ev}_\omega$. Using these identities, we compute
\[ \langle L_+^Q u \, , u \rangle = \sqrt{\omega} \langle L_+^{Q0} ( \text{Ev}_\omega u ) \, , ( \text{Ev}_\omega u ) \rangle \geqslant \sqrt{\omega} \left [ c_1 || \text{Ev}_\omega u ||_{H^1}^2 - c_2 | \langle \text{Ev}_\omega u \, , Q^2 \rangle |^2 - c_3 | \langle \text{Ev}_\omega u \, , Q' \rangle |^2 \right ] \]
where $\langle \text{Ev}_\omega u \, , Q^2 \rangle = \omega^{-1/2} \langle u \, , Q_\omega^2 \rangle$, $\langle \text{Ev}_\omega u \, , Q' \rangle = \omega^{-1/2} \langle u \, , Q_\omega ' \rangle$ and $|| \text{Ev}_\omega u ||_{H^1}^2 = \omega^{-1/2} || u ||_{H_\omega^1}^2$ with $||u||_{H_\omega^1}^2 := \omega ||u||^2 + ||u'||^2$. Hence, the following lower bound holds for all $u \in H^1 ( \R )$,
\[ \langle L_+^Q u \, , u \rangle \geqslant c_1 ||u||_{H_\omega^1}^2 - \frac{c_2}{\sqrt{\omega}} | \langle u \, , Q_\omega^2 \rangle |^2 - \frac{c_3}{\sqrt{\omega}} | \langle u \, , Q_\omega ' \rangle |^2. \]
Now, take $\varepsilon > 0$ which we will fix later. We take $\omega_0> 0$ small enough (to be fixed later) and $\omega > 0$ close enough to $\omega_0$ (we ask that $| \omega - \omega_0 | \leqslant \varepsilon \omega_0$). We denote $\tau := \frac{\phi_\omega - \phi_{\omega_0}}{\omega - \omega_0}$ that satisfies the equation
\[ \begin{array}{rl} & \displaystyle{\tau '' = \phi_\omega + \omega_0 \tau - ( \phi_\omega^2 + \phi_\omega \phi_{\omega_0}^2 + \phi_{\omega_0}^2 ) \tau + \phi_\omega \, \frac{g(\phi_\omega^2) - g(\phi_{\omega_0}^2)}{\omega - \omega_0} + g(\phi_{\omega_0}^2) \tau} \\
\\ \text{i.e.} & \displaystyle{L_+^{Q} \tau = - \phi_\omega + (\phi_\omega^2 + \phi_\omega \phi_{\omega_0} + \phi_{\omega_0}^2 - 3 Q_{\omega_0}^2 ) \tau - \phi_\omega \, \frac{g(\phi_\omega^2) - g(\phi_{\omega_0}^2)}{\omega - \omega_0} - g(\phi_{\omega_0}^2) \tau,} \end{array} \]
where $L_+^Q$ is the previous operator with the pulsation $\omega_0$. We take $\omega_0$ small enough such that the bounds in Lemma 2 hold. Moreover, we see that $| Q_\omega - Q_{\omega_0} | \leqslant C | \omega - \omega_0 | \omega_0^{-1/2} \leqslant C \varepsilon \sqrt{\omega_0}$. To see that, recall that $Q_\omega$ is known explicitly, thus we can compute $\Lambda_\omega^Q := \sqrt{\frac{\omega}{2}} \left ( 1 - \sqrt{\omega} x \tanh ( \sqrt{\omega} x) \right ) \, \frac{1}{\cosh ( \sqrt{\omega} x )}$ which gives $| \Lambda_\omega^Q | \leqslant C \sqrt{\omega_0}$ and then $| \partial_\omega Q_\omega | \leqslant C \omega_0^{-1/2}$. This proves the upper bound on $| Q_\omega - Q_{\omega_0} |$. Now, let us estimate $\langle L_+^Q \tau \, , \tau \rangle$. First,
\[ | \langle \phi_\omega \, , \tau \rangle | \leqslant || \phi_\omega || \, || \tau || \leqslant C \omega_0^{1/4} || \tau ||. \]
Now, about the second term, writing
\[ \begin{array}{rcl} | \phi_\omega^2 + \phi_\omega \phi_{\omega_0} + \phi_{\omega_0}^2 - 3 Q_{\omega_0}^2 | & \leqslant & (\phi_\omega + Q_\omega ) | \phi_\omega - Q_\omega | + (Q_\omega + Q_{\omega_0} ) | Q_\omega - Q_{\omega_0} | \\ \\ & & \, \, \, + \, \phi_\omega | \phi_{\omega_0} - Q_{\omega_0} | + Q_{\omega_0} | \phi_\omega - Q_\omega | + Q_{\omega_0} | Q_\omega - Q_{\omega_0} | \\
\\ & & \, \, \, + \, (\phi_{\omega_0} + Q_{\omega_0} ) | \phi_{\omega_0} - Q_{\omega_0} |, \end{array}  \]
we get $| \phi_\omega^2 + \phi_\omega \phi_{\omega_0} + \phi_{\omega_0}^2 - 3Q_{\omega_0}^2 | \leqslant C \varepsilon \omega_0$. Thus, 
\[ \left | \left \langle ( \phi_\omega^2 + \phi_\omega \phi_{\omega_0} + \phi_{\omega_0}^2 ) \tau \, , \tau \right \rangle \right | \leqslant C \varepsilon \omega_0 || \tau ||^2. \]
Now, about the third term, we take $\omega_0$ (and $\omega$) small enough such that $| g'(s) | \leqslant \varepsilon$ for all $s \in [0 \, , 5 \omega_0 ]$. This implies $|g(\phi_\omega^2) - g(\phi_{\omega_0}^2)| \leqslant \varepsilon | \phi_\omega^2 - \phi_{\omega_0}^2 |$, which leads to
\[ \left | \phi_\omega \, \frac{g(\phi_\omega^2) - g(\phi_{\omega_0}^2)}{\omega - \omega_0} \right | \leqslant \phi_\omega \varepsilon | \tau | ( \phi_\omega + \phi_{\omega_0} ) \leqslant C \varepsilon \omega_0 | \tau |. \]
Thus, 
\[ \left | \left \langle \phi_\omega \, \frac{g(\phi_\omega^2) - g(\phi_{\omega_0}^2)}{\omega - \omega_0} \, \tau \, , \tau \right \rangle \right | \leqslant C \varepsilon \omega_0 || \tau ||^2. \]
Finally, about the last term, $| g(\phi_{\omega_0}^2) | \leqslant \varepsilon \phi_{\omega_0}^2 \leqslant C \varepsilon \omega_0$, thus $\left | \langle g(\phi_{\omega_0}^2) \tau \, , \tau \rangle \right | \leqslant C \varepsilon \omega_0 || \tau ||^2$. Gathering these estimates, we have
\[ | \langle L_+^Q \tau \, , \tau \rangle | \leqslant C \omega_0^{1/4} || \tau || + C \varepsilon \omega_0 || \tau ||^2. \]
Using the spectral inequality, and since $\tau \in H^1 ( \R )$, we know that $\langle L_+^Q \tau \, , \tau \rangle \geqslant c_1 || \tau ||_{H_{\omega_0}^1}^2 - \frac{c_2}{\sqrt{\omega_0}} | \langle \tau \, , Q_{\omega_0}^2 \rangle |^2 - \frac{c_3}{\sqrt{\omega_0}} | \langle \tau \, , Q_{\omega_0}' \rangle |^2$. Since $\tau$ is even and $Q_{\omega_0} '$ is odd, $\langle \tau \, , Q_{\omega_0} ' \rangle = 0$. We estimate the other scalar product as follows, using both that $L_+^Q Q_{\omega_0}^2 = -3 \omega_0 Q_{\omega_0}^2$ and that $L_+^Q$ is self-adjoint:
\[ \begin{array}{rcl} | \langle \tau \, , Q_{\omega_0}^2 \rangle | &=& \displaystyle{\frac{1}{3 \omega_0} | \langle \tau \, , L_+^Q Q_{\omega_0}^2 \rangle | \, \, = \, \, \frac{1}{3 \omega_0} | \langle L_+^Q \tau \, , Q_{\omega_0}^2 \rangle |} \\
\\ & \leqslant & \displaystyle{\frac{1}{3 \omega_0} \left [ | \langle \phi_\omega \, , Q_{\omega_0}^2 \rangle | + | \langle (\phi_\omega^2 + \phi_\omega \phi_{\omega_0} + \phi_{\omega_0}^2 -3 Q_{\omega_0}^2 ) \tau \, , Q_{\omega_0}^2 \rangle | + \left | \left \langle \phi_\omega \, \frac{g(\phi_\omega^2) - g(\phi_{\omega_0}^2)}{\omega - \omega_0} \, , Q_{\omega_0}^2 \right \rangle \right | \right.} \\ \\ & & \displaystyle{\left. \, \, \, \, \, \, \, \, \, \, \, \, \, \, \, \, \, \, \textcolor{white}{\frac{1}{2}} \, \, \, \, + \, | \langle g(\phi_{\omega_0}^2) \tau \, , Q_{\omega_0}^2 \rangle | \right ].} \end{array} \]
Directly using the exponential controls, we find $| \langle \phi_\omega \, , Q_{\omega_0}^2 \rangle | \leqslant C \omega_0$. In order to control the other terms, we recall the estimates proved above:
\[ | \phi_\omega^2 + \phi_\omega \phi_{\omega_0} + \phi_{\omega_0}^2 - 3 Q_{\omega_0}^2 | \leqslant C \varepsilon \omega_0 , \, \, \, \, \left | \phi_\omega \, \frac{g(\phi_\omega^2) - g(\phi_{\omega_0}^2)}{\omega - \omega_0} \right | \leqslant C \varepsilon \omega_0 | \tau | \, \, \, \, \text{and} \, \, \, \, | g(\phi_{\omega_0}^2) | \leqslant C \varepsilon \omega_0. \]
This leads to: first,
\[ | \langle (\phi_\omega^2 + \phi_\omega \phi_{\omega_0} + \phi_{\omega_0}^2 -3 Q_{\omega_0}^2 ) \tau \, , Q_{\omega_0}^2 \rangle | \leqslant C \varepsilon \omega_0 \int_{\R} | \tau | Q_{\omega_0}^2 \leqslant C \varepsilon \omega_0 || \tau || \, || Q_{\omega_0}^2 || \leqslant C \varepsilon \omega_0^{7/4} || \tau ||; \]
second,
\[ \left | \left \langle \phi_\omega \, \frac{g(\phi_\omega^2) - g(\phi_{\omega_0}^2)}{\omega - \omega_0} \, , Q_{\omega_0}^2 \right \rangle \right | \leqslant C \varepsilon \omega_0 \int_{\R} |\tau| Q_{\omega_0}^2 \leqslant C \varepsilon \omega_0^{7/4} || \tau ||; \]
and third,
\[ | \langle g(\phi_{\omega_0}^2) \tau \, , Q_{\omega_0}^2 \rangle | \leqslant C \varepsilon \omega_0 \int_{\R} | \tau | Q_{\omega_0}^2 \leqslant C \varepsilon \omega_0^{7/4} || \tau ||. \]
Overall, we obtain $\langle \tau \, , Q_{\omega_0}^2 \rangle | \leqslant C + C \varepsilon \omega_0^{3/4} || \tau ||$ which leads to $| \langle \tau \, , Q_{\omega_0}^2 \rangle |^2 \leqslant C + C \varepsilon \omega_0^{3/4} || \tau || + C \varepsilon^2 \omega_0^{3/2} || \tau ||^2$. Henceforth, going back to the spectral inequality, we obtain
\[ \begin{array}{rl} & \displaystyle{|| \tau ||_{H_{\omega_0}^1}^2 \leqslant C | \langle L_+^Q \tau \, , \tau \rangle | + \frac{C}{\sqrt{\omega_0}} | \langle \tau \, , Q_{\omega_0}^2 \rangle |^2 \leqslant C \omega_0^{1/4} || \tau || + C \varepsilon \omega_0 || \tau||^2 + C \omega_0^{-1/2} + C \varepsilon \omega_0^{1/4} || \tau || + C \varepsilon^2 \omega_0 || \tau ||^2} \\
\\ \text{thus} & \displaystyle{\omega_0 || \tau ||^2 \leqslant || \tau ||_{H_{\omega_0}^1}^2 \leqslant C \omega_0^{1/4} || \tau || + C \varepsilon \omega_0 || \tau ||^2 + C \omega_0^{-1/2}} \\
\\ \text{thus} & \displaystyle{\omega_0 ( 1 - C \varepsilon ) || \tau ||^2 - C \omega_0^{1/4} || \tau || - C \omega_0^{-1/2} \leqslant 0.} \end{array} \]
Choosing $\varepsilon > 0$ small enough, we may assume $1 - C \varepsilon \geqslant \frac{1}{2}$ and thus
\[ \frac{\omega_0}{2} || \tau ||^2 - C \omega_0^{1/4} || \tau || - C \omega_0^{-1/2} \leqslant 0. \]
The positive root of the polynomial $\frac{\omega_0}{2} X^2 - C \omega_0^{1/4} X - C \omega_0^{-1/2}$ being $C \omega_0^{-3/4}$ (where the constant $C$ is different), we have $|| \tau || \leqslant C \omega_0^{-3/4}$. \\
\\ Now, recalling that $|| \tau ' ||^2 \leqslant || \tau ||_{H_{\omega_0}^1}^2 \leqslant C \omega_0^{1/4} || \tau ||  + C \varepsilon \omega_0 || \tau ||^2 + C \omega_0^{-1/2}$ and using the upper bound above about $|| \tau ||$, we get $|| \tau ' ||^2 \leqslant C \omega_0^{-1/2}$. This leads to $|| \tau ||_{L^\infty}^2 \leqslant 2 || \tau || \, || \tau' || \leqslant C \omega_0^{-3/4} \omega_0^{-1/4} = C \omega_0^{-1}$ and thus $|| \tau ||_{L^\infty} \leqslant C \omega_0^{-1/2}$. \\
\\ Now, take $x \in \R$ fixed. We have $\left | \frac{\phi_\omega (x) - \phi_{\omega_0} (x)}{\omega - \omega_0} \right | = | \tau (x) | \leqslant C \omega_0^{-1/2}$ for $\omega$ taken as before. Letting $\omega \to \omega_0$, we obtain $\left | ( \partial_\omega \phi_\omega )_{\omega = \omega_0} (x) \right | \leqslant C \omega_0^{-1/2}$ and thus $| \Lambda_{\omega_0} (x) | \leqslant C \sqrt{\omega_0}$. As we will see in the next lemma, we could not hope for a better estimate. The constant $C$ is uniform (it does not depend on $x$), showing that $\Lambda_{\omega_0}$ is indeed bounded. This is the result announced. \hfill \qedsymbol

\textcolor{white}{a} \\ \noindent Now let us give more precise bounds about $\Lambda_\omega$.

\begin{leftbar}
\noindent \textbf{Lemma 5.} Assume $g$ to be $\mathscr{C}^{5} ( (0 \, , + \infty ))$, $\mathscr{C}^1 ([0 \, , \infty ))$ and such that $g(0)=g'(0)=0$. For any $k \in [\![ 0 \, , 6 ]\!]$, there exists $C_k > 0$ such that, for any $\omega >0$ small enough and any $x \in \R$,
\[ | \Lambda_\omega^{(k)} (x) | \leqslant C_k \omega^{\frac{1+k}{2}} (1 + \sqrt{\omega} |x| ) e^{- \sqrt{\omega} |x|}. \]
Moreover, for every $\varepsilon > 0$, for any $\omega >0$ small enough, 
\[ | \Lambda_\omega (x) - \Lambda_\omega^Q (x) | \leqslant \varepsilon \sqrt{\omega} (1 + \sqrt{\omega} |x| ) e^{- \sqrt{\omega} |x|}.  \]
At last, for $\omega$ small enough, $\langle \phi_\omega \, , \Lambda_\omega \rangle \geqslant C \sqrt{\omega}$. 
\end{leftbar}

\noindent \textit{Proof.} The condition $\langle W \, , \phi_\omega ' \rangle = 0$ is in particular satisfied by $W= - \omega \phi_\omega$ since $\phi_\omega \phi_\omega '$ is odd. We know that $L_+ \Lambda_\omega = - \omega \phi_\omega$. Hence, there exists some constants $c_\omega^A,c_\omega^\phi$ (possibly depending on $\omega$) such that $\Lambda_\omega = I_+ [ - \omega \phi_\omega ] + c_\omega^A A_\omega + c_\omega^\phi \phi_\omega '$. Since $I_+ [ - \omega \phi_\omega ]$, $A_\omega$ and $\Lambda_\omega$ are even while $\phi_\omega '$ is odd, we obtain $c_\omega^\phi = 0$. Moreover, since $\Lambda_\omega$ is bounded on $\R$ (see Lemma 4), $c_\omega^A = 0$. Hence $\Lambda_\omega = I_+ [ - \omega \phi_\omega ]$. We also easily check that, using the bounds on $\phi_\omega$, $\phi_\omega '$ and $A_\omega$, we have $| I_+ [ - \omega \phi_\omega ] (x) | \leqslant C \sqrt{\omega} ( 1 + \sqrt{\omega} |x| ) e^{- \sqrt{\omega} |x|}$. The term $\omega |x| e^{- \sqrt{\omega} |x|}$ comes from the first integral in the definition of $I_+$. Thus,
\[ | \Lambda_\omega (x) | \leqslant C \sqrt{\omega} ( 1 + \sqrt{\omega} |x| ) e^{- \sqrt{\omega} |x|}. \]
Differentiating the formula $\Lambda_\omega = I_+ [ - \omega \phi_\omega ]$, we similarly get the estimates on the derivatives of $\Lambda_\omega$. Now consider the second point of the lemma: let $\varepsilon > 0$ and $\delta > 0$ which will be fixed later (depending on $\varepsilon$). The proof is similar to the one of the analogous result in Lemma 2. Let us denote $\Theta_\omega := \Lambda_\omega - \Lambda_\omega^Q$. Recalling that $P_\omega = \phi_\omega - Q_\omega$, the equation satisfied by $\Theta_\omega$ is
\[ \Theta_\omega '' = \omega P_\omega + \omega \Theta_\omega - 3 \phi_\omega^2 \Theta_\omega -3 \Lambda_\omega^Q P_\omega ( \phi_\omega + Q_\omega ) + 2 \Lambda_\omega \phi_\omega^2 g'(\phi_\omega^2) + \Lambda_\omega \phi_\omega^2. \]
Taking $\omega$ small enough, we can assume that $|P_\omega| \leqslant \delta \sqrt{\omega}$, $\phi_\omega^2 \leqslant \zeta_\omega^2 \leqslant 3 \omega$, $|g'(\phi_\omega^2)| \leqslant \delta$ and $|g(\phi_\omega^2)| \leqslant \delta \phi_\omega^2 \leqslant C \delta \omega$. We also see, from the bound above about $\Lambda_\omega$, that $| \Lambda_\omega | \leqslant C \sqrt{\omega}$ (for example, observe that $x \mapsto (1 + \sqrt{\omega} x ) e^{- \sqrt{\omega} x}$ is nonincreasing on $[0 \, , + \infty )$). Gathering these bounds we obtain
\[ | \Theta_\omega '' | \leqslant C \delta \omega^{3/2} + 10 \omega | \Theta_\omega |. \]
We can assume $\omega$ small enough such that $\left | \omega \, \frac{\text{d} \zeta_\omega}{\text{d} \omega} - \sqrt{\frac{\omega}{2}} \right | \leqslant \delta \sqrt{\omega}$ i.e. $| \Theta_\omega (0)| \leqslant \delta \sqrt{\omega}$. By Grönwall's lemma, we get that, for any $x>0$,
\[ | \Theta_\omega (x) | \leqslant \sqrt{\omega} \left [ \frac{C \delta}{10} + e^{10 \sqrt{\omega} x} \left ( \delta + \frac{C \delta}{10} \right ) \right ] \leqslant C \delta \sqrt{\omega} (1 + e^{10 \sqrt{\omega} x}). \]
We also know that $| \Theta_\omega (x)| \leqslant C \sqrt{\omega} ( 1 + \sqrt{\omega} x ) e^{- \sqrt{\omega} x} \leqslant C \sqrt{\omega} \, e^{- \sqrt{\omega} x /2}$. Denoting $x_\omega := 2 \omega^{-1/2} \ln (C/\varepsilon)$, we see that, for any $x \geqslant x_\omega$, $| \Theta_\omega (x) | \leqslant C \sqrt{\omega} \, e^{- \sqrt{\omega} x_\omega / 2} = \varepsilon \sqrt{\omega}$. On the other hand, for any $x \in [0 \, , x_\omega ]$, $| \Theta_\omega (x) | \leqslant C \delta \sqrt{\omega} \left ( 1 + C \varepsilon^{-20} \right ) \leqslant \varepsilon \sqrt{\omega}$, provided we take $\delta$ small enough (depending on $\varepsilon$ only, not depending on $\omega$). Therefore, we have proved that $|| \Theta_\omega ||_\infty \leqslant \varepsilon \sqrt{\omega}$. \\
\\ Now, consider $\widetilde{T}_\omega := -3 \phi_\omega^2 \Theta_\omega -3 \Lambda_\omega^Q P_\omega ( \phi_\omega + Q_\omega ) + 2 \Lambda_\omega \phi_\omega^2 g'(\phi_\omega^2) + \Lambda_\omega g(\phi_\omega^2)$ and $T_\omega := \omega P_\omega + \widetilde{T}_\omega$, in order that $\Theta_\omega '' - \omega \Theta_\omega = T_\omega$. The method of the variation of the constants and the initial condition $\Theta_\omega '(0)=0$ show that, for $x > 0$,
\[ \Theta_\omega (x) = \left ( \frac{\Theta_\omega (0)}{2} + \frac{\text{IT}^-}{2 \sqrt{\omega}} \right ) e^{\sqrt{\omega} x} + \frac{\Theta_\omega (0)}{2} e^{- \sqrt{\omega} x} - \frac{e^{\sqrt{\omega} x}}{2 \sqrt{\omega}} \int_x^{+ \infty} T_\omega (y) e^{- \sqrt{\omega} y} \, \text{d}y - \frac{e^{- \sqrt{\omega} x}}{2 \sqrt{\omega}} \int_0^x (\omega P_\omega (y) + \widetilde{T}_\omega (y)) e^{\sqrt{\omega} y} \, \text{d}y, \]
where $\text{IT}^- = \int_0^{+ \infty} T_\omega (y) e^{- \sqrt{\omega} y} \, \text{d}y$. The previous bounds on $\phi_\omega$ and $\Lambda_\omega$ assure the existence of $\text{IT}^-$ and of all the integral terms in the expression of $\Theta_\omega (x)$. Since $\Theta_\omega (x) \, \underset{x \to + \infty}{\longrightarrow} \, 0$, $\frac{\Theta_\omega (0)}{2} + \frac{\text{IT}^-}{2 \sqrt{\omega}} = 0$. Moreover, using the bounds on $\phi_\omega$ and $\Lambda_\omega$, we see that
\[ \left | \int_x^{+ \infty} T_\omega (y) e^{- \sqrt{\omega} y} \, \text{d}y \right | \leqslant \varepsilon \omega e^{-2 \sqrt{\omega} x}, \, \, \, \left | \int_0^x \omega P_\omega (y) e^{\sqrt{\omega} y} \, \text{d}y \right | \leqslant \varepsilon \omega^{3/2} x , \, \, \, \left | \int_0^x \widetilde{T}_\omega (y) e^{\sqrt{\omega} y} \, \text{d}y \right | \leqslant \varepsilon \omega. \]
Gathering these estimates in the expression of $\Theta_\omega$, we obtain
\[ | \Theta_\omega (x) | \leqslant C \varepsilon \sqrt{\omega} ( 1 + \sqrt{\omega} x) e^{- \sqrt{\omega} x}, \]
which is the desired result. \\
\\ For the last point of the lemma, we take $\varepsilon > 0$ that we will fix later. Providing we take $\omega$ small enough, we have
\[ | \phi_\omega (x) - Q_\omega (x) | \leqslant \varepsilon \sqrt{\omega} \, e^{- \sqrt{\omega} |x|} \, \, \, \, \, \text{and} \, \, \, \, \, | \Lambda_\omega (x) - \Lambda_\omega^Q (x) | \leqslant \varepsilon \sqrt{\omega} ( 1 + \sqrt{\omega} |x| ) e^{- \sqrt{\omega} |x|} \]
where we recall that
\[ Q_\omega (x) = \frac{\sqrt{2 \omega}}{\cosh ( \sqrt{\omega} \, x )} \, \, \, \, \, \text{and} \, \, \, \, \, \Lambda_\omega^Q (x) = \sqrt{\frac{\omega}{2}} \left ( 1 - \sqrt{\omega} \, x \tanh ( \sqrt{\omega} \, x) \right ) \, \frac{1}{\cosh ( \sqrt{\omega} \, x )}. \]
We write that $\langle \phi_\omega \, , \Lambda_\omega \rangle = \langle Q_\omega \, , \Lambda_\omega^Q \rangle + \langle \phi_\omega - Q_\omega \, , \Lambda_\omega^Q \rangle + \langle \phi_\omega \, , \Lambda_\omega - \Lambda_\omega^Q \rangle$, where
\[ \langle Q_\omega \, , \Lambda_\omega^Q \rangle = 2 \sqrt{\omega} \int_0^{+ \infty} (1 - y \tanh y ) \frac{\text{d}y}{\cosh^2 (y)} = \sqrt{\omega} \int_0^{+ \infty} \frac{\text{d}y}{\cosh^2 (y)} \geqslant \frac{\sqrt{\omega}}{2}, \]
integrating by parts. Using the control on $\phi_\omega - Q_\omega$ we find
\[ | \langle \phi_\omega - Q_\omega \, , \Lambda_\omega^Q \rangle | \leqslant 2 \varepsilon \sqrt{2 \omega} \int_0^{+ \infty} e^{-2y} (1+y) \, \text{d}y = C \varepsilon \sqrt{\omega}. \]
Using the control on $\Lambda_\omega - \Lambda_\omega^Q$ we similarly find that $| \langle \phi_\omega \, , \Lambda_\omega - \Lambda_\omega^Q \rangle | \leqslant C \varepsilon \sqrt{\omega}$. Gathering these estimates we find $\langle \phi_\omega \, , \Lambda_\omega \rangle \geqslant \left ( \frac{1}{2} - C \varepsilon \right ) \sqrt{\omega} \geqslant \frac{\sqrt{\omega}}{4}$ provided we take $\varepsilon$ small enough (and thus $\omega$ small enough). \hfill \qedsymbol

\subsection{Conjugate identity}
\noindent Let $S = \phi_\omega \cdot \partial_x \cdot \frac{1}{\phi_\omega}$ so that $S^* = - \frac{1}{\phi_\omega} \cdot \partial_x \cdot \phi_\omega$. Let us define
\[ \begin{array}{rl} & M_+ = - \partial_x^2 + \omega - g ( \phi_\omega^2 ) + 2 \, \frac{G ( \phi_\omega^2 )}{\phi_\omega^2} \\ \\ \text{and}& M_- = - \partial_x^2 + \omega - 5 g ( \phi_\omega^2) + 2 \phi_\omega^2 \, g'( \phi_\omega^2) + 6 \, \frac{G ( \phi_\omega^2)}{\phi_\omega^2}. \end{array} \]

\begin{leftbar}
\noindent \textbf{Lemma 6.} We have $S^2 L_+ L_- = M_+ M_- S^2$.
\end{leftbar}

\noindent \textit{Proof.} From (3.25)-(3.26) of \cite{Ch} we recall the following general formula: for any nonvanishing function $R$, denoting $V_\pm = R^2 \pm 3R' + \frac{R''}{R}$, we have
\begin{equation}
    ( \partial_x - R ) ( \partial_x^2 - V_+ ) ( \partial_x + R) = ( \partial_x + R) ( \partial_x^2 - V_-) ( \partial_x - R).
    \label{eqR}
\end{equation}
Let us apply this identity with $R = \phi_\omega ' / \phi_\omega$. Thanks to \eqref{eqphi} and the identity $(\phi_\omega ')^2 = \omega \phi_\omega^2 - \frac{1}{2} \phi_\omega^4 + G ( \phi_\omega^2 )$ that is itself derived from \eqref{eqphi}, we find that
\[ \begin{array}{rl} & \displaystyle{R^2 = \omega - \frac{1}{2} \phi_\omega^2 + \frac{G ( \phi_\omega^2 )}{\phi_\omega^2},} \\
\\ & \displaystyle{R' = - \frac{1}{2} \phi_\omega^2 + g ( \phi_\omega^2 ) - \frac{G ( \phi_\omega^2 )}{\phi_\omega^2},} \\
\\ \text{and} & \displaystyle{\frac{R''}{R} = - \phi_\omega^2 + 2 \left ( \phi_\omega^2 \, g'( \phi_\omega^2) - g ( \phi_\omega^2 ) + \frac{G ( \phi_\omega^2 )}{\phi_\omega^2} \right ).} \end{array} \]
The last expression shows that, even though $R$ vanishes at $x=0$, $R''/R$ can be extended by continuity without any complication. Hence, \eqref{eqR} remains valid and we get
\[ V_+ = \omega - 3 \phi_\omega^2 + 2 \phi_\omega^2 \, g'( \phi_\omega^2) + g ( \phi_\omega^2 ) \, \, \, \, \, \text{and} \, \, \, \, \, V_- = \omega - 5 g( \phi_\omega^2) + 2 \phi_\omega^2 g'( \phi_\omega^2) + 6 \, \frac{G ( \phi_\omega^2 )}{\phi_\omega^2}. \]
We easily check that $\partial_x - R= S$, $\partial_x + R = S^*$, $\partial_x^2 - V_+ = -L_+$ and $\partial_x^2 - V_- = -M_-$. We also check that $S^*S = L_-$ and $SS^* = M_+$. Thus the identity we have started with gives $-SL_+S^* = -S^* M_- S$. Composing by $S$ on the left and $S$ on the right, we get $S^2 L_+L_- = M_+M_-S^2$. \hfill \qedsymbol
\noindent \\ \textcolor{white}{a} \\ In what follows, we will denote $a_\omega^- = -5 g(\phi_\omega^2) + 2 \phi_\omega^2 g'(\phi_\omega^2) + 6 \frac{G(\phi_\omega^2)}{\phi_\omega^2}$ and $a_\omega^+ = -g(\phi_\omega^2) + 2 \frac{G(\phi_\omega^2)}{\phi_\omega^2}$ (in order that $M_\pm = - \partial_x^2 + \omega + a_\omega^{\pm}$). These potentials are crucial in our proof. 

\subsection{Invertibility of \texorpdfstring{$M_-$}{M-}}
\noindent In this section we assume that $\text{Ker} (M_-) = \{ 0 \}$. In the next section, Corollary 1 will show that hypotheses $(H_1)$ and $(H_2)$ are sufficient to ensure that this assumption is true. We follow the same reasoning as \cite{Ma1}. Denoting by $B_1$ and $B_2$ two solutions of $M_- B_1 = M_- B_2 = 0$ satisfying
\[ |B_1^{(k)} (x)| \leqslant C_k \omega^{- \frac{1}{4} + \frac{k}{2}} e^{- \sqrt{\omega} x} , \, \, \, \, \, \, | B_2^{(k)} (x) | \leqslant C_k \omega^{- \frac{1}{4} + \frac{k}{2}} e^{\sqrt{\omega} x} \]
for $C_k > 0$ and $B_1 B_2 ' - B_1 ' B_2 = 1$ on $\R$. These estimates are proved as in Lemma 3. Two such independent solutions exist because $\text{Ker} (M_-) = 0$. For any bounded continuous function $W$, the formula
\[ J_- [W] (x) := B_1 (x) \int_{- \infty}^x B_2 W + B_2 (x) \int_x^{+ \infty} B_1 W \]
defines a solution to $M_- U = W$.

\section{Non-existence of internal modes}
\noindent As explained in the introduction, we seek hypotheses on $g$ that will ensure that the equation \eqref{NLS} does not have internal modes. An internal mode is a solution $(X \, , Y \, , \lambda ) \in H^1 ( \R )^2 \times \C$ to the following system:
\[ \left \{ \begin{array}{ccl} L_- X &=& \lambda Y \\ L_+ Y &=& \lambda X. \end{array} \right. \]
For $\omega$ small enough, let us denote $P_B^{\pm} = - (a_\omega^\pm)' \frac{\Phi_B}{\zeta_B^2}$ and $P_B = \frac{P_B^+ + P_B^-}{2}$. We recall the definition of $\varepsilon_\omega := \sup\limits_{0 \leqslant s \leqslant 3 \omega} |sg''(s)|$. We recall that $\omega$ is always assumed small enough so that $\phi_\omega \leqslant \zeta_\omega \leqslant \sqrt{3 \omega}$. Under the hypothesis $(H_1)$, Taylor's formula gives 
\[ \left | \frac{G ( \phi_\omega^2 )}{\phi_\omega^4} \right | \leqslant \varepsilon_\omega , \, \, \, \, \, \, \, \left | \frac{g( \phi_\omega^2 )}{\phi_\omega^2} \right | \leqslant \varepsilon_\omega , \, \, \, \, \, \, \, \left | g'( \phi_\omega^2 ) \right | \leqslant \varepsilon_\omega , \, \, \, \, \, \, \, \left | \phi_\omega^2 g '' ( \phi_\omega^2 ) \right | \leqslant \varepsilon_\omega. \]
Therefore, using the expressions of $P_B^+$ and $P_B^-$, and also using that $| \phi_\omega ' / \phi_\omega | \leqslant C \sqrt{\omega}$ we see that
\[ |P_B (x)| \leqslant C \varepsilon_\omega \phi_\omega^2 (x) \left | \frac{\phi_\omega '(x)}{\phi_\omega (x)} \right | \frac{|\Phi_B (x)|}{\zeta_B^2 (x)} \leqslant C \sqrt{\omega} \varepsilon_\omega |x| \phi_\omega^2 (x) \leqslant C \sqrt{\omega} \varepsilon_\omega x \omega e^{- \sqrt{\omega} |x|} \leqslant C \varepsilon_\omega \omega e^{- \sqrt{\omega} |x| /10}.  \]
From now on, in everything that follows, we assume the hypothesis $(H_1)$ to be satisfied. \\
\\ The following lemma is a coercivity result about the quadratic form $u \mapsto \int_{\R} P_B u^2$. It is a weaker version of a theorem from Simon, see \cite{Si}. The proof given here is elementary. This result will intervene both in the proof of the spectral question we study here, and in the proof of the main theorem that will take place later.

\begin{leftbar}
\noindent \textbf{Lemma 7.} Assume that $\displaystyle{\int_{\R} \frac{a_{\omega}^+ + a_{\omega}^-}{2}} > 0$. For $\omega > 0$ small enough and $B>0$ large enough, for any $u \in H^1 ( \R )$,
\[ \int_{\R} P_B u^2 \geqslant C \gamma_B \varepsilon_\omega \sqrt{\omega} \int_{\R} \rho u^2 - \frac{C \varepsilon_\omega \sqrt{\omega}}{\gamma_B} \int_{\R} (u')^2  \]
where $\displaystyle{P_B = \frac{P_B^+ + P_B^-}{2} = - \frac{(a_{\omega}^+ + a_{\omega}^-)'}{2} \, \frac{\Phi_B}{\zeta_B}}$ and $\displaystyle{\gamma_B := \int_{\R} \frac{P_B}{\varepsilon_\omega} \in \, ]0 \, , C \sqrt{\omega} [}$. \\
\\ Setting  $P_\infty := - \frac{x(a_\omega^+ + a_\omega^- )'}{2}$ and $\gamma_\infty := \varepsilon_\omega^{-1} \int_{\R} P_\infty$, the same result holds replacing $B$ by $\infty$ everywhere: for $\omega > 0$ small enough and any $u \in H^1 ( \R )$,
\[ \int_{\R} P_\infty u^2 \geqslant C \gamma_\infty \varepsilon_\omega \sqrt{\omega} \int_{\R} \rho u^2 - \frac{C \varepsilon_\omega \sqrt{\omega}}{\gamma_\infty} \int_{\R} (u')^2. \]
\end{leftbar}

\noindent \textit{Proof.} We start by writing that, for $x,y \in \R$, $\displaystyle{u^2 (x) = u^2(y) - 2 \int_x^y u'(z) u(z) \, \text{d}z}$. In what follows let us denote $\widetilde{P_B} (y) := \frac{P_B (y)}{C \omega \varepsilon_\omega}$ such that $|\widetilde{P_B} (y)| \leqslant e^{- \sqrt{\omega} |y| /10}$. We multiply the previous identity by $\widetilde{P_B} (y)$ and integrate in $y$, leading to
\[ \left ( \int_{\R} \widetilde{P_B} \right ) u^2 (x) = \int_{\R} u^2 \widetilde{P_B} -2 \int_x^{+ \infty} \widetilde{P_B} (y) \int_x^y u'(z) u(z) \, \text{d}z \, \text{d}y + 2 \int_{- \infty}^x \widetilde{P_B} (y) \int_y^x u'(z) u(z) \, \text{d}z. \]
We now multiply by $e^{- \sqrt{\omega} |x|/10}$ and integrate in $x$, using $\int_{\R} e^{- \sqrt{\omega} |x|/2} \, \text{d}x = \frac{C}{\sqrt{\omega}}$:
\[  \begin{array}{rcl} \displaystyle{\left ( \int_{\R} \widetilde{P_B} \right ) \int_{\R} u^2 (x) e^{- \sqrt{\omega} |x|/10} \, \text{d}x} &=& \displaystyle{\frac{C}{\sqrt{\omega}} \int_{\R} u^2 \widetilde{P_B} - 2 \int_{\R} e^{- \sqrt{\omega} |x|/10} \int_x^{+ \infty} \widetilde{P_B} (y) \int_x^y u'(z) u(z) \, \text{d}z \, \text{d}y \, \text{d}x} \\
\\ & & \, \, \, \, \, \, \displaystyle{+ 2 \int_{\R} e^{- \sqrt{\omega} |x|/10} \int_{- \infty}^x \widetilde{P_B} (y) \int_y^x u'(z) u(z) \, \text{d}z \, \text{d}y \, \text{d}x}. \end{array} \]
By the Fubini theorem,
\[ \int_{\R} e^{- \sqrt{\omega} |x| /10} \int_x^{+ \infty} \widetilde{P_B} (y) \int_x^y u'(z) u(z) \, \text{d}z \, \text{d}y \, \text{d}x = \int_{\R} \left ( \int_{- \infty}^z e^{- \sqrt{\omega} |x|/10} \, \text{d}x \right ) \left ( \int_z^{+ \infty} \widetilde{P_B} (y) \, \text{d}y \right ) u'(z) u(z) \, \text{d}z. \]
We notice that
\[ \int_{- \infty}^z e^{- \sqrt{\omega} |x|/10} \, \text{d}x \leqslant \frac{C}{\sqrt{\omega}} \, \, \, \, \text{if $z>0$} \, \, \, \, \, \, \, \, \, \, \, \text{and} \, \, \, \, \, \, \, \, \, \, \, \int_{- \infty}^z e^{- \sqrt{\omega} |x|/10} \, \text{d}x \leqslant \frac{C}{\sqrt{\omega}} e^{- \sqrt{\omega} |x|/10} \, \, \, \, \text{if $z<0$}. \]
Similarly, since $| \widetilde{P_B} (y)| \leqslant e^{- \sqrt{\omega} |y| /10}$,
\[ \left | \int_z^{+ \infty} \widetilde{P_B} (y) \, \text{d}y \right | \leqslant \frac{C}{\sqrt{\omega}} e^{- \sqrt{\omega} |x|/10} \, \, \, \, \text{if $z>0$} \, \, \, \, \, \, \, \, \, \, \, \text{and} \, \, \, \, \, \, \, \, \, \, \, \left | \int_z^{+ \infty} \widetilde{P_B} (y) \, \text{d}y \right | \leqslant \frac{C}{\sqrt{\omega}}  \, \, \, \, \text{if $z<0$}. \]
Thus, for all $z \in \R$,
\[ \left | \left ( \int_{- \infty}^z e^{- \sqrt{\omega} |x|/10} \, \text{d}x \right ) \left ( \int_z^{+ \infty} \widetilde{P_B} (y) \right ) \right | \leqslant \frac{C}{\omega} e^{- \sqrt{\omega} |x|/10}. \]
By the Cauchy-Schwarz inequality, we get
\[ \left | \int_{\R} e^{- \sqrt{\omega} |x|/10} \int_x^{+ \infty} \widetilde{P_B} (y) \int_x^y u'(z) u(z) \, \text{d}z \, \text{d}y \, \text{d}x \right | \leqslant \frac{C}{\omega} \left ( \int_{\R} u'(x)^2 e^{- \sqrt{\omega} |x|/10} \, \text{d}x \right )^{1/2} \left ( \int_{\R} u(x)^2 e^{- \sqrt{\omega} |x|/10} \, \text{d}x \right )^{1/2}. \]
Hence,
\[ \begin{array}{rcl} \displaystyle{\left ( \int_{\R} \widetilde{P_B} \right ) \int_{\R} u(x)^2 e^{- \sqrt{\omega} |x|/10} \, \text{d}x} & \leqslant & \displaystyle{\frac{C}{\sqrt{\omega}} \int_{\R} u^2 \widetilde{P_B} + \frac{C}{\omega} \left ( \int_{\R} u'(x)^2 e^{- \sqrt{\omega} |x|/10} \, \text{d}x \right )^{1/2} \left ( \int_{\R} u(x)^2 e^{- \sqrt{\omega} |x|/10} \, \text{d}x \right )^{1/2}} \\
\\ & \leqslant & \displaystyle{\frac{C}{\sqrt{\omega}} \int_{\R} u^2 \widetilde{P_B} + \frac{C}{\omega^2 \int_{\R} \widetilde{P_B}} \int_{\R} u'(x)^2 e^{- \sqrt{\omega}|x|/10} \, \text{d}x + \frac{\int_{\R} \widetilde{P_B}}{2} \int_{\R} u(x)^2 e^{- \sqrt{\omega} |x|/10} \, \text{d}x}, \end{array} \]
using Young's inequality in the last line. We finally get that
\[ \left ( \int_{\R} \widetilde{P_B} \right ) \int_{\R} u(x)^2 e^{- \sqrt{\omega} |x|/10} \, \text{d}x \leqslant \frac{C}{\sqrt{\omega}} \int_{\R} u^2 \widetilde{P_B} + \frac{C}{\omega^2 \int_{\R} \widetilde{P_B}} \int_{\R} u'(x)^2 e^{- \sqrt{\omega}|x|/10} \, \text{d}x. \]
Now recalling the definition of $\widetilde{P_B}$, we see that $\displaystyle{\int_{\R} \widetilde{P_B} = \frac{\gamma_B}{C \omega}}$. Also writing that $e^{- \sqrt{\omega} |x|/10} \leqslant 1$ in the second integral of the right side, and that $e^{- \sqrt{\omega} |x|/10} \geqslant \rho (x) /2$ in the integral of the left side, and multiplying the inequality above by $\varepsilon_\omega \omega^{3/2}$, we obtain the desired inequality. The proof for the analogous result with $B=\infty$ is identical. \hfill \qedsymbol

\noindent \textcolor{white}{a} \\ Now we prove that hypotheses $(H_1)$ and $(H_2)$ are sufficient to ensure there does not exist an internal mode in our problem.

\begin{leftbar}
\noindent \textbf{Proposition 2.} Assume that hypotheses $(H_1)$ and $(H_2)$ hold. Then, for $\omega$ small enough, there does not exist $V,W \in H^1 ( \R )$ and $\lambda \in \C$ such that 
\begin{equation}
    \left \{ \begin{array}{ccc} M_- V &=& \lambda W \\ M_+ W &=& \lambda V. \end{array} \right.
    \label{sysM}
\end{equation}
other than $V=W=0$.
\end{leftbar}

\noindent \textit{Proof.} Note that the hypothesis $(H_2)$ implies that, for $K_0>0$ any fixed positive constant, and for $\omega$ small enough (which is the case we will consider in what follows),
\[ \varepsilon_\omega \gamma_\infty = - \int_{\R} \frac{x (a_\omega^+ + a_\omega^-)'}{2} = \int_{\R} \frac{a_{\omega}^+ + a_{\omega}^-}{2} > K_0 \varepsilon_\omega^2 \sqrt{\omega}. \]
Starting with the system $\left \{ \begin{array}{ccc} M_- V &=& \lambda W \\ M_+ W &=& \lambda V \end{array} \right.$, we multiply the first line by $(2 \Phi_B V' + \Phi_B ' V)$, the second by $(2 \Phi_B W' + \Phi_B ' W)$, we integrate on $\R$ and we sum:
\[ \begin{array}{rcl} \displaystyle{\int_{\R} ( M_- V)(2 \Phi_B V' + \Phi_B ' V) + \int_{\R} ( M_+ W) (2 \Phi_B W' + \Phi_B ' W)} & = & \displaystyle{\lambda \int_{\R} ((WV'+VW') 2 \Phi_B + 2 \Phi_B ' VW)} \\
\\ &=& \displaystyle{\lambda \left [ 2VW \Phi_B \right ]_{- \infty}^{+ \infty} \, \, = \, \, 0.} \end{array} \]
Now, following virial computations (basically integrating by parts),
\[  \begin{array}{rcl} \displaystyle{\int_{\R} (M_- V)(2 \Phi_B V' + \Phi_B ' V)} &=& \displaystyle{\int_{\R} -V''(2 \Phi_B V' + \Phi_B ') + \omega \underbrace{\int_{\R} V(2 \Phi_B V' + \Phi_B 'V)}_{= \, 0} + \int_{\R} a_\omega^- V (2 \Phi_B V' + \Phi_B ' V)} \\
\\ &=& \displaystyle{\int_{\R} 2 ((\zeta_B V)')^2 + \int_{\R} ( \ln \zeta_B)'' ( \zeta_B V)^2 - \int_{\R} (a_\omega^-)' \Phi_B V^2.} \end{array} \]
Now let $B \to + \infty$. We recall that $V \in H^1 ( \R ) \subset L^\infty ( \R )$. First, $| ( \ln \zeta_B )'' (x) | \leqslant \frac{C \sqrt{\omega}}{B} \, \mathbbm{1}_{[1,2]} ( \sqrt{\omega} |x| ) \leqslant \frac{C}{B} \rho (x)$, thus $\int_{\R} ( \ln \zeta_B)'' ( \zeta_B V)^2 \longrightarrow 0$ as $B \to + \infty$. Moreover, since $\Phi_B (x) \longrightarrow x$ as $B \to + \infty$, the dominated convergence theorem shows that $\int_{\R} (a_\omega^-)' \Phi_B V^2 \longrightarrow \int_{\R} x(a_\omega^-)'V^2$ as $B \to + \infty$. Finally, note that $\zeta_B (x) \longrightarrow 1$ as $B \to + \infty$, $| \zeta_B '(x)| \leqslant \frac{C}{B} e^{-|x|/B}$ and $| \zeta_B ''(x)| \leqslant \frac{C}{B^2} e^{-|x| / B} + \frac{C}{B} \theta (x)$ where $\theta$ has a compact support that does not depend on $B$. Using these estimates and the dominated convergence theorem, we see that
\[ \int_{\R} ((\zeta_B V)')^2 = \int_{\R} \zeta_B^2 (V')^2 - \int_{\R} \zeta_B \zeta_B '' V^2 \, \underset{B \to + \infty}{\longrightarrow} \, \int_{\R} (V')^2. \]
Hence, 
\[ \int_{\R} (M_- V)(2 \Phi_B V' + \Phi_B ' V) \, \underset{B \to + \infty}{\longrightarrow} \, \int_{\R} 2 (V')^2 - \int_{\R} x (a_\omega^-)' V^2. \]
We have a similar formula involving $M_+ W$. Combining these two identities, we get
\begin{equation}
    0 = 2 \int_{\R} ((V')^2 + (W')^2) - \int_{\R} x(a_\omega^-)' V^2 - \int_{\R} x(a_\omega^+)' W^2. 
    \label{sp1}
\end{equation}
Now, let us take $R_\infty$ a bounded function that we will define later. Taking the initial system \eqref{sysM}, we multiply the first line by $R_\infty V$ and the second line by $R_\infty W$, before again integrating on $\R$ and taking the difference; this leads to
\[ \int_{\R} M_- V \cdot R_\infty V - \int_{\R} M_+ W \cdot R_\infty W = \lambda \int_{\R} R_\infty VW - \lambda \int_{\R} R_\infty VW = 0.  \]
We compute 
\[ \begin{array}{rcl} \displaystyle{\int_{\R} M_- V \cdot R_\infty V} &=& \displaystyle{\int_{\R} -V'' R_\infty V + \omega \int_{\R} R_\infty V^2 + \int_{\R} a_\omega^- R_\infty V^2} \\
\\ &=& \displaystyle{ \int_{\R} R_\infty (V')^2 - \int_{\R} \frac{R_\infty ''}{2} V^2 + \omega \int_{\R} R_\infty V^2 + \int_{\R} a_\omega^- R_\infty V^2.} \end{array} \]
Here too, we have a similar formula involving $M_+ W$. Taking the difference, we find
\begin{equation}
    0 = \int_{\R} \left ( \omega R_\infty - \frac{R_\infty ''}{2} \right ) (V^2 - W^2) + \int_{\R} R_\infty ((V')^2 - (W')^2) + \int_{\R} a_\omega^- R_\infty V^2 - \int_{\R} a_\omega^+ R_\infty W^2.
    \label{sp2}
\end{equation}
Now summing \eqref{sp1} and \eqref{sp2}, we get
\[ \begin{array}{rcl} 0 &=& \displaystyle{2 \int_{\R} ((V')^2 + (W')^2) + \int_{\R} \left ( -x(a_\omega^-)' + \omega R_\infty - \frac{R_\infty''}{2} \right ) V^2 + \int_{\R} \left ( - x( a_\omega^+)'  - \omega R_\infty + \frac{R_\infty''}{2} \right ) W^2} \\
\\ & & \, \, \, \, \displaystyle{+ \, \, \int_{\R} R_\infty ((V')^2 - (W')^2) + \int_{\R} a_\omega^- R_\infty V^2 - \int_{\R} a_\omega^+ R_\infty W^2.} \end{array} \]
Now, let us define $R_\infty$ as the bounded solution of the ordinary differential equation $- \frac{R_\infty''}{2} + \omega R_\infty = D_\infty$ where $D_\infty := - \frac{x(a_\omega^+ - a_\omega^-)'}{2}$. We finally obtain
\[ 0 = 2 \int_{\R} ((V')^2 + (W')^2) + \int_{\R} P_\infty (V^2 + W^2) + K_{2a} + K_{2b} \]
where $K_{2a} := \int_{\R} R_\infty ((V')^2 - (W')^2)$ and $K_{2b} := \int_{\R} a_\omega^- R_\infty V^2 - \int_{\R} a_\omega^+ R_\infty W^2$. \\
\\ By Lemma 7, we can assume $\omega$ small enough so that
\[ \int_{\R} P_\infty V^2 \geqslant C \gamma_\infty \varepsilon_\omega \sqrt{\omega} \int_{\R} \rho V^2 - \frac{C \varepsilon_\omega \sqrt{\omega}}{\gamma_\infty} \int_{\R} (V')^2 \]
and that the same inequality holds taking $W$ instead of $V$. Let us now control the error terms $J_1$, $K_{2a}$ and $K_{2b}$. \\
\\ About $K_{2a}$, we first see that $R_\infty$ is bounded and we can control this aspect. Indeed, the explicit expression of $R_\infty$ is given by the variation of the constants:
\[ R_\infty (x) = \frac{1}{\sqrt{2 \omega}} \left ( \int_{- \infty}^x e^{\sqrt{2 \omega} (y-x)} D_\infty (y) \, \text{d}y + \int_x^{+ \infty} e^{\sqrt{2 \omega} (x-y)} D_\infty (y) \, \text{d}y \right ). \]
Using this expression and the estimate $|D_\infty(x)| \leqslant C \varepsilon_{\omega} \omega^{3/2} |x| e^{- \sqrt{\omega} |x|}$, we show that $|R_\infty| \leqslant \frac{C}{\omega} |D_\infty| \leqslant C \varepsilon_\omega \rho^2$. This leads to
\[ |K_{2a}| \leqslant C \varepsilon_\omega \int_{\R} ((V')^2 + (W')^2). \]
About $K_{2b}$, we first recall that $|a_\omega^{\pm}| \leqslant \varepsilon_\omega \phi_\omega^2 \leqslant \varepsilon_\omega \omega \rho$. This and the estimate $||R_B||_{\infty} \leqslant C \varepsilon_\omega$ lead to
\[ |K_{2b}| \leqslant C \varepsilon_\omega^2 \omega \int_{\R} \rho (V^2 + W^2). \]
Putting all this together, we find that
\[ \begin{array}{rcl} 0 &=& \displaystyle{2 \int_{\R} ((V')^2 + (W')^2) + \int_{\R} P_B (V^2 + W^2) + K_{2a} + K_{2b}} \\
\\ & \geqslant & \displaystyle{2 \int_{\R} ((V')^2 + (W')^2) + C \varepsilon_\omega \gamma_\infty \sqrt{\omega} \int_{\R} \rho (V^2 + W^2) - C \varepsilon_\omega \frac{\sqrt{\omega}}{\gamma_\infty} \int_{\R} ((V')^2 + (W')^2)} \\
\\ & & \displaystyle{\, \, \, \, \, \, - \, C \varepsilon_\omega \int_{\R} ((V')^2 + (W')^2) - C \varepsilon_\omega^2 \omega \int_{\R} \rho (V^2 + W^2)} \\
\\ & \geqslant & \displaystyle{\left ( 2 - C \varepsilon_\omega \frac{\sqrt{\omega}}{\gamma_\infty} - C \varepsilon_\omega \right ) \int_{\R} ((V')^2 + (W')^2) + \left ( C \varepsilon_\omega \gamma_\infty \sqrt{\omega} - C \varepsilon_\omega^2 \omega \right ) \int_{\R} \rho (V^2 + W^2).} \end{array} \]
We first see that $2 - C \varepsilon_\omega \frac{\sqrt{\omega}}{\gamma_\infty} - C \varepsilon_\omega \geqslant 2 - \frac{C}{K_0} - C \varepsilon_\omega$. Thus we can assume $\omega$ small enough and $K_0$ large enough such that $2 - \frac{C}{K_0} - C \varepsilon_\omega \geqslant 1$. Note that $K_0$ does not depend on $\omega$. On the other hand, we see that
\[ C \varepsilon_\omega \gamma_\infty \sqrt{\omega} - C \varepsilon_\omega^2 \omega \geqslant K_0 \varepsilon_\omega^2 \omega - C \varepsilon_\omega^2 \omega = \varepsilon_\omega^2 \omega \left ( K_0 - C \right ). \]
We can assume $\omega$ small enough and $K_0$ large enough (still not depending on $\omega$) such that $K_0 - C \geqslant 1$ for instance. Putting all this together, we get
\[ 0 \geqslant \int_{\R} ((V')^2 + (W')^2) + \varepsilon_\omega^2 \omega \int_{\R} \rho (V^2 + W^2) \]
which leads to $V=W=0$. \hfill \qedsymbol

\noindent \textcolor{white}{a} \\ \noindent Before concluding the proof of Theorem 1, let us check, as announced in the previous section, that hypotheses $(H_1)$ and $(H_2)$ ensure that $\text{Ker} (M_-) = \{ 0 \}$. 

\begin{leftbar}
\noindent \textbf{Corollary 1.} Assume that hypotheses $(H_1)$ and $(H_2)$ hold. Then, for $\omega$ small enough, $\text{Ker} (M_-) = \{ 0 \}$.
\end{leftbar}

\noindent \textit{Proof.} Take $V \in \text{Ker} (M_-)$, $\lambda =0$ and $W=0$. We have $M_- V = \lambda W$ and $M_+ W = \lambda V$, thus Proposition 2 gives $V=0$. \hfill \qedsymbol

\noindent \textcolor{white}{a} \\ \noindent Now we can give the proof of Theorem 1. Let $X,Y \in H^1 ( \R )$ and $\lambda \in \C$ be solutions of the system \eqref{IntM} that we recall here:
\[ \left \{ \begin{array}{ccl} L_- X &=& \lambda Y \\ L_+ Y &=& \lambda X. \end{array} \right. \]
Thanks to this system we see that $X,Y \in H^6 ( \R )$. Then $M_+ M_- S^2 X = S^2 L_+ L_- X = \lambda^2 S^2 X$. Let $V := S^2 X$. First, assume $\lambda \neq 0$. Denoting $W := \lambda^{-1} M_- V$, we have
\[ \left \{ \begin{array}{ccl} M_- V &=& \lambda W \\ M_+ W &=& \lambda V. \end{array} \right. \]
Therefore we know from Proposition 2 that, providing $\omega$ is small enough, $V=W=0$. As $\text{Ker} (S^2) = \text{span} ( \phi_\omega \, , x \phi_\omega )$, the relation $S^2 X = 0$ gives $X = c_1 \phi_\omega + c_2 x \phi_\omega$. This gives $L_- X = -2 c_2 \phi_\omega '$. Hence, $Y = -2 c_2 \lambda^{-1} \phi_\omega '$. This leads to $L_+ Y = 0$ i.e. $X=0$ and then $Y=0$. \\
\\ Now, assume $\lambda = 0$. We have $L_- X = L_+ Y = 0$. Since $\text{Ker} (L_-) = \text{span} ( \phi_\omega )$ and $\text{Ker} (L_+) = \text{span} ( \phi_\omega ' )$, we get $X = c_1 \phi_\omega$ and $Y = c_2 \phi_\omega '$. Reciprocally, all of these are solutions of the system. This completes the proof of Theorem 1. \hfill \qedsymbol
\textcolor{white}{a} \\ \\ Theorem 1, which is now proved, shows that there does not exist internal modes under hypotheses $(H_1)$ and $(H_2)$. We can go a little further and show, with the same proof, that there does not exist \textit{resonances} under the same hypotheses, in the sense below. See \cite{Ge} for similar arguments on the Klein-Gordon equation.

\begin{leftbar}
\noindent \textbf{Corollary 2.} Assume that hypotheses $(H_1)$ and $(H_2)$ are satisfied and that $\omega$ is small enough. Let $(X \, , Y \, , \lambda)$ be a solution to the system \eqref{IntM}. Assume that $X,Y$ belong to $L^\infty$ and that $X',Y'$ belong to $H^1$. Such a solution is called a \textit{resonance}. Then, either $X=Y=0$; or $\lambda = 0$, $X \in \text{span} ( \phi_\omega )$ and $Y \in \text{span} ( \phi_\omega ')$. 
\end{leftbar}

\noindent \textit{Proof.} In Proposition 2, one can assume $V$ and $W$ to be $L^\infty$ with derivatives in $L^2$, the result remains true. Indeed, the integrals $\int_{\R} (V')^2$ and $\int_{\R} (W')^2$ still have a sense, and so have the other integrals since $V^2$ and $W^2$ are always integrated after multiplication by an appropriate weight. For instance, the virial computations hold thanks to the presence of $\zeta_B$ and $\Phi_B$; and the integrals $\int_{\R} P_\infty V^2$, $\int_{\R} R_\infty V^2$ or $\int_{\R} \rho V^2$ exist since $P_\infty$, $R_\infty$ and $\rho$ are $L^1$ while $V^2$ (and $W^2$) are $L^\infty$. Hence, Proposition 2 remains true after this change. \\
\\ Now, take $(X \, , Y \, , \lambda)$ a resonance in our problem. As in the proof of Theorem 1, assume first that $\lambda \neq 0$ and let $V := S^2X$ and $W := \lambda^{-1} M_- V$. We can compute
\[ S^2 = \partial_x^2 - 2 \frac{\phi_\omega '}{\phi_\omega} \cdot \partial_x + \omega - g(\phi_\omega^2) + 2 \frac{G(\phi_\omega^2)}{\phi_\omega^2}. \]
We know that $X' \in H^1 \subset L^\infty$, thus $V=S^2 X \in L^\infty$. Besides, deriving the relation $\lambda Y = L_- X$ we see that $X''' \in L^2$, which shows that
\[ V' = X''' - 2 \frac{\phi_\omega '}{\phi_\omega} X '' - 2 \left ( \frac{\phi_\omega '}{\phi_\omega} \right ) ' X' + \left ( \omega - g(\phi_\omega^2) + 2 \frac{G(\phi_\omega^2)}{\phi_\omega^2} \right ) X' + \left ( - g(\phi_\omega^2) + 2 \frac{G(\phi_\omega^2)}{\phi_\omega^2} \right ) ' \in L^2. \]
Similarly, we show that $W \in L^\infty$ and $W' \in L^2$. Now, thanks to the new version of Proposition 2, we obtain $V=W=0$. The relation $S^2X=0$ is nothing but a second order ordinary differential equation, therefore here too we find $X=c_1 \phi_\omega + c_2 x \phi_\omega$, then $Y = -2c_2 \lambda^{-1} \phi_\omega '$ and finally $X=Y=0$. \\
\\ Now assume $\lambda = 0$. We have $L_- X= L_+ Y=0$ but this time $X,Y$ are not supposed to be in $H^1$. However, $L_+ Y = 0$ leads to $Y \in \text{span} ( \phi_\omega ' \, , A_\omega )$ where $A_\omega$ is defined just before Lemma 3. Since $Y$ and $\phi_\omega '$ are bounded while $A_\omega$ is not, we get $Y \in \text{span} ( \phi_\omega ')$. The same argument holds for $X$ and we find that $X \in \text{span} ( \phi_\omega )$. This completes the proof of Corollary 2. \hfill \qedsymbol

\section{Asymptotic stability}

\subsection{Modulation decomposition}
\noindent We fix an initial data $\phi_\omega \in H^1 ( \R )$ close to $\phi_{\omega_0}$. By the orbital stability property we know that the global solution $\psi$ of \eqref{NLS} remains close to the family of solitary waves for all time. It is standard to decompose $\psi$ as
\[ \psi (t \, , y) = e^{i ( \beta (t) (y- \sigma (t)) +  \gamma (t))} \left [ \phi_{\omega (t)} (y - \sigma (t)) + u(t \, , y - \sigma (t)) \right ] \]
where the functions $\beta$, $\sigma$, $\gamma$ and $\omega$ are of class $\mathscr{C}^1$ (as functions of time) and uniquely fixed so that, for all $t \geqslant 0$, the following orthogonality relations hold:
\[ \langle u \, , \phi_\omega \rangle = \langle u \, , x \phi_\omega \rangle = \langle u \, , i \Lambda_\omega \rangle = \langle u \, , i \phi_\omega ' \rangle = 0. \]
This choice of orthogonality relations is known to lead to the following inequality, satisfied for all $t \geqslant 0$,
\begin{equation}
    \frac{| \dot{\beta} |}{\sqrt{\omega}} + \frac{| \dot{\omega} |}{\omega} + \sqrt{\omega} | \dot{\sigma} - 2 \beta | + | \dot{\gamma} - \omega - \beta^2| \leqslant C \sqrt{\omega} \left \| \sech ( \sqrt{ \omega} x/2) u \right \|^2 \leqslant C \sqrt{\omega} || \rho^2 u ||^2. 
    \label{orth}
\end{equation}
See \cite{We2}. Furthermore, the orbital stability result can be written as follows: for $\epsilon$ small and for all $t \geqslant 0$, 
\begin{equation} 
    || \partial_x u || + ||u|| + | \beta | + | \omega - \omega_0 | \leqslant \epsilon
    \label{orbstab}
\end{equation}
for $\psi_0$ taken sufficiently close to $\phi_{\omega_0}$. \\
\\ Write $u = u_1 + i u_2$. The equation \eqref{NLS} satisfied by $\psi$ leads to the following system satisfied by $(u_1 \, , u_2)$:
\begin{equation}
    \left \{ \begin{array}{ccl} \partial_t u_1 &=& L_- u_2 + \theta_2 + m_2 - q_2 \\ \partial_t u_2 &=& - L_+ u_1 - \theta_1 - m_1 + q_1 \end{array} \right.
    \label{Su}
\end{equation}
where
\[ \begin{array}{l} \theta_1 \, = \, \dot{\beta} x \phi_\omega + ( \dot{\gamma} - \omega - \beta^2 ) \phi_\omega - \beta ( \dot{\sigma} - 2 \beta ) \phi_\omega, \\ 
\\ \theta_2 \, = \, - \frac{\dot{\omega}}{\omega} \Lambda_\omega + ( \dot{\sigma} - 2 \beta ) \phi_\omega ', \\
\\ m_1 \, = \, \dot{\beta} x u_1 + ( \dot{\gamma} - \omega - \beta^2 ) u_1 - ( \dot{\sigma} - 2 \beta ) \partial_x u_2 - \beta ( \dot{\sigma} - 2 \beta ) u_1, \\ 
\\ m_2 \, = \, \dot{\beta} x u_2 + ( \dot{\gamma} - \omega - \beta^2 ) u_2 + ( \dot{\sigma} - 2 \beta ) \partial_x u_1 - \beta ( \dot{\sigma} - 2 \beta ) u_2, \\ 
\\ q_1 \, = \, \text{Re} \left [ h(\phi_\omega + u) - h ( \phi_\omega ) -  h' ( \phi_\omega ) u \right ], \\ 
\\ q_2 \, = \, \text{Im} \left [ h ( \phi_\omega + u ) - \frac{h ( \phi_\omega )}{\phi_\omega} \, u \right ] \end{array} \]
where $h(r) = |r|^2 r - g(|r|^2)r$.

\subsection{First virial estimate}
\noindent Since $| \omega - \omega_0 | \leqslant \epsilon$, we get, for $\epsilon < \frac{\omega_0}{2}$, that $\frac{\omega_0}{2} \leqslant \omega \leqslant \frac{3 \omega_0}{2}$. This enables to control $\phi_\omega$, $\Lambda_\omega$ and their derivatives by powers of $\rho$. More precisely, for instance, $\phi_\omega \leqslant C \sqrt{\omega} \, \rho^N$, $| \phi_\omega ' | \leqslant C \omega \rho^N$, $| \Lambda_\omega | \leqslant C \sqrt{\omega} \, \rho^N$ and $| \Lambda_\omega ' | \leqslant C \omega \rho^N$ for any $N \in [\![ 0 \, , 7 ]\!]$. 

\begin{leftbar}
\noindent \textbf{Proposition 3.} There exists $C>0$ such that, for $\epsilon$ small enough and any $T \geqslant 0$,
\[ \int_0^T \left ( || \eta_A \partial_x u ||^2 + \frac{1}{A^2} ||\eta_A u||^2 \right ) \, \text{d}t \leqslant C \epsilon + C \omega_0 \int_0^T || \rho^2 u ||^2 \, \text{d}t. \]
\end{leftbar}

\noindent \textit{Proof.} We will use a virial argument. Let $w = \zeta_A u$ and 
\[ \mathcal{I} = \int_{\R} u_1 \left ( 2 \Phi_A \partial_x u_2 + \Phi_A ' u_2 \right ). \]
From the equation \eqref{Su} and noticing that $\int_{\R} (2 \Phi_A \partial_x u_1 + \Phi_A ' u_1) u_1 = \int_{\R} (2 \Phi_A \partial_x u_2 + \Phi_A ' u_2 ) u_2 = 0$ (by integration by parts), we get that
\[ \begin{array}{rcl} \dot{\mathcal{I}} &=& \displaystyle{- \int_{\R} (2 \Phi_A \partial_x u_1 + \Phi_A ' u_1) \partial_x^2 u_1 - \int_{\R} (2 \Phi_A \partial_x u_2 + \Phi_A ' u_2 ) \partial_x^2 u_2} \\
\\ & & \, \, \, \displaystyle{+ \int_{\R} (2 \Phi_A \partial_x u_1 + \Phi_A ' u_1) ( \theta_1 + m_1 ) + \int_{\R} (2 \Phi_A \partial_x u_2 + \Phi_A ' u_2 ) ( \theta_2 + m_2 )} \\
\\ & & \, \, \, \displaystyle{- \, \, \text{Re} \left [ \int_{\R} \left ( 2 \Phi_A \partial_x \overline{u} + \Phi_A ' \overline{u} \right ) \left ( h( \phi_\omega + u) - h ( \phi_\omega ) \right ) \right ].} \end{array} \]
Integrating by parts, we get that, for $k \in \{ 1 \, , 2 \}$,
\[ - \int_{\R} \left ( 2 \Phi_A \partial_x u_k + \Phi_A ' u_k \right ) \partial_x^2 u_k = 2 \int_{\R} ( \partial_x w_k)^2 + \int_{\R} ( \ln \zeta_A )'' w_k^2 \]
where, after computations, $(\ln \zeta_A )'' = - \frac{|x|}{A} \left ( 1 - \chi ( \sqrt{\omega_0} x) \right ) \mathbbm{1}_{[1,2]} (\sqrt{\omega_0} x)$. We see that
\[ | ( \ln \zeta_A )'' (x) | \leqslant \frac{C \sqrt{\omega_0}}{A} \, \mathbbm{1}_{[1,2]} ( \sqrt{\omega_0} |x| ) \leqslant \frac{C \sqrt{\omega_0}}{A} \rho^4 (x).  \]
Thus, the first part of $\dot{\mathcal{I}}$ is controlled as follows:
\[ - \int_{\R} \left ( 2 \Phi_A \partial_x u_k + \Phi_A ' u_k \right ) \partial_x^2 u_k \geqslant 2 \int_{\R} ( \partial_x w_k )^2 - \frac{C \sqrt{\omega_0}}{A} ||\rho^2 w_k||^2. \]
Now, about the second term in $\dot{\mathcal{I}}$, we notice that, denoting $H(r)=\frac{|r|^4}{4} - \frac{G(|r|^2)}{2}$,
\[ \partial_x \text{Re} \left [ H( \phi_\omega + u) - H ( \phi_\omega ) - h( \phi_\omega ) u \right ] = \text{Re} \left [ ( \partial_x \overline{u} ) \left ( h ( \phi_\omega + u) - h ( \phi_\omega ) \right ) \right ] + \text{Re} \left [ \phi_\omega ' \left ( h ( \phi_\omega + u) - h ( \phi_\omega ) - h ' ( \phi_\omega ) u \right ) \right ]. \]
Now integrating by parts, we decompose
\[ - \text{Re} \left [ \int_{\R} \left ( 2 \Phi_A \partial_x \overline{u} + \Phi_A ' \overline{u} \right ) \left ( h( \phi_\omega + u) - h ( \phi_\omega ) \right ) \right ] = I_1 + I_2 + I_3 \]
with
\[ \begin{array}{rl} & \displaystyle{I_1 = 2 \int_{\R} \Phi_A ' \, \text{Re} \left [ H ( \phi_\omega + u ) - H ( \phi_\omega ) - h( \phi_\omega )u \right ],} \\
\\ & \displaystyle{I_2 = 2 \int_{\R} \Phi_A \, \text{Re} \left [ \phi_\omega ' \left ( h ( \phi_\omega + u ) - h ( \phi_\omega ) - h ' ( \phi_\omega ) u \right ) \right ],} \\
\\ \text{and} & \displaystyle{I_3 = - \int_{\R} \Phi_A ' \, \text{Re} \left [ \overline{u} \left ( h ( \phi_\omega + u ) - h ( \phi_\omega ) \right ) \right ].} \end{array} \]
We recall that $\Phi_A ' = \zeta_A^2$. We note that $0 < \Phi_A ' \leqslant 1$ and $| \Phi_A (x) | \leqslant |x|$ on $\R$. Therefore, 
\[ |\Phi_A (x) \phi_\omega (x) | \leqslant \sqrt{\omega} |x| \sech ( \sqrt{\omega} x ) \leqslant C \rho^4 (x). \]
Now, about $I_1$, using the definitions of $H$ and $h$ and developing $|\phi_\omega + u|^4$ we compute that
\[ \begin{array}{rcl} \text{Re} \left [ H ( \phi_\omega + u ) - H ( \phi_\omega ) - h( \phi_\omega )u \right ] &=& \displaystyle{\frac{|u|^4}{4} + \phi_\omega^2 \text{Re} (u)^2 + \frac{\phi_\omega^2 |u|^2}{2} + \frac{\phi_\omega |u|^2 \text{Re} (u)}{2}} \\
\\ & & \displaystyle{\, \, \, \, \, - \frac{G ( |\phi_\omega + u|^2 )}{2} - \frac{| \phi_\omega |^4}{4} + \frac{G ( \phi_\omega^2 )}{2} + \phi_\omega g ( \phi_\omega^2 ) \text{Re} (u).} \end{array} \]
Now, $G$ is real-valued and we can write Taylor's expansion:
\[ G ( | \phi_\omega + u|^2 ) = G ( \phi_\omega^2 + |u|^2 + 2 \phi_\omega \text{Re} (u)) = G ( \phi_\omega^2) + \left ( |u|^2 + 2 \phi_\omega \text{Re}(u) \right ) g ( \phi_\omega^2 ) + \int_{\phi_\omega^2}^{| \phi_\omega + u|^2} \left ( | \phi_\omega +u|^2 - t \right ) g'(t) \, \text{d}t \]
where $\displaystyle{\left | \int_{\phi_\omega^2}^{| \phi_\omega + u|^2} \left ( | \phi_\omega +u|^2 - t \right ) g'(t) \, \text{d}t \right |} \leqslant C \left | |\phi_\omega + u|^2 - \phi_\omega^2 \right | \leqslant C |u|^4 + C \phi_\omega |u|^2 |\text{Re} (u)| + C \phi_\omega^2 | \text{Re} (u)|^2$. Putting these estimations together and using the inequalities $| \text{Re} (u) | \leqslant |u|$ and $\phi_\omega |u|^3 = ( \phi_\omega |u|) (|u|^2) \leqslant \frac{\phi_\omega^2 |u|^2}{2} + \frac{|u|^4}{2}$, we ultimately find that
\[ |I_1| \leqslant C \int_{\R} \Phi_A ' ( \phi_\omega^2 |u|^2 + |u|^4 ) = C \int_{\R} \zeta_A^2 ( \phi_\omega^2 |u|^2 + |u|^4 ) \leqslant C \omega_0 \int_{\R} \rho^4 |u|^2 + C \int_{\R} \zeta_A^2 |u|^4, \]
using $\phi_\omega \leqslant C \sqrt{\omega} \, \rho^2 \leqslant C \sqrt{\omega_0} \, \rho^2$ and $\zeta_A^2 \leqslant 1$ to control the first term. The control of the third term is similar, writing
\[ h ( \phi_\omega + u) - h ( \phi_\omega) = |u|^2 \phi_\omega + |u|^3 + \phi_\omega^2 u + 2 \phi_\omega^2 \text{Re} (u) + 2 \phi_\omega u \text{Re} (u) - u g ( |\phi_\omega + u|^2 ) - \phi_\omega \left ( g (|\phi_\omega + u|^2 ) - g ( \phi_\omega^2 ) \right ). \]
Using the inequalities $|u g ( |\phi_\omega +u|^2)| \leqslant C |u| \, |\phi_\omega + u|^2 \leqslant C|u| \phi_\omega^2 + C |u|^3$ and $|g(|\phi_\omega +u|^2) - g(\phi_\omega^2)| \leqslant C \left | |\phi_\omega + u|^2 - \phi_\omega^2 \right | = C |u|^2 + 2C \phi_\omega \text{Re} (u)$, we find that $\left | \text{Re} \left [ \overline{u} \left ( h ( \phi_\omega + u ) - h ( \phi_\omega ) \right ) \right ] \right | \leqslant C ( \phi_\omega^2 |u|^2 + |u|^4 )$ and then we can control $I_3$ the same way we controlled $I_1$:
\[ |I_3| \leqslant C \omega_0 \int_{\R} \rho^4 |u|^2 + C \int_{\R} \zeta_A^2 |u|^4. \]
About $I_2$, we compute that
\[ \begin{array}{rcl} h(\phi_\omega + u) - h ( \phi_\omega) - h'(\phi_\omega ) u &=& \phi_\omega |u|^2 + |u|^2 u + 2 \phi_\omega u \text{Re} (u) - (\phi_\omega + u) g(\phi_\omega^2 + |u|^2 + 2 \phi_\omega \text{Re}(u) ) \\ \\ & & \, \, \, \, + (\phi_\omega + u) g(\phi_\omega^2) + 2 \phi_\omega^2 u g'(\phi_\omega^2). \end{array} \]
Using Taylor's expansion formula, we write that
\[ g ( \phi_\omega^2 + |u|^2 + 2 \phi_\omega \text{Re} (u)) = g ( \phi_\omega^2 ) + (|u|^2 + 2 \phi_\omega \text{Re}(u)) g'(\phi_\omega^2) + \underbrace{\int_{\phi_\omega^2}^{\phi_\omega^2 + |u|^2 + 2 \phi_\omega \text{Re}(u)} (\phi_\omega^2 + |u|^2 + 2 \phi_\omega \text{Re}(u) - s) g''(s) \, \text{d}s}_{=: \, \text{IR}} \]
where we control the integral term IR as follows, recalling that $g''(s) = \mathcal{O} (1/s)$ since $(H_1)$ holds,
\[ \begin{array}{rcl} \displaystyle{\left | \int_{\phi_\omega^2}^{\phi_\omega^2 + |u|^2 + 2 \phi_\omega \text{Re}(u)} (\phi_\omega^2 + |u|^2 + 2 \phi_\omega \text{Re}(u) - s) g''(s) \, \text{d}s \right |} & \leqslant & \displaystyle{\left | |u|^2 + 2 \phi_\omega \text{Re} (u) \right | \, \left | \int_{\phi_\omega^2}^{\phi_\omega^2 + |u|^2 + 2 \phi_\omega \text{Re}(u)} \frac{C \, \text{d}s}{s} \right |} \\
\\ & \leqslant & \displaystyle{C \left | |u|^2 + 2 \phi_\omega \text{Re} (u) \right | \, \left | \ln \left ( 1 + \frac{|u|^2}{\phi_\omega^2} + \frac{2 \text{Re} (u)}{\phi_\omega} \right ) \right |.} \end{array} \]
We know that $\ln (1 + \cdot )$ is $1$-Lipschitz on $[0 \, , + \infty )$. We can say that this function is $C$-Lipschitz on $[-1/2 \, , + \infty )$ for example. We shall separate two cases. First, assume that $\left | \frac{u}{\phi_\omega} \right | \leqslant \frac{1}{4}$. Then $\frac{|u|^2}{\phi_\omega^2} + \frac{2 \text{Re}(u)}{\phi_\omega} \geqslant - \frac{1}{2}$ and we have 
\[ \begin{array}{rcl} \displaystyle{\text{IR}} & \leqslant & \displaystyle{C \left | |u|^2 + 2 \phi_\omega \text{Re} (u) \right | \, \left | \frac{|u|^2}{\phi_\omega^2} + \frac{2 \text{Re} (u)}{\phi_\omega} \right | \, \, \leqslant \, \, \frac{C}{\phi_\omega^2} \left | |u|^2 + 2 \phi_\omega \text{Re} (u) \right |^2} \\
\\ & \leqslant & \displaystyle{\frac{C}{\phi_\omega^2} \left ( |u|^4 + \phi_\omega^2 |u|^2 \right ) \, \, \leqslant \, \, C |u|^2} \end{array} \]
recalling, for the last inequality, that $|u / \phi_\omega | \leqslant C$. This gives
\[ \begin{array}{rl} & \left | - (\phi_\omega + u) g(\phi_\omega^2 + |u|^2 + 2 \phi_\omega \text{Re}(u) ) + (\phi_\omega + u) g(\phi_\omega^2) + 2 \phi_\omega^2 u g'(\phi_\omega^2) \right | \\
\\ = & \left | - |u|^2 ( \phi_\omega + \text{Re}(u)) g'(\phi_\omega^2) - 2 \phi_\omega \text{Re} (u)^2 g'(\phi_\omega^2) + \text{IR} \cdot (\phi_\omega + \text{Re} (u)) \right | \\
\\ \leqslant & C ( |u|^2 \phi_\omega + |u|^3 ) + ( \phi_\omega + |u| ) |\text{IR}| \\
\\ \leqslant & C ( |u|^2 \phi_\omega + |u|^3 ). \end{array} \]
This leads to $\left | \text{Re} \left ( h(\phi_\omega + u) - h ( \phi_\omega ) - h'(\phi_\omega ) u \right ) \right | \leqslant C ( \phi_\omega |u|^2 + |u|^3 )$. \\
\\ Now, assume that $\left | \frac{u}{\phi_\omega} \right | > \frac{1}{4}$. We have $\phi_\omega \leqslant C |u|$ and everything is easier. Using $|g(s)| \leqslant Cs$ and $|g'(s)| \leqslant C$, we see that $\left | \text{Re} \left ( h(\phi_\omega + u) - h ( \phi_\omega ) - h'(\phi_\omega ) u \right ) \right | \leqslant C ( \phi_\omega |u|^2 + |u|^3 )$ in this case too. \\
\\ Hence, whatever case we are in, we have the inequality above and thus,
\[ |I_2| \leqslant C \int_{\R} | \Phi_A \phi_\omega ' | ( \phi_\omega |u|^2 + |u|^3 ) \leqslant C \omega_0 \int_{\R} \rho^4 |u|^2, \]
using the inequalities $|\Phi_A \phi_\omega '| \leqslant |x \phi_\omega '| \leqslant C \sqrt{\omega} \rho^4$, $\phi_\omega \leqslant C \sqrt{\omega_0}$ and $|u| \leqslant C \omega_0 \leqslant C \sqrt{\omega_0}$. This last inequality follows from Sobolev embedding. Indeed, by the orbital stability property, we have $||u||_{H^1 ( \R )} \leqslant C \epsilon$ and thus, by Sobolev embedding, $||u||_{L^\infty} \leqslant C ||u||_{H^1 ( \R )} \leqslant C \epsilon \leqslant C \omega_0$. \\
\\ Now, we put the estimates on $I_1$, $I_2$ and $I_3$ together and we use the following inequality (see \cite{Ma2} or \cite{Ma1}):
\[  \int_{\R} \zeta_A^2 |u|^4 \leqslant CA^2 ||u||_{L^\infty}^2 \int_{\R} | \partial_x w|^2 \leqslant CA^2 \epsilon^2 \int_{\R} | \partial_x w|^2. \]
We then obtain that
\[ |I_1| + |I_2| + |I_3| \leqslant C \omega_0 \int_{\R} \rho^4 |u|^2 + CA^2 \epsilon^2 \int_{\R} | \partial_x w |^2. \]
Now, we integrate by parts to see that, for $k \in \{ 1 \, , 2 \}$,
\[ \left | \int_{\R} (2 \Phi_A \partial_x u_k + \Phi_A ' u_k ) \theta_k \right | = \left | \int_{\R} u_k ( 2 \Phi_A \partial_x \theta_k + \Phi_A ' \theta_k ) \right | \leqslant C ||u||_{L^\infty} \int_{\R} ( |x| \, |\partial_x \theta_k| + |\theta_k| ), \]
using that $|\Phi_A (x)| \leqslant |x|$ and $|\Phi_A'| \leqslant 1$. Now, recalling the expressions of $\theta_k$, we see that $\partial_x \theta_1 = \dot{\beta} \phi_\omega + \dot{\beta} x \phi_\omega ' + ( \dot{\gamma} - \omega - \beta^2) \phi_\omega ' - \beta ( \dot{\sigma} - 2 \beta ) \phi_\omega '$ and $\partial_x \theta_2 = - \frac{\dot{\omega}}{\omega} \Lambda_\omega ' + ( \dot{\sigma} - 2 \beta ) \phi_\omega ''$. Using that all of the functions $\phi_\omega$, $x \phi_\omega '$, $\phi_\omega '$, $\phi_\omega ''$ and $\Lambda_\omega '$ are bounded (by $C$, independent of $\omega$ and $\epsilon$), we see that
\[ \int_{\R} (|x| \, |\partial_x \theta_k| + | \theta_k | ) \leqslant C ||\rho^2 u||^2 \]
using \eqref{orth} and the fact that $\beta$ is bounded. Thus we get
\[ \left | \int_{\R} (2 \Phi_A \partial_x u_k + \Phi_A ' u_k ) \theta_k \right | \leqslant C \epsilon ||\rho^2 u||^2 \leqslant C \omega_0 ||\rho^2 u||^2. \]
The last terms remaining in the expression of $\dot{\mathcal{I}}$ are $\displaystyle{\int_{\R} (2 \Phi_A \partial_x u_k + \Phi_A ' u_k)m_k}$. Integrating by parts using the expression of $m_1$, we get
\[ - \int_{\R} (2 \Phi_A \partial_x u_1 + \Phi_A ' u_1) m_1 = \dot{\beta} \int_{\R} \Phi_A u_1^2 + ( \dot{\sigma} - 2 \beta ) \int_{\R} (2 \Phi_A \partial_x u_1 + \Phi_A ' u_1 ) \partial_x u_2. \]
Combining this identity with the corresponding identity for $\int_{\R} ( 2 \Phi_A \partial_x u_2 + \Phi_A ' u_2 )m_2$, we get
\[ - \int_{\R} (2 \Phi_A \partial_x u_1 + \Phi_A ' u_1) m_1 - \int_{\R} ( 2 \Phi_A \partial_x u_2 + \Phi_A ' u_2 )m_2 = \dot{\beta} \int_{\R} \Phi_A |u|^2 + ( \dot{\sigma} - 2 \beta ) \int_{\R} \Phi_A ' ( u_2 \partial_x u_1 - u_1 \partial_x u_2 ). \]
Therefore, using the upper bounds $||\Phi_A||_{L^\infty} \leqslant C A$, $|\Phi_A '| \leqslant 1$, $||u||, || \partial_x u|| \leqslant C \epsilon$, \eqref{orth} and the fact that $A > \frac{1}{\sqrt{\omega_0}}$, we find that
\[ \left | \int_{\R} (2 \Phi_A \partial_x u_1 + \Phi_A ' u_1) m_1 + \int_{\R} ( 2 \Phi_A \partial_x u_2 + \Phi_A ' u_2 )m_2 \right | \leqslant C A \epsilon^2 \sqrt{\omega_0} ||\rho^2 u||^2. \]
Putting all these estimates together, noticing that $|| \rho^2 w || \leqslant || \rho^2 u||$ and taking $\epsilon$ small enough so that $C A^2 \epsilon^2 \leqslant \frac{1}{2}$ (which also implies that $CA \epsilon^2 \sqrt{\omega_0} \leqslant CA^2 \epsilon^2 \frac{\sqrt{\omega_0}}{A} \leqslant \frac{\sqrt{\omega_0}}{2A} \leqslant \frac{\omega_0}{2}$), we get that
\[ \dot{\mathcal{I}} \geqslant \left ( 2 - CA^2 \epsilon^2 \right ) \int_{\R} | \partial_x w|^2 - C \left ( \omega_0 + \frac{\sqrt{\omega_0}}{A} + A \epsilon^2 \sqrt{\omega_0} \right ) ||\rho^2 u||^2 \geqslant \int_{\R} | \partial_x w|^2 - C \omega_0 || \rho^2 u ||^2. \]
This being established, we can conclude the proof. For any $T \geqslant 0$, the above estimates for $\Phi_A$ and \eqref{orbstab} give, by definition of $\mathcal{I}$,
\[ | \mathcal{I}(T) | \leqslant C ( || \Phi_A ||_{L^\infty} + || \Phi_A ' ||_{L^\infty} ) ||u(T)||_{H^1 ( \R )}^2 \leqslant CA \epsilon^2 \leqslant C \epsilon \]
providing we take $\epsilon$ small enough (which we assume from now on). Integrating on $[0 \, , T]$ the inequality satisfied by $\dot{\mathcal{I}}$, we get
\[ \int_0^T \int_{\R} | \partial_x w |^2 \leqslant \underbrace{\int_0^T \dot{\mathcal{I}}}_{\leqslant \, | \mathcal{I} (T)| + | \mathcal{I} (0)|} + \, \, C \omega_0 \int_0^T || \rho^2 u ||^2 \leqslant C \epsilon + C \omega_0 \int_0^T ||\rho^2 u||^2. \]
Now recall the following inequality from \cite{Ma2} or \cite{Ma1}:
\[ \int_{\R} \eta_A |w|^2 \leqslant CA^2 \int_{\R} | \partial_x w|^2 + CA \sqrt{\omega_0} \int_{\R} \rho^4 |w|^2, \]
which implies
\[ \frac{1}{A^2} \int_0^T \int_{\R} \eta_A^2 |u|^2 \leqslant C \epsilon + C \omega_0 \int_0^T || \rho^2 u ||^2 \]
using $\eta_A \leqslant C \zeta_A^2$ and $1/A < \sqrt{\omega_0}$. Now, recalling $w= \zeta_A u$ and writing that $| \zeta_A|^3 | \zeta_A '| \, |u \partial_x u| \leqslant \frac{1}{4} \zeta_A^4 | \partial_x u |^2 + 4 \zeta_A^2 | ( \zeta_A')^2 |u|^2$, we find that
\[ \begin{array}{rcl} \displaystyle{\int_{\R} \zeta_A^2 | \partial_x w|^2} &=& \displaystyle{\int_{\R} \zeta_A^2 \left | \zeta_A \partial_x u + \zeta_A ' u \right |^2} \\ \\ & \geqslant & \displaystyle{\int_{\R} \zeta_A^4 | \partial_x u|^2 - 2 \int_{\R} \zeta_A^3 | \zeta_A ' | \, |u \partial_x u | - \int_{\R} \zeta_A^2 ( \zeta_A ')^2 |u|^2} \\ \\ & \geqslant & \displaystyle{\frac{1}{2} \int_{\R} \zeta_A^4 | \partial_x u |^2 - 9 \int_{\R} \zeta_A^2 ( \zeta_A')^2 |u|^2} \end{array} \]
and thus, using the inequalities $\frac{1}{C} \eta_A \leqslant \zeta_A^2 \leqslant C \eta_A$ and $| \zeta_A'| \leqslant \frac{C}{A} \zeta_A$, we obtain
\[ \int_{\R} \eta_A^2 | \partial_x u|^2 \leqslant C \int_{\R} | \partial_x w|^2 + \frac{C}{A^2} \int_{\R} \eta_A^2 |u|^2. \]
Integrating over $[0 \, , T]$ and combining with the previous inequalities, we finally find that
\[ \begin{array}{rcl} \displaystyle{\int_0^T \left ( || \eta_A \partial_x u ||^2 + \frac{1}{A^2} || \eta_A u ||^2 \right ) \, \text{d}t} & \leqslant & \displaystyle{ C \int_0^T \int_{\R} | \partial_x w |^2 + \frac{C}{A^2} \int_0^T \int_{\R} \eta_A^2 |u|^2} \\
\\ & \leqslant & \displaystyle{C \epsilon + C \omega_0 \int_0^T || \rho^2 u ||^2 \, \text{d}t,} \end{array} \]
which is the desired result. \hfill \qedsymbol

\subsection{Transformed problem}
\noindent We will later fix a certain $\alpha > 0$, chosen small. For this $\alpha$ we introduce $v_1 = X_\alpha^2 M_- S^2 u_2$, $v_2 = - X_\alpha^2 S^2 L_+ u_1$ and $v = v_1 + iv_2$. We recall that
\[ S^2 = \partial_x^2 - 2 \, \frac{\phi_\omega '}{\phi_\omega} \cdot \partial_x + \omega - g ( \phi_\omega^2 ) + 2 \, \frac{G ( \phi_\omega^2 )}{\phi_\omega^2}. \]
We then compute
\[ \begin{array}{rcl} M_- S^2 &=& - \partial_x^4 + 2 \partial_x^2 \cdot \frac{\phi_\omega '}{\phi_\omega} \cdot \partial_x + \partial_x \cdot \left ( 2 \phi_\omega^2 g'(\phi_\omega^2) - 4 g ( \phi_\omega^2 ) + 4 \, \frac{G ( \phi_\omega^2)}{\phi_\omega^2} \right ) \cdot \partial_x \\
\\ & & \, \, + \left ( 4 \phi_\omega \phi_\omega ' g' ( \phi_\omega^2 ) - 6 \, \frac{\phi_\omega '}{\phi_\omega} \, g( \phi_\omega^2 ) - 4 \phi_\omega ' \phi_\omega^3 g''( \phi_\omega^2 ) + 4 \frac{\phi_\omega'}{\phi_\omega} \frac{G ( \phi_\omega^2 )}{\phi_\omega^2} - 2 \omega \, \frac{\phi_\omega '}{\phi_\omega} \right ) \cdot \partial_x \\
\\ & & \, \, + \, \omega^2 + 2 \omega \left ( g ( \phi_\omega^2 ) - \phi_\omega^2 g'(\phi_\omega^2) + 2 \phi_\omega^4 g''( \phi_\omega^2) \right ) \\
\\ & & \, \, - \, 2 g'(\phi_\omega^2) G ( \phi_\omega^2 ) + \phi_\omega^4 g'(\phi_\omega^2) - 2 \phi_\omega^6 g''( \phi_\omega^2) + 4 \phi_\omega^2 G(\phi_\omega^2) g''(\phi_\omega^2) \\
\\ & & \, \, - \, 2 \phi_\omega^2 g ( \phi_\omega^2 ) + 2 G ( \phi_\omega^2 ) + g ( \phi_\omega^2 )^2 \end{array} \]
and
\[ \begin{array}{rcl} S^2 L_+ &=& - \partial_x^4 + 2 \partial_x^2 \cdot \frac{\phi_\omega '}{\phi_\omega} \cdot \partial_x + \partial_x \cdot \left ( - \phi_\omega^2 - 2 g ( \phi_\omega^2 ) + 2 \, \frac{G ( \phi_\omega^2 )}{\phi_\omega^2} + 2 \phi_\omega^2 g'( \phi_\omega^2 ) \right ) \cdot \partial_x \\
\\ & & \, \, + \left ( -2 \phi_\omega \phi_\omega ' + 4 \phi_\omega \phi_\omega ' g' ( \phi_\omega^2 ) - 2 \, \frac{\phi_\omega '}{\phi_\omega} \, g( \phi_\omega^2 ) + 4 \phi_\omega ' \phi_\omega^3 g''( \phi_\omega^2 ) - 2 \omega \, \frac{\phi_\omega '}{\phi_\omega} \right ) \cdot \partial_x \\
\\ & & \, \, + \, \omega^2 + \omega \left ( -3 \phi_\omega^2 + 20 \phi_\omega^4 g''(\phi_\omega^2) + 8 \phi_\omega^6 g'''(\phi_\omega^2) + 2 \phi_\omega^2 g'(\phi_\omega^2) + 2 \, \frac{G ( \phi_\omega^2)}{\phi_\omega^2} \right ) \\
\\ & & \, \, + \, 3 \phi_\omega^4 - 3 \phi_\omega^2 g ( \phi_\omega^2) - 3 \phi_\omega^4 g'(\phi_\omega^2) + 4 \phi_\omega^2 g ( \phi_\omega^2 ) g'( \phi_\omega^2 ) -2 g'( \phi_\omega^2 ) G( \phi_\omega^2 ) \\
\\ & & \, \, - \, 12 \phi_\omega^6 g''(\phi_\omega^2) + 16 \phi_\omega^2 G( \phi_\omega^2 ) g''(\phi_\omega^2) + 4 \phi_\omega^4 g( \phi_\omega^2 ) g''( \phi_\omega^2 ) - 4 \phi_\omega^8 g'''(\phi_\omega^2) \\
\\ & & \, \, + \, 8 \phi_\omega^4 G ( \phi_\omega^2 ) g''' ( \phi_\omega^2 ) - g ( \phi_\omega^2 )^2 + 2 g ( \phi_\omega^2 ) \, \frac{G ( \phi_\omega^2 )}{\phi_\omega^2}. \end{array} \]
We introduce the operators $Q_-$ and $Q_+$, obtained respectively from $M_- S^2$ and $S^2 L_+$ by differentiation with respect to $\omega$ and then multiplication by $\omega$. Their exact expressions are given below.
\[ \begin{array}{rcl} Q_- &=& 2 \partial_x^2 \cdot \left ( \frac{\Lambda_\omega ' \phi_\omega - \phi_\omega ' \Lambda_\omega}{\phi_\omega^2} \right ) \cdot \partial_x + \partial_x  \cdot \left ( -4 \phi_\omega \Lambda_\omega g'( \phi_\omega^2 ) + 4 \phi_\omega^3 \Lambda_\omega g'' ( \phi_\omega^2 ) + 8 \frac{\Lambda_\omega g( \phi_\omega^2 )}{\phi_\omega} - 8 \frac{\Lambda_\omega G ( \phi_\omega^2 )}{\phi_\omega^3} \right ) \cdot \partial_x \\
\\ & & \, \, + \, \left ( 4 \Lambda_\omega ' \phi_\omega g'( \phi_\omega^2 ) -8 \Lambda_\omega \phi_\omega ' g'( \phi_\omega^2 ) -4 \Lambda_\omega \phi_\omega ' \phi_\omega^2 g''( \phi_\omega^2 ) -4 \Lambda_\omega ' \phi_\omega^3 g'' ( \phi_\omega^2 ) -8 \Lambda_\omega \phi_\omega ' \phi_\omega^4 g'''(\phi_\omega^2) \right. \\
\\ & & \, \, \, \, \, \, \, \, \, \, \, \, \, \, \left. -6 \frac{\Lambda_\omega ' g ( \phi_\omega^2 )}{\phi_\omega} + 4 \frac{\Lambda_\omega ' G ( \phi_\omega^2 )}{\phi_\omega^3} + 14 \frac{\Lambda_\omega \phi_\omega ' g( \phi_\omega^2 )}{\phi_\omega^2} -12 \frac{\Lambda_\omega \phi_\omega ' G (\phi_\omega^2)}{\phi_\omega^4} -2 \omega \frac{\phi_\omega '}{\phi_\omega} -2 \omega \frac{\Lambda_\omega ' \phi_\omega - \Lambda_\omega \phi_\omega '}{\phi_\omega^2} \right ) \cdot \partial_x \\
\\ & & \, \, + \, 2 \omega^2 + 2 \omega \left ( g ( \phi_\omega^2 ) - \phi_\omega^2 g'( \phi_\omega^2 ) + 2 \phi_\omega^4 g''( \phi_\omega^2 ) \right ) + 4 \omega \left ( 3 \Lambda_\omega \phi_\omega^3 g''( \phi_\omega^2 ) + 2 \Lambda_\omega \phi_\omega^5 g'''(\phi_\omega^2 ) \right ) \\
\\ & & \, \, + \, 4 \Lambda_\omega \phi_\omega g''( \phi_\omega^2 ) G ( \phi_\omega^2 ) -10 \Lambda_\omega \phi_\omega^5 g''( \phi_\omega^2 ) -4 \Lambda_\omega \phi_\omega^7 g'''( \phi_\omega^2 ) + 8 \Lambda_\omega \phi_\omega^3 g( \phi_\omega^2) g''( \phi_\omega^2) + 8 \Lambda_\omega \phi_\omega^3 G ( \phi_\omega^2 ) g'''( \phi_\omega^2 ) \end{array} \]
and
\[ \begin{array}{rcl} Q_+ &=& 2 \partial_x^2 \cdot \left ( \frac{\Lambda_\omega ' \phi_\omega - \phi_\omega ' \Lambda_\omega}{\phi_\omega^2} \right ) \cdot \partial_x + \partial_x  \cdot \left ( -2 \Lambda_\omega \phi_\omega + 4 \frac{\Lambda_\omega g( \phi_\omega^2 )}{\phi_\omega} - 4 \frac{\Lambda_\omega G ( \phi_\omega^2 )}{\phi_\omega^3} + 4 \Lambda_\omega \phi_\omega^3 g'' ( \phi_\omega^2 ) \right ) \cdot \partial_x \\
\\ & & \, \, + \, \left ( -2 \Lambda_\omega \phi_\omega ' -2 \Lambda_\omega ' \phi_\omega + 4 \Lambda_\omega ' \phi_\omega g'( \phi_\omega^2 ) + 20 \Lambda_\omega \phi_\omega ' \phi_\omega^2 g''( \phi_\omega^2 ) + 4 \Lambda_\omega ' \phi_\omega^3 g''(\phi_\omega^2) \right. \\
\\ & & \, \, \, \, \, \, \, \, \, \, \, \, \, \, \left.  + 8 \Lambda_\omega \phi_\omega ' \phi_\omega^4 g'''(\phi_\omega^2 ) -2 \frac{\Lambda_\omega ' g(\phi_\omega^2)}{\phi_\omega} + 2 \frac{\Lambda_\omega \phi_\omega ' g ( \phi_\omega^2 )}{\phi_\omega^2} - 2 \omega \frac{\phi_\omega '}{\phi_\omega} - 2 \omega \frac{\Lambda_\omega ' \phi_\omega - \Lambda_\omega \phi_\omega '}{\phi_\omega^2} \right ) \cdot \partial_x \\
\\ & & \, \, + \, 2 \omega^2 + \omega \left ( -3 \phi_\omega^2 + 20 \phi_\omega^4 g''(\phi_\omega^2) + 8 \phi_\omega^6 g'''( \phi_\omega^2 ) + 2 \phi_\omega^2 g'(\phi_\omega^2) + 2 \frac{G ( \phi_\omega^2 )}{\phi_\omega^2} \right ) \\
\\ & & \, \, + \, 2 \omega \left ( -3 \Lambda_\omega \phi_\omega + 42 \Lambda_\omega \phi_\omega^3 g''(\phi_\omega^2) + 44 \Lambda_\omega \phi_\omega^5 g'''(\phi_\omega^2) + 8 \Lambda_\omega \phi_\omega^7 g''''(\phi_\omega^2) \right. \\
\\ & & \, \, \, \, \, \, \, \, \, \, \, \, \, \, \left. + 2 \Lambda_\omega \phi_\omega g'(\phi_\omega^2) + 2 \frac{\Lambda_\omega g(\phi_\omega^2)}{\phi_\omega} - 4 \frac{\Lambda_\omega G ( \phi_\omega^2 )}{\phi_\omega^3} \right ) \\
\\ & & \, \, + \, 12 \Lambda_\omega \phi_\omega^3 -6 \Lambda_\omega \phi_\omega g ( \phi_\omega^2 ) -18 \Lambda_\omega \phi_\omega^3 g'(\phi_\omega^2) -78 \Lambda_\omega \phi_\omega^5 g''(\phi_\omega^2) + 8 \Lambda_\omega \phi_\omega^3 g'(\phi_\omega^2)^2 \\
\\ & & \, \,  + \, 56 \Lambda_\omega \phi_\omega^3 g(\phi_\omega^2) g''(\phi_\omega^2) + 28 \Lambda_\omega \phi_\omega g''(\phi_\omega^2) G(\phi_\omega^2) - 56 \Lambda_\omega \phi_\omega^7 g'''(\phi_\omega^2) + 64 \Lambda_\omega \phi_\omega^3 G(\phi_\omega^2) g'''(\phi_\omega^2) \\
\\ & & \, \, + \, 8 \Lambda_\omega \phi_\omega^5 g'(\phi_\omega^2) g''(\phi_\omega^2) + 24 \Lambda_\omega \phi_\omega^5 g ( \phi_\omega^2 ) g'''(\phi_\omega^2 ) -8 \Lambda_\omega \phi_\omega^9 g''''(\phi_\omega^2) + 16 \Lambda_\omega \phi_\omega^5 G(\phi_\omega^2) g''''(\phi_\omega^2) \\
\\ & & \, \, + \, 4 \frac{\Lambda_\omega g(\phi_\omega^2)^2}{\phi_\omega} - 4 \frac{\Lambda_\omega g ( \phi_\omega^2 ) G ( \phi_\omega^2 )}{\phi_\omega^3} + 4 \frac{\Lambda_\omega G(\phi_\omega^2) g'(\phi_\omega^2)}{\phi_\omega}. \end{array} \]
We give without proof several estimates about the operators $X_\alpha$ that can be found in \cite{Ma3} or \cite{Ma1}.

\begin{leftbar}
\noindent \textbf{Lemma 8.} There exists $C>0$ such that, for $\alpha > 0$ small enough and any $q \in L^2 ( \R )$,
\[ \begin{array}{ll} ||X_\alpha q|| \leqslant ||q|| , & || \partial_x X_\alpha^{1/2} q || \leqslant \alpha^{-1/2} ||q||, \\ || \rho X_\alpha q|| \leqslant C ||X_\alpha ( \rho q) || , & ||\rho^{-1} X_\alpha ( \rho q ) || \leqslant C ||X_\alpha q|| , \\ || \eta_A X_\alpha q || \leqslant C || X_\alpha ( \eta_A q ) || \leqslant C || \eta_A q || , & || \eta_A^{-1} X_\alpha ( \eta_A q) || \leqslant C || X_\alpha q || , \\ || \rho^{-1} X_\alpha \partial_x^2 ( \rho q) || \leqslant C \alpha^{-1} ||q|| , & || \rho^{-1} X_\alpha \partial_x ( \rho q ) || \leqslant C \alpha^{-1/2} ||q|| , \\ || \eta_A X_\alpha \partial_x^2 q || \leqslant C \alpha^{-1} || \eta_A q || , & || \eta_A X_\alpha \partial_x q || \leqslant C \alpha^{-1/2} || \eta_A q ||. \end{array} \]
\end{leftbar}

\noindent We then obtain the following estimates, about $M_-$ and $L_+$.

\begin{leftbar}
\noindent \textbf{Lemma 9.} There exists $C>0$ such that, for $\alpha > 0$ small enough and any $q \in L^2 ( \R )$,
\[ \begin{array}{ll} || \eta_A X_\alpha^2 M_- S^2 q || + || \eta_A X_\alpha^2 S^2 L_+ q || \leqslant C \left ( \alpha^{-3/2} || \eta_A \partial_x q || + \omega_0^2 || \eta_A q || \right ), \\ \\ || \eta_A \partial_x X_\alpha^2 M_- S^2 q || + || \eta_A \partial_x X_\alpha^2 S^2 L_+ q || \leqslant C \left ( \alpha^{-2} || \eta_A \partial_x q || + \omega_0^{5/2} || \rho^2 g || \right ). \end{array} \]
\end{leftbar}

\noindent \textit{Proof.} Let us start with $X_\alpha^2 M_- S^2$, whose explicit expression is given before. We have to analyse each term constituting $M_- S^2$. To do so, notice that $X_\alpha$ and $\partial_x$ commute. First,
\[ || \eta_A X_\alpha^2 \partial_x^4 q || = || \eta_A X_\alpha \partial_x^2 (X_\alpha \partial_x \partial_x q) || \leqslant C \alpha^{-1} || \eta_A X_\alpha \partial_x (\partial_x q) || \leqslant C \alpha^{-3/2} || \eta_A \partial_x q ||. \]
We also have
\[ || \eta_A \partial_x X_\alpha^2 \partial_x^4 q || = || \eta_A X_\alpha \partial_x^2 (X_\alpha \partial_x^2 \partial_x q ) || \leqslant C \alpha^{-2} || \eta_A \partial_x q || \]
for the same reason. Now, let $R = \phi_\omega ' / \phi_\omega$ (as in the proof of Lemma 6). We recall for what follows that $| R | \leqslant C \sqrt{\omega}$. Thus,
\[ || \eta_A X_\alpha^2 \partial_x^2 \cdot R \cdot \partial_x q || = || \eta_A X_\alpha \partial_x^2 (X_\alpha \cdot R \cdot \partial_x q) || \leqslant C \alpha^{-1} || \eta_A X_\alpha \cdot R \cdot \partial_x q || \leqslant C \alpha^{-1} || \eta_A R \partial_x q || \leqslant C \alpha^{-1} \sqrt{\omega} || \eta_A \partial_x q ||. \]
And also
\[ || \eta_A \partial_x X_\alpha^2 \partial_x^2 \cdot R \cdot \partial_x q || = || \eta_A X_\alpha \partial_x^2 (X_\alpha \partial_x \cdot R \cdot \partial_x q) || \leqslant C \alpha^{-3/2} \sqrt{\omega} || \eta_A \partial_x q || \]
for the same reason. Then, denoting $b_\omega^1 := 2 \phi_\omega^2 g'( \phi_\omega^2 ) - 4 g ( \phi_\omega^2 ) + 4 \frac{G ( \phi_\omega^2 )}{\phi_\omega^2}$, we find that $| b_\omega^1 | \leqslant C \phi_\omega^2 \leqslant C \omega$. Therefore,
\[ || \eta_A X_\alpha^2 \partial_x \cdot b_\omega^1 \cdot \partial_x q || \leqslant || \eta_A X_\alpha \partial_x ( X_\alpha \cdot b_\omega^1 \cdot \partial_x q) || \leqslant C \alpha^{-1/2} || \eta_A X_\alpha \cdot b_\omega^1 \cdot \partial_x q || \leqslant C \alpha^{-1/2} || \eta_A b_\omega^1 \partial_x q || \leqslant C \alpha^{-1/2} \omega || \eta_A \partial_x q ||. \]
And also
\[ || \eta_A \partial_x X_\alpha^2 \partial_x \cdot b_\omega^1 \cdot \partial_x q || \leqslant || \eta_A X_\alpha \partial_x^2 ( X_\alpha \cdot b_\omega^1 \cdot \partial_x q) || \leqslant C \alpha^{-1} \omega || \eta_A \partial_x q || \]
for the same reason. Now, denoting $b_\omega^2 := 4R \phi_\omega^2 g'(\phi_\omega^2) - 6 R g ( \phi_\omega^2) - 4 R \phi_\omega^4 g''(\phi_\omega^2) + 4 R \frac{G ( \phi_\omega^2)}{\phi_\omega^2} - 2 \omega R$, we see that $|b_\omega^2| \leqslant C |R| \phi_\omega^2 \leqslant C \omega^{3/2}$. Consequently,
\[ || \eta_A X_\alpha^2 b_\omega^2 \cdot \partial_x q || \leqslant C || \eta_A b_\omega^2 \partial_x q || \leqslant C \omega^{3/2} || \eta_A \partial_x q ||. \]
And also
\[ || \eta_A \partial_x X_\alpha^2 b_\omega^2 \cdot \partial_x q || = || \eta_A X_\alpha \partial_x (X_\alpha b_\omega^2 \cdot \partial_x q) || \leqslant C \alpha^{-1/2} \omega^{3/2} || \eta_A \partial_x q || \]
for the same reason. Finally, we denote $b_\omega^3 := \omega^2 + 2 \omega \left ( g ( \phi_\omega^2 ) - \phi_\omega^2 g'(\phi_\omega^2) + 2 \phi_\omega^4 g''( \phi_\omega^2) \right ) - 2 g'(\phi_\omega^2) G ( \phi_\omega^2 ) + \phi_\omega^4 g'(\phi_\omega^2) - 2 \phi_\omega^6 g''( \phi_\omega^2) + 4 \phi_\omega^2 G(\phi_\omega^2) g''(\phi_\omega^2) - 2 \phi_\omega^2 g ( \phi_\omega^2 ) + 2 G ( \phi_\omega^2 ) + g ( \phi_\omega^2 )^2$. We see that $| b_\omega^3 | \leqslant \omega^2 + C \omega \phi_\omega^2 + C \phi_\omega^4 \leqslant C \omega^2$. This gives
\[ || \eta_A X_\alpha^2  ( b_\omega^3 q) || \leqslant C || \eta_A b_\omega^3 q || \leqslant C \omega^2 || \eta_A q ||. \]
On the other hand, $\partial_x (b_\omega^3 q) = (b_\omega^3)' q + b_\omega^3 \partial_x q$ where
\[ \begin{array}{rcl} (b_\omega^3)' &=& 4 \omega \left ( 3 \phi_\omega^3 g''(\phi_\omega^2) + 2 \phi_\omega^5 g'''(\phi_\omega^2) \right ) \phi_\omega ' -10 \phi_\omega ' \phi_\omega^5 g''(\phi_\omega^2) - 4 \phi_\omega ' \phi_\omega^7 g'''(\phi_\omega^2) +4 \phi_\omega ' \phi_\omega G ( \phi_\omega^2 ) g''(\phi_\omega^2) \\ \\ & & \, \, \, \, + \, 8 \phi_\omega ' \phi_\omega^3 g ( \phi_\omega^2 ) g''(\phi_\omega^2) + 8 \phi_\omega ' \phi_\omega^3 G ( \phi_\omega^2 ) g'''(\phi_\omega^2). \end{array} \]
Recalling that $| \phi_\omega ' | \leqslant C \omega \rho^2$, we find that $|(b_\omega^3)'| \leqslant C \omega \phi | \phi ' | + C \phi^3 | \phi'| \leqslant C \omega^{5/2} \rho^2$. This leads to
\[ || \eta_A \partial_x X_\alpha^2 (b_\omega^3 q ) || \leqslant || \eta_A X_\alpha^2 (b_\omega^3)'q|| + || \eta_A X_\alpha^2 (b_\omega^3 \partial_x q)|| \leqslant C || \eta_A (b_\omega^3)' q|| + C \omega^2 || \eta_A \partial_x q || \leqslant C \omega^{5/2} || \rho^2 q || + C \omega^2 || \eta_A \partial_x q ||. \]
We conclude simply by noticing that $\omega \leqslant 1$. The proof for $X_\alpha^2 S^2 L_+ q$ is identical and does not add any complication to the proof above. \hfill \qedsymbol
\\ \textcolor{white}{a} \\ \noindent Applying this lemma to $u_2$ and $u_1$, we obtain the following estimates. 

\begin{leftbar}
\noindent \textbf{Lemma 10.} There exists $C>0$ such that, for $\alpha > 0$ small enough,
\[ \begin{array}{ll} || \eta_A v || \leqslant C \left ( \alpha^{-3/2} || \eta_A \partial_x u || + \omega_0^2 || \eta_A u || \right ), \\ \\ || \eta_A \partial_x v || \leqslant C \left ( \alpha^{-2} || \eta_A \partial_x u || + \omega_0^{5/2} || \rho^2 u || \right ). \end{array} \]
\end{leftbar}

\noindent We have to check similar estimates on the operators $Q_-$ and $Q_+$. 

\begin{leftbar}
\noindent \textbf{Lemma 11.} There exists $C>0$ such that, for $\alpha > 0$ small enough and any $q \in L^2 ( \R )$,
\[ \begin{array}{ll} || \eta_A X_\alpha^2 Q_- q || + || \eta_A X_\alpha^2 Q_+ q || \leqslant C \left ( \alpha^{-1} \sqrt{\omega_0} || \eta_A \partial_x q || + \omega_0^2 || \eta_A q || \right ), \\ \\ || \eta_A \partial_x X_\alpha^2 Q_- q || + || \eta_A \partial_x X_\alpha^2 Q_+ q || \leqslant C \left ( \alpha^{-3/2} \sqrt{\omega_0} || \eta_A \partial_x q || + \omega_0^{5/2} || \rho^2 g || \right ). \end{array} \]
\end{leftbar}

\noindent \textit{Proof.} The proof is similar to the one of the previous lemma. We first show that
\[ \left | \frac{\Lambda_\omega ' \phi_\omega - \Lambda_\omega \phi_\omega '}{\phi_\omega^2} \right | \leqslant C \sqrt{\omega}. \]
Indeed, we first see that
\[ ( \Lambda_\omega ' \phi_\omega - \Lambda_\omega \phi_\omega ' )' = \Lambda_\omega '' \phi_\omega - \Lambda_\omega \phi_\omega '' = \omega \phi_\omega^2 -2 \Lambda_\omega \phi_\omega^3 + 2 \Lambda_\omega \phi_\omega^3 g'(\phi_\omega^2), \]
using the equations satisfied by $\phi_\omega$ and $\Lambda_\omega$. Therefore, writing that $|g'(\phi_\omega^2)| \leqslant 1$, we see that, for any $x \geqslant 0$,
\[ | \Lambda_\omega ' \phi_\omega - \Lambda_\omega \phi_\omega ' | (x) = \left | - \int_x^{+ \infty} \left ( \omega \phi_\omega^2 - 2 \Lambda_\omega \phi_\omega^3 + 2 \Lambda_\omega \phi_\omega^3 g'(\phi_\omega^2) \right ) \right | \leqslant C \omega \int_x^{+ \infty} \phi_\omega^2 + C \int_x^{+ \infty} | \Lambda_\omega | \phi_\omega^3. \]
Now using the estimates on $\Lambda_\omega$ and $\phi_\omega$ we get
\[ | \Lambda_\omega ' \phi_\omega - \Lambda_\omega \phi_\omega ' | (x) \leqslant C \omega^{3/2} e^{-2 \sqrt{\omega} x} + C \omega^{3/2} e^{-4 \sqrt{\omega} x} \leqslant C \omega^{3/2} e^{-2 \sqrt{\omega} x}. \]
We recall that $\phi_\omega (x) \geqslant c \sqrt{\omega} e^{- \sqrt{\omega} |x|}$. Thus,
\[ \left | \frac{\Lambda_\omega ' \phi_\omega - \Lambda_\omega \phi_\omega '}{\phi_\omega^2} \right | \leqslant C \sqrt{\omega}. \]
We also see, thanks to the estimates on $\Lambda_\omega$ and its derivatives, that $| \Lambda_\omega | \leqslant C \sqrt{\omega}$ and $| \Lambda_\omega ' | \leqslant C \omega$. Now let us write the operator $Q_-$ as
\[ Q_- = \partial_x^2 \cdot c_\omega^1 \cdot \partial_x + \partial_x \cdot c_\omega^2 \cdot \partial_x + c_\omega^3 \cdot \partial_x + c_\omega^4. \]
Using $(H_1)$, we see that $|c_\omega^1| \leqslant C \sqrt{\omega}$, $|c_\omega^2| \leqslant C \omega$, $|c_\omega^3| \leqslant C \omega^{3/2}$, $|c_\omega^4| \leqslant C \omega^2$ and $|(c_\omega^4)'| \leqslant C \omega^{5/2}$. Reasoning as in the previous proof, we obtain the desired result. The same estimates and the same proof hold for $Q_+$. It is for this proof that we use $(H_1)$ in its entirety: we indeed have to control $g$ up to its fifth derivative (because of the expression of $Q_+$). \hfill \qedsymbol

\noindent \textcolor{white}{a} \\ \noindent Now let us prove a last estimate, more elementary (in the sense that it does not involve any derivative of $q$) but that will be useful.

\begin{leftbar}
\noindent \textbf{Lemma 12.} There exists $C>0$ such that, for $\alpha > 0$ small enough and any $q \in L^2 ( \R )$,
\[ || \eta_A X_\alpha^2 M_- S^2 q || + || \eta_A X_\alpha^2 S^2 L_+ q || \leqslant C \alpha^{-2} || \eta_A q ||. \]
\end{leftbar}

\noindent \textit{Proof.} The proof is analogous to the one of Lemma 9. For example, see that
\[ || \eta_A X_\alpha^2 \partial_x^4 q || = || \eta_A X_\alpha \partial_x^2 ( X_\alpha \partial_x^2 q) || \leqslant C \alpha^{-1} || \eta_A X_\alpha \partial_x^2 q || \leqslant C \alpha^{-2} || \eta_A \partial_x q ||. \]
For the other terms it is similar and easier; for instance the last term is controlled as follows:
\[ || \eta_A \partial_x X_\alpha^2 (b_\omega^3 q) || \leqslant C \alpha^{-1/2} || \eta_A X_\alpha ( b_\omega^3 q) || \leqslant C \alpha^{-1/2} || \eta_A b_\omega^3 q || \leqslant C \alpha^{-1/2} || \eta_A q ||. \]
This completes the proof. \hfill \qedsymbol

\subsection{Second virial estimate}
\noindent Using the system \eqref{Su} satisfied by $u$ and the identity of Lemma 6, we find the following system satisfied by $v$:
\begin{equation}
    \left \{ \begin{array}{ccl} \partial_t v_1 &=& M_- v_2 + Y_\alpha^- v_2 + X_\alpha^2 n_2 - X_\alpha^2 r_2 \\ \partial_t v_2 &=& - M_+ v_1 - Y_\alpha^+ v_1 - X_\alpha^2 n_1 + X_\alpha^2 r_1 \end{array} \right.
    \label{Sv}
\end{equation}
where
\[ \begin{array}{lll} n_1 \, = \, S^2 L_+ m_2 + \frac{\dot{\omega}}{\omega} Q_+ u_1,  & r_1 \, = \, S^2 L_+ q_2 , & Y_\alpha^- \, = \, X_\alpha^2 \cdot a_\omega^- \cdot X_\alpha^{-2} - a_\omega^-, \\
\\ n_2 \, = \, -M_- S^2 m_1 + \frac{\dot{\omega}}{\omega} Q_- u_2 , & r_2 \, = \, -M_- S^2 q_1 , & Y_\alpha^+ \, = \, X_\alpha^2 \cdot a_\omega^+ \cdot X_\alpha^{-2} - a_\omega^+. \end{array} \]

\begin{leftbar}
\noindent \textbf{Proposition 4.} Suppose hypotheses $(H_1)$ and $(H_2)$ are satisfied. Assume that $\omega_0$ is small enough. There exists $C>0$ such that, for $B>0$ large enough, $\alpha > 0$ and $\epsilon > 0$ small enough, and for any $T>0$,
\[ \varepsilon_{\omega_0}^2 \omega_0 \int_0^T || \rho v ||^2 \, \text{d}t \leqslant C \epsilon^2 + C \int_0^T \left ( \frac{1}{A \sqrt{\omega_0}} || \eta_A \partial_x u ||^2 + \frac{\omega_0^{5/2}}{A^3} || \eta_A u ||^2 + \frac{\omega_0^5}{A} || \rho^2 u ||^2 \right ) \, \text{d}t. \]
\end{leftbar}

\noindent \textit{Proof.} We use another virial argument. Let $z = \chi_A \zeta_B v$ and
\[ \mathcal{J} = \int_{\R} v_1 \left ( 2 \Psi_{A,B} \partial_x v_2 + \Psi_{A,B} ' v_2 \right ). \]
Using the equation \eqref{Sv} and integrating by parts (following computations from \cite{Ma3} and \cite{Ma1}), we get that
\[ \dot{\mathcal{J}} = \int_{\R} \left ( 2 ( \partial_x z_1)^2 + P_B^+ z_1^2 \right ) + \int_{\R} \left ( 2 ( \partial_x z_2 )^2 + P_B^- z_2^2 \right ) + J_1 + J_2 + J_3 + J_4 + J_5 \]
where $\displaystyle{P_B^{\pm} := - (a_{\omega_0}^{\pm})' \frac{\Phi_B}{\zeta_B^2}}$ and
\[ \begin{array}{l} \displaystyle{J_1 \, = \, \sum_{k=1}^2 \int_{\R} (\ln \zeta_B)'' z_k^2} , \\
\\ \displaystyle{J_2 \, = \, - \sum_{k=1}^2 \int_{\R} \left ( \frac{1}{2} ( \chi_A^2 )' (\zeta_B^2) ' + \left ( 3 ( \chi_A ')^2 + \chi_A '' \chi_A \right ) \zeta_B^2 + \frac{1}{2} ( \chi_A^2 )''' \Phi_B \right ) v_k^2 + 2 \sum_{k=1}^2 \int_{\R} (\chi_A^2)' \Phi_B ( \partial_x v_k)^2} , \\
\\ \displaystyle{J_3 \, = \, \int_{\R} ( 2 \Psi_{A,B} \partial_x v_1 + \Psi_{A,B} ' v_1 ) Y_\alpha^+ v_1 + \int_{\R} (2 \Psi_{A,B} \partial_x v_2 + \Psi_{A,B} ' v_2 ) Y_\alpha^- v_2} , \\
\\ \displaystyle{J_4 \, = \, \sum_{k=1}^2 \int_{\R} ( 2 \Psi_{A,B} \partial_x v_k + \Psi_{A,B} ' v_k ) ( X_\alpha^2 n_k - X_\alpha^2 r_k )} , \\
\\ \displaystyle{J_5 \, = \, \int_{\R} \frac{\Phi_B}{\zeta_B^2} \left ( \left ( a_{\omega_0}^- - a_{\omega}^- \right )' z_1^2 + \left ( a_{\omega_0}^+ - a_{\omega}^+ \right )' z_2^2 \right ).} \end{array} \]
Notice the obvious similarities with the notation in Lemma 7 and Proposition 2; however, the pulsation involved in $P_B$ is $\omega_0$ (not $\omega$). Setting $\mathcal{K} := - \int_{\R} z_1 z_2 R_B$ where $R_B$ is a bounded function to be defined later, we find that
\[ \dot{\mathcal{J}} + \dot{\mathcal{K}} = \int_{\R} \left [ 2 ( \partial_x z_1)^2 + \left ( P_B^+ + \omega_0 R_B - \frac{R_B ''}{2} \right ) z_1^2 \right ] + \int_{\R} \left [ 2 ( \partial_x z_2)^2 + \left ( P_B^- - \omega_0 R_B + \frac{R_B''}{2} \right ) z_2^2 \right ] + \sum_{j=1}^5 (J_j + K_j ) \]
where
\[ \begin{array}{l} \displaystyle{K_1 \, = \, \sum_{k=1}^2 (-1)^k \int_{\R} \left ( (\chi_A \zeta_B)' \chi_A \zeta_B R_B' + ((\chi_A \zeta_B)')^2 R_B \right ) v_k^2} , \\
\\ \displaystyle{K_2 \, = \, \int_{\R} \left ( ( \partial_x z_1)^2 - ( \partial_x z_2)^2 \right ) R_B - \int_{\R} (a_{\omega_0}^- z_2^2 - a_{\omega_0}^+ z_1^2) R_B} , \\
\\ \displaystyle{K_3 \, = \, \int_{\R} \left ( (Y_\alpha^+ v_1)v_1 - (Y_\alpha^- v_2)v_2 \right ) \chi_A^2 \zeta_B^2 R_B} , \\
\\ \displaystyle{K_4 \, = \, \sum_{k=1}^2 (-1)^{k-1} \int_{\R} (X_\alpha^2 n_k - X_\alpha^2 r_k) v_k \chi_A^2 \zeta_B^2 R_B} , \\
\\ \displaystyle{K_5 \, = \, ( \omega - \omega_0 ) \int_{\R} (z_1^2 - z_2^2) R_B + \int_{\R} \left [ (a_\omega^+ - a_{\omega_0}^+ ) z_1^2 - (a_\omega^- - a_{\omega_0}^-) z_2^2 \right ] R_B.} \end{array} \]
Let us define $R_B$ as the bounded solution of the ordinary differential equation $- \frac{R_B ''}{2} + \omega_0 R_B = D_B$ where $D_B := \frac{P_B^+ - P_B^-}{2}$. Here also, notice the similarities with the notation in Lemma 7 and Proposition 2. We have the control $|R_B| \leqslant C \varepsilon_{\omega_0} \rho$. Such a choice leads to
\[ \dot{\mathcal{J}} + \dot{\mathcal{K}} = \int_{\R} \left [ 2 ( \partial_x z_1)^2 + P_B z_1^2 \right ] + \int_{\R} \left [ 2 ( \partial_x z_2)^2 + P_B z_2^2 \right ] + \sum_{j=1}^5 (J_j + K_j ).  \]
We will need a result that enables us to control $|| \rho \partial_x v||$ and $|| \rho v ||$ in terms of $||\partial_x z||$ and $||\rho z||$, plus error terms involving $u$. This is the following lemma.

\begin{leftbar}
\noindent \textbf{Lemma 13.} There exists $C>0$ such that, for $A,B>0$ large enough (depending on $\omega_0$) and $\alpha > 0$ small enough,
\begin{align*} 
& || \rho v ||^2 \leqslant C \int_{\R} \rho |z|^2 + \frac{C}{A^3 \omega_0^{3/2}} \left ( \alpha^{-4} || \eta_A \partial_x u ||^2 + \omega_0^4 || \eta_A u ||^2 \right ) \\
\text{and} \quad & || \rho \partial_x v ||^2 + \leqslant C \int_{\R} \left ( | \partial_x z |^2 + \frac{1}{B^2} \rho |z|^2 \right ) + \frac{C}{A^3 \omega_0^{3/2}} \left ( \alpha^{-4} || \eta_A \partial_x u ||^2 + \omega_0^4 || \eta_A u ||^2 \right ).
\end{align*}
\end{leftbar}

\noindent \textit{Proof.} First, for $|x| \leqslant A$, $z = \zeta_B v$ and we write that
\[ \int_{|x| \leqslant A} \rho^2 |v|^2 \leqslant C \int_{|x| \leqslant A} \rho \zeta_B^2 |v|^2 = C \int_{|x| \leqslant A} \rho |z|^2 \]
using that $\rho \leqslant C \zeta_B^2$. Now, we have $\partial_x z = \zeta_B ' v + \zeta_B \partial_x v$ and $| \zeta_B ' | \leqslant \frac{C}{B} \zeta_B$ which lead to
\[ \rho^2 | \partial_x v |^2 \leqslant C \rho \zeta_B^2 | \partial_x v |^2 \leqslant C \rho | \partial_x z|^2 + C \rho \frac{\zeta_B^2}{B^2} |v|^2 \leqslant C | \partial_x z |^2 + \frac{C}{B^2} \, \rho |z|^2. \]
Therefore,
\[ \int_{|x| \leqslant A} \rho^2 | \partial_x v |^2 \leqslant C \int_{|x| \leqslant A} | \partial_x z |^2 + \frac{C}{B^2} \int_{|x| \leqslant A} \rho | z|^2 . \]
Now, for $|x| > A$, we see that $\rho (x)^2 \leqslant C e^{\left ( \frac{4}{A} - \frac{\sqrt{\omega_0}}{5} \right ) |x|} \eta_A (x)^2$. If we take $A$ large enough such that $\frac{4}{A} < \frac{\sqrt{\omega_0}}{5}$, we see that 
\[ \rho^2 \leqslant C e^{- \frac{A \sqrt{\omega_0}}{5}} \eta_A^2 \leqslant \frac{C}{A^3 \omega_0^{3/2}} \eta_A^2, \]
the last inequality being true if $A \sqrt{\omega_0}$ is large enough, i.e. if $A$ is large enough (depending on $\omega_0$). Then, using Lemma 10, we obtain
\[  \begin{array}{rcl} \displaystyle{\int_{|x| > A} \rho^2 ( | \partial_x v |^2 + |v|^2 )} & \leqslant & \displaystyle{\frac{C}{A^3 \omega_0^{3/2}} \left ( || \eta_A \partial_x v ||^2 + || \eta_A v ||^2 \right )} \\
\\ & \leqslant & \displaystyle{\frac{C}{A^3 \omega_0^{3/2}} \left ( \alpha^{-3} || \eta_A \partial_x u ||^2 + \omega_0^4 || \eta_A u ||^2 + \alpha^{-4} || \eta_A \partial_x u ||^2 + \omega_0^5 || \rho^2 u ||^2 \right )} \\
\\ & \leqslant & \displaystyle{\frac{C}{A^3 \omega_0^{3/2}} \left ( \alpha^{-3} || \eta_A \partial_x u ||^2 + \omega_0^4 || \eta_A u ||^2 + \alpha^{-4} || \eta_A \partial_x u ||^2 \right ).} \end{array} \]
Putting these estimates together, we get the desired result. \hfill \qedsymbol 

\textcolor{white}{a} \\ \noindent We now get back to the proof of Proposition 4 and in the first place we control the terms $J_j$,$K_j$. \\
\\ \textit{(About $J_1$.)} We write that
\[ | ( \ln \zeta_B )''| \leqslant \frac{C \sqrt{\omega_0}}{B} \, \mathbbm{1}_{[1,2]} ( \sqrt{\omega_0} |x| ) \leqslant \frac{C \sqrt{\omega_0}}{B} \rho, \]
which leads to
\[ |J_1| \leqslant \frac{C \sqrt{\omega_0}}{B} \int_{\R} \rho |z|^2. \]
\textit{(About $K_1$.)} We start by writing that $| \chi_A ' | \leqslant \frac{C}{A} \leqslant \frac{C}{B}$, $| \zeta_B ' | \leqslant \frac{C}{B} \zeta_B$, $|R_B| \leqslant C \varepsilon_{\omega_0} \rho^2$ and $|R_B'| \leqslant C \varepsilon_{\omega_0} \sqrt{\omega_0} \, \rho^2$. The estimates on $R_B$ are shown similarly as the estimates on $R_\infty$ in the proof of Proposition 2. Recalling that $B > \omega_0^{-1/2}$, this leads to
\[ \left | ( \chi_A \zeta_B )' \chi_A \zeta_B R_B ' + ( ( \chi_A \zeta_B )')^2 R_B \right | \leqslant \frac{C \varepsilon_{\omega_0} \sqrt{\omega_0}}{B} \, \rho^2 \]
and then 
\[ \begin{array}{rcl} | K_1 | & \leqslant & \displaystyle{\frac{C \varepsilon_{\omega_0} \sqrt{\omega_0}}{B} \int_{\R} \rho^2 |v|^2} \\
\\ & \leqslant & \displaystyle{ \frac{C \varepsilon_{\omega_0} \sqrt{\omega_0}}{B}  \left [ || \partial_x z ||^2 + \int_{\R} \rho |z|^2 + \frac{1}{A^3 \omega_0^{3/2}} \left ( \alpha^{-4} || \eta_A \partial_x u ||^2 + \omega_0^4 || \eta_A u ||^2 \right ) \right ],} \end{array} \]
using Lemma 13. \\
\\ \textit{(About $J_2$.)} We start by recalling that $| \chi_A ' | \leqslant \frac{C}{A} \mathbbm{1}_{A < |x| < 2A}$, $| \chi_A '' | \leqslant \frac{C}{A^2} \mathbbm{1}_{A < |x| < 2A}$ and $| \chi_A ''' | \leqslant \frac{C}{A^3}\mathbbm{1}_{A < |x| < 2A}$. Moreover, for $|x|>A$, $| \zeta_B (x)| \leqslant C e^{-A/B}$ and $| \zeta_B '(x)| \leqslant \frac{C}{B} e^{-A/B}$. Thus, using the fact that $\zeta_B \leqslant C \eta_A^2$ (since $A \gg B$),
\[ \begin{array}{l} \displaystyle{| ( \chi_A^2 ) ' ( \zeta_B^2 )' | \leqslant \frac{C e^{-A/B}}{AB} \zeta_B^2 \leqslant \frac{C e^{-A/B}}{AB} \eta_A^2 \leqslant \frac{CB}{A^{3}} \eta_A^2 \, , } \\ \\ \displaystyle{ \left ( ( \chi_A')^2 + | \chi_A '' \chi_A | \right ) \zeta_B^2 \leqslant \frac{C e^{-A/B}}{A^2} \zeta_B^2 \leqslant \frac{C e^{-A/B}}{A^2} \eta_A^2 \leqslant \frac{CB}{A^3} \eta_A^2 \, ,} \end{array} \]
for $A/B$ large enough (we recall that $A \gg B$). We also know that $| \Phi_B | \leqslant CB$. Using the fact that $\mathbbm{1}_{|x| < 2A} \leqslant C \eta_A^2$, we obtain
\[ | ( \chi_A^2)' \Phi_B | \leqslant \frac{CB}{A} \, \eta_A^2 \, , \, \, \, \, \, \, | (\chi_A^2)''' \Phi_B | \leqslant \frac{CB}{A^3} \, \eta_A^2. \]
Putting these estimates together we get
\[ \begin{array}{rcl} |J_2| & \leqslant & \displaystyle{\frac{CB}{A} || \eta_A \partial_x v||^2 + \frac{CB}{A^3} || \eta_A v ||^2} \\
\\ & \leqslant & \displaystyle{\frac{CB}{A} \left ( \alpha^{-4} || \eta_A \partial_x u ||^2 + \omega^{5} || \rho^2 u ||^2 \right ) + \frac{CB}{A^3} \left ( \alpha^{-3} || \eta_A \partial_x u ||^2 + \omega^4 || \eta_A u ||^2 \right )} \\
\\ & \leqslant & \displaystyle{\frac{CB \alpha^{-4}}{A} || \eta_A \partial_x u ||^2 + \frac{CB \omega_0^4}{A} \left ( \frac{1}{A^2} || \eta_A u ||^2 + \omega_0 || \rho^2 u ||^2 \right ). } \end{array} \]
\textit{(About $K_2$.)} We know that $R_B$ is bounded and that $||R_B||_\infty \leqslant C \varepsilon_{\omega_0}$. Moreover, $|a_{\omega_0}^{\pm}| \leqslant C \varepsilon_{\omega_0} \phi_{\omega_0}^2 \leqslant C \varepsilon_{\omega_0} \omega_0 \rho$. This gives
\[ |K_2| \leqslant C \varepsilon_{\omega_0} || \partial_x z ||^2 + C_2 \varepsilon_{\omega_0}^2 \omega_0 \int_{\R} \rho |z|^2. \]
Here, an explicit name has been given to the constant $C_2$ in order to be clear a little later. \\
\\ \textit{(About $J_3$.)} We have $| \Psi_{A,B} | \leqslant CB$ and $| \Psi_{A,B} ' | \leqslant C$ (thanks to the bounds $| \chi_A ' | \leqslant C/B$ and $| \Phi_B | \leqslant CB$). Using the Cauchy-Schwarz inequality, we find
\[ \begin{array}{rcl} \displaystyle{\left | \int_{\R} ( 2 \Psi_{A,B} \partial_x v_1 + \Psi_{A,B} ' v_1 ) Y_\alpha^+ v_1 \right |} &=& \displaystyle{\left | \int_{\R} ( 2 \Psi_{A,B} \partial_x v_1 + \Psi_{A,B} ' v_1 ) \rho \cdot \rho^{-1} Y_\alpha^+ v_1 \right |} \\
\\ & \leqslant & \displaystyle{\left \| (2 \Psi_{A,B} \partial_x v_1 + \Psi_{A,B} ' v_1 ) \rho \right \| + || \rho^{-1} Y_\alpha^+ v_1 ||} \\
\\ & \leqslant & C \left ( B || \rho \partial_x v_1 || + || \rho v_1 || \right ) || \rho^{-1} Y_\alpha^+ v_1 || \end{array} \]
where we recall that
\[ \begin{array}{rcl} Y_\alpha^{\pm} &=& X_\alpha^2 (a_\omega^{\pm} \cdot X_\alpha^{-2} - X_\alpha^{-2} \cdot a_\omega^{\pm} ) \\
\\ &=& X_\alpha^2 \cdot \left [ 2 \alpha ( 2 \partial_x \cdot (a_\omega^\pm) ' - ( a_\omega^{\pm} ) '' ) + \alpha^2 \left (  -4 \partial_x^3 \cdot (a_\omega^{\pm} )' +6 \partial_x^2 \cdot (a_\omega^{\pm} )'' -4 \partial_x \cdot (a_\omega^{\pm}) ''' + (a_\omega^{\pm})'''' \right ) \right ]. \end{array} \]
Using Lemma 8 and the bounds on $a_\omega^{\pm}$ and its derivative, we find
\[ \begin{array}{rcl} || \alpha \rho^{-1} X_\alpha^2 \partial_x ( (a_\omega^\pm)'v_k) || & = & \alpha || \rho^{-1} X_\alpha \partial_x ( \rho \rho^{-1} X_\alpha ((a_\omega^{\pm})'v_k)) || \, \, \leqslant \, \, \alpha \cdot C \alpha^{-1/2} || \rho^{-1} X_\alpha ((a_\omega^\pm)'v_k)|| \\
\\ & \leqslant & C \sqrt{\alpha} || \rho^{-1} (a_\omega^\pm)' v_k || \, \, \leqslant \, \, C \sqrt{\alpha} \, \omega^{3/2} || \rho v_k ||. \end{array} \]
Similarly, we find for instance $|| \alpha^2 \rho^{-1} X_\alpha^2 \partial_x^3 ( (a_\omega^\pm)' v_k ) || \leqslant \alpha^2 \cdot C \alpha^{-3/2} || \rho^{-1} (a_\omega^\pm)' v_k || \leqslant C \sqrt{\alpha} \, \omega^{3/2} || \rho v_k ||$. All the other terms are smaller, for example $|| \alpha \rho^{-1} X_\alpha^2 ( (a_\omega^\pm)'' v_k) || \leqslant C \alpha \omega^2 || \rho v_k ||$. We obtain the following estimate:
\[ || \rho^{-1} Y_\alpha^+ v_1 || \leqslant C \sqrt{\alpha} \, \omega^{3/2} || \rho v_1 || \]
and a similar estimate holds for $Y_\alpha^- v_2$. These estimates lead to
\[  \begin{array}{rcl} |J_3| & \leqslant & \displaystyle{C \sum_{k=1}^2 \left ( B || \rho \partial_x v_k || + || \rho v_k || \right ) \sqrt{\alpha} \, \omega^{3/2}  || \rho v_k ||} \\
\\ & \leqslant & \displaystyle{C \sqrt{\alpha} \,  \omega^{3/2}  \sum_{k=1}^2 \left ( B^2 || \rho \partial_x v_k ||^2 + || \rho v_k ||^2 \right ) } \\
\\ & \leqslant & \displaystyle{C \sqrt{\alpha} \, \omega^{3/2} \left ( B^2 || \rho \partial_x v ||^2 + || \rho v ||^2 \right )} \\
\\ & \leqslant & \displaystyle{C \sqrt{\alpha} \, \omega_0^{3/2} B^2 \left [ || \partial_x z ||^2 +  \int_{\R} \rho |z|^2 + \frac{1}{A^3 \omega_0^{3/2}} \left ( \alpha^{-4} || \eta_A \partial_x u ||^2 + \omega_0^4 || \eta_A u ||^2 \right ) \right ].} \end{array} \]
\textit{(About $K_3$.)} The estimate is quite similar to $J_3$. We use the bounds $\chi_A^2 \zeta_B^2 \leqslant 1$ and $|R_B| \leqslant C \varepsilon_{\omega_0} \leqslant C$, as well as the Cauchy-Schwarz inequality again:
\[ \begin{array}{rcl} \displaystyle{\left | \int_{\R} ( Y_\alpha^+ v_1) v_1 \chi_A^2 \zeta_B^2 R_B \right |} & \leqslant & \displaystyle{C \int_{\R} | \rho^{-1} Y_\alpha^+ v_1 | \, | \rho v_1 |} \\
\\ & \leqslant & \displaystyle{C || \rho^{-1} Y_\alpha^+ v_1 || \, || \rho v_1 ||} \\
\\ & \leqslant & \displaystyle{C \sqrt{\alpha} \, \omega^{3/2} \varepsilon_\omega || \rho v_1 ||^2} \\
\\ & \leqslant & \displaystyle{C \sqrt{\alpha} \, \omega_0^{3/2} \varepsilon_{3 \omega_0 /2} \left [ || \partial_x z ||^2 + \int_{\R} \rho |z|^2 + \frac{1}{A^3 \omega_0^{3/2}} \left ( \alpha^{-4} || \eta_A \partial_x u ||^2 + \omega_0^4 || \eta_A u ||^2 \right ) \right ],} \end{array} \]
using the estimates obtained previously. \\
\\ \textit{(About $J_4$.)} First, we recall from the proof of Proposition 3 that 
\[ \left | \text{Re} \left [ h(\phi_\omega + u) - h ( \phi_\omega ) - h'(\phi_\omega) u \right ] \right | \leqslant C \left ( \phi_\omega |u|^2 + |u|^3 \right ). \]
This shows that $|q_1| \leqslant C |u|^2 \leqslant C \epsilon |u|$. Now, to control $q_2$, let us write
\[ \text{Im} \left [ h(\phi_\omega + u) - \frac{h ( \phi_\omega )}{\phi_\omega} \, u \right ] = |u|^2 u_2 + 2 \phi_\omega u_1 u_2 - u_2 \left ( g(\phi_\omega^2 + |u|^2 + 2 \phi_\omega u_1) - g(\phi_\omega^2) \right ). \]
Here we notice that
\[ | g(\phi_\omega^2 + |u|^2 + 2 \phi_\omega u_1 ) - g(\phi_\omega^2) | = \left | \int_{\phi_\omega^2}^{\phi_\omega^2 + |u|^2 + 2 \phi_\omega u_1} g'(s) \, \text{d}s \right | \leqslant \left | |u|^2 + 2 \phi_\omega u_1 \right | \leqslant C |u| \]
which gives $|q_2| \leqslant C |u|^2 \leqslant C \epsilon |u|$. Using the definitions of $r_1$ and $r_2$, we find that, for $k \in \{ 1 \, , 2 \}$,
\[ || \eta_A X_\alpha^2 r_k || \leqslant C \alpha^{-2} || \eta_A q_k || \leqslant C \alpha^{-2} \epsilon || \eta_A u ||. \]
Hence, using the Cauchy-Schwarz inequality and the upper bounds $|\Psi_{A,B}| \leqslant CB \eta_A^2$ and $|\Psi_{A,B}'| \leqslant C \eta_A^2$,
\[ \begin{array}{rcl} \displaystyle{\left | \int_{\R} ( 2 \Psi_{A,B} \partial_x v_k + \Psi_{A,B} ' v_k ) X_\alpha^2 r_k \right |} & \leqslant & \displaystyle{C \left ( B || \eta_A \partial_x v_k || \, || \eta_A X_\alpha^2 r_k || + || \eta_A v_k || \, || \eta_A X_\alpha^2 r_k || \right )} \\
\\ & \leqslant & C \alpha^{-2} \epsilon || \eta_A u || \left [ B \left ( \alpha^{-2} || \eta_A \partial_x u || + \omega_0^{5/2} || \rho^2 u || \right ) + \alpha^{-3/2} || \eta_A \partial_x u || + \omega_0^2 || \eta_A u || \right ] \\
\\ & \leqslant & C \alpha^{-2} \epsilon || \eta_A u || B \left ( \alpha^{-2} || \eta_A \partial_x u || + \omega_0^{5/2} || \eta_A u || \right )  \\
\\ & \leqslant & C \alpha^{-2} B \epsilon \left ( || \eta_A u ||^2 + \alpha^{-4} || \eta_A \partial_x u ||^2 \right ) \end{array} \]
where we have used that $|| \rho^2 u || \leqslant || \eta_A u ||$ and $B > \omega_0^{-1/2}$. Now, let us control the other term in $J_4$. We write that
\[ || \eta_A X_\alpha^2 n_k || \leqslant C \alpha^{-2} || \eta_A m_k || + \left | \frac{\dot{\omega}}{\omega} \right | \, || \eta_A X_\alpha^2 Q_{\pm} u_k ||. \]
Gathering the estimates $| x \eta_A | \leqslant CA$ and \eqref{orth}, we see that
\[ \begin{array}{rcl} || \eta_A m_k || & = & || \dot{\beta} x \eta_A u_k + ( \dot{\gamma} - \omega - \beta^2 ) \eta_A u_k \pm ( \dot{\sigma} - 2 \beta ) \eta_A \partial_x u_{3-k} - \beta ( \dot{\sigma} - 2 \beta ) \eta_A u_k || \\
\\ & \leqslant & C \left ( \omega_0 A ||u_k|| + \sqrt{\omega} || \eta_A u_k || + || \eta_A \partial_x u_{3-k} || + || \eta_A u_k || \right ) \underbrace{|| \rho^2 u ||^2}_{\leqslant \, C \epsilon || \eta_A u ||} \\
& \leqslant & C \sqrt{\omega_0} A \epsilon^2 || \eta_A u ||, \end{array} \]
using \eqref{orbstab} and the fact that $\sqrt{\omega_0} A > 1$. Besides,
\[ || \eta_A X_\alpha^2 Q_{\pm} u_k || \leqslant C \left ( \alpha^{-1} \sqrt{\omega_0} || \eta_A \partial_x u_k|| + \omega_0^2 || \eta_A u_k || \right ) \]
which leads to
\[ \begin{array}{rcl} || \eta_A X_\alpha^2 n_k || & \leqslant & C \alpha^{-2} \sqrt{\omega_0} A \epsilon^2 || \eta_A u || + C |\dot{\omega}| \left ( \alpha^{-1} \omega_0^{-1/2} || \eta_A \partial_x u || + \omega_0 || \eta_A u || \right ) \\
\\ & \leqslant & C \alpha^{-2} \sqrt{\omega_0} A \epsilon^2 || \eta_A u || + C \omega_0^{3/2} \epsilon^2 \left ( \alpha^{-1} \omega_0^{-1/2} || \eta_A \partial_x u || + \omega_0 || \eta_A u|| \right ) \\
\\ & \leqslant & C \alpha^{-2} \sqrt{\omega_0} A \epsilon^2 || \eta_A u || + C \alpha^{-1} \omega_0 \epsilon^2 || \eta_A \partial_x u ||. \end{array} \]
Hence, using the same arguments as previously,
\[ \begin{array}{rcl} \displaystyle{\left | \int_{\R} (2 \Psi_{A,B} \partial_x v_k + \Psi_{A,B} ' v_k ) X_\alpha^2 n_k \right |} & \leqslant & C \left ( B || \eta_A \partial_x v_k || + || \eta_A v_k || \right ) || \eta_A X_\alpha^2 n_k || \\
\\ & \leqslant & C \left ( B \alpha^{-2} || \eta_A \partial_x u || + B \omega_0^{5/2} || \eta_A u || \right ) \left ( \alpha^{-2} \sqrt{\omega_0} A \epsilon^2 || \eta_A u || + \alpha^{-1} \omega_0 \epsilon^2 || \eta_A \partial_x u || \right ) \\
\\ & \leqslant & C (AB \sqrt{\omega_0} \alpha^{-4} \epsilon^2 + B \alpha^{-3} \omega_0 \epsilon^2 ) || \eta_A \partial_x u ||^2 + C(AB \alpha^{-2} \omega_0^{3} \epsilon^2 + AB \alpha^{-4} \sqrt{\omega_0} \epsilon^2) || \eta_A u ||^2 \\
\\ & \leqslant & CAB \sqrt{\omega_0} \alpha^{-4} \epsilon^2 \left ( || \eta_A \partial_x u||^2 + || \eta_A u ||^2 \right ) \end{array} \]
after computations. Gathering these estimates we find
\[ | J_4 | \leqslant C (AB \sqrt{\omega_0} \alpha^{-4} \epsilon^2 + \alpha^{-2} B \epsilon ) || \eta_A u ||^2 + C (AB \sqrt{\omega_0} \alpha^{-4} \epsilon^2 + \alpha^{-6} B \epsilon ) || \eta_A \partial_x u ||^2. \]
\textit{(About $K_4$.)} The estimates we use are the same as for $J_4$ and the integral upper bounds are slightly easier. We recall that $|\chi_A^2 \zeta_B^2 R_B| \leqslant C \varepsilon_{\omega_0} \eta_A^2$. We find
\[ \begin{array}{rcl} \displaystyle{\left | \int_{\R} (X_\alpha^2 n_k) v_k \chi_A^2 \zeta_B^2 R_B \right |} & \leqslant & \displaystyle{C \varepsilon_{\omega_0} || \eta_A v_k || \, || \eta_A X_\alpha^2 n_k ||} \\
\\ & \leqslant & \displaystyle{C \varepsilon_{\omega_0} \left ( \alpha^{-3/2} || \eta_A \partial_x u || + \omega_0^2 || \eta_A u || \right ) \left ( \alpha^{-2} \sqrt{\omega_0} A \epsilon^2 || \eta_A u || + \alpha^{-1} \omega_0 \epsilon^2 || \eta_A \partial_x u || \right )} \\
\\ & \leqslant & \displaystyle{C \varepsilon_{\omega_0} \alpha^{-7/2} \sqrt{\omega_0} A \epsilon^2 \left ( || \eta_A u ||^2 + || \eta_A \partial_x u ||^2 \right )} \end{array} \]
after computations. And on the other hand,
\[ \begin{array}{rcl} \displaystyle{\left | \int_{\R} (X_\alpha^2 r_k) v_k \chi_A^2 \zeta_B^2 R_B \right |} & \leqslant & C \varepsilon_{\omega_0} || \eta_A v_k || \, || \eta_A X_\alpha^2 r_k || \\
\\ & \leqslant & C \varepsilon_{\omega_0} \left ( \alpha^{-3/2} || \eta_A \partial_x u || + \omega_0^2 || \eta_A u || \right ) \alpha^{-2} \epsilon || \eta_A u || \\
\\ & \leqslant & C \varepsilon_{\omega_0} \alpha^{-7/2} \epsilon \left ( || \eta_A u ||^2 + || \eta_A \partial_x u ||^2 \right ). \end{array} \]
This leads to
\[ |K_4| \leqslant C \varepsilon_{\omega_0} \alpha^{-7/2} \epsilon \left ( 1 + \sqrt{\omega_0} A \epsilon \right ) \left ( || \eta_A u ||^2 + || \eta_A \partial_x u ||^2 \right ). \]
\textit{(About $J_5$.)} We first notice that
\[ \begin{array}{rcl} \partial_\omega ( a_\omega^+ ) ' &=& \displaystyle{-2 \partial_\omega \phi_\omega ' \phi_\omega g'(\phi_\omega^2) - 4 \phi_\omega ' \phi_\omega^2 \partial_\omega \phi_\omega g''(\phi_\omega^2) + 4 \frac{\partial_\omega \phi_\omega '}{\phi_\omega} g(\phi_\omega^2) + 6 \phi_\omega ' \partial_\omega \phi_\omega g'(\phi_\omega^2)} \\
\\ & & \, \, \, \displaystyle{- \, 4 \frac{\partial_\omega \phi_\omega '}{\phi_\omega^3} G(\phi_\omega^2) + 12 \frac{\phi_\omega ' \partial_\omega \phi_\omega}{\phi_\omega^4} G(\phi_\omega^2) - 12 \frac{\phi_\omega ' \partial_\omega \phi_\omega}{\phi_\omega^2} g'(\phi_\omega^2).} \end{array} \]
We recall that $\partial_\omega \phi_\omega = \omega^{-1} \Lambda_\omega$ and we know estimates on $\Lambda_\omega$. More precisely, we recall that $| \partial_\omega \phi_\omega ' | \leqslant C \rho^4$, $| \phi_\omega | \leqslant C \sqrt{\omega_0}$, $|g'(\phi_\omega^2)| \leqslant \varepsilon_{3 \omega_0/2}$, $| \partial_\omega \phi_\omega | \leqslant \frac{C \rho^4}{\sqrt{\omega_0}}$, $| g''(\phi_\omega^2)| \leqslant \frac{\varepsilon_{3 \omega_0/2}}{\phi_\omega^2}$, $| g(\phi_\omega^2) | \leqslant \varepsilon_{3 \omega_0/2} \phi_\omega^2$, $|G (\phi_\omega^2)| \leqslant \varepsilon_{3 \omega_0/2} \phi_\omega^4$ and $|\phi_\omega '| \leqslant C \omega_0$. This gives $| \partial_\omega ( a_\omega^+ )' | \leqslant C \varepsilon_{3 \omega_0 /2} \sqrt{\omega_0} \rho^4$. Thus, integrating this inequality on $[\omega_0 \, , \omega]$, we get
\[ \left | \frac{\Phi_B}{\zeta_B^2} \left ( (a_\omega^+)' - (a_{\omega_0}^+)' \right ) \right | \leqslant C |x| \varepsilon_{3 \omega_0/2} \sqrt{\omega_0} \rho^4 | \omega - \omega_0 | \leqslant C \varepsilon_{3 \omega_0 /2} |\omega - \omega_0| \rho. \]
The same proof holds for $a_\omega^-$ with a minor difference. Indeed, $\partial_\omega (a_\omega^-)'$ involves $g'''$ (not only $G$, $g$, $g'$ and $g''$) and this derivative is not controlled by $\varepsilon_{\omega_0}$. We thus have to introduce $\widetilde{\varepsilon}_{\omega} := \sup\limits_{|s| \leqslant 3 \omega} | s^2 g'''(s)|$. We cannot be sure that $\varepsilon_\omega \leqslant \widetilde{\varepsilon}_\omega$, since $g''(0)$ is possibly not zero (it possibly does not even exist). With the same arguments as $a_{\omega}^+$, we find that
\[ \left | \frac{\Phi_B}{\zeta_B^2} \left ( (a_\omega^+)' - (a_{\omega_0}^+)' \right ) \right | \leqslant C ( \varepsilon_{3 \omega_0 / 2} + \widetilde{\varepsilon}_{3 \omega_0 /2} ) | \omega - \omega_0 | \rho. \]
Using the upper bound $|\omega - \omega_0 | \leqslant \epsilon$, we finally obtain the following estimate:
\[ |J_5| \leqslant C ( \varepsilon_{3 \omega_0 / 2} + \widetilde{\varepsilon}_{3 \omega_0 /2} ) \epsilon \int_{\R} \rho |z|^2. \]
\textit{(About $K_5$.)} This estimate is similar. The first part is easier. Using the estimate $|R_B| \leqslant C \varepsilon_{\omega_0} \rho$ (which is analogous to the estimate on $R_\infty$ given in the proof Proposition 2), we have
\[ | \omega - \omega_0 | \int_{\R} |R_B| \, |z|^2 \leqslant C \varepsilon_{\omega_0} \epsilon \int_{\R} \rho |z|^2. \]
For the second part of $K_5$, similarly as $J_5$ we write that $| \partial_\omega a_\omega^{\pm} | \leqslant C \varepsilon_{3 \omega_0 /2} \rho$ and thus $| a_\omega^\pm - a_{\omega_0}^\pm | \leqslant C \varepsilon_{3 \omega_0 / 2} \epsilon \rho$. Then we get
\[ |K_5| \leqslant C ( \varepsilon_{\omega_0} + \varepsilon_{\omega_0} \varepsilon_{3 \omega_0 /2} ) \epsilon \int_{\R} \rho |z|^2 \leqslant C \varepsilon_{\omega_0} \epsilon \int_{\R} \rho |z|^2. \]
\textit{(Conclusion.)} We first recall from Lemma 7 that
\[ \int_{\R} P_B |z|^2 \geqslant C \varepsilon_{\omega_0} \gamma_B \sqrt{\omega_0} \int_{\R} \rho |z|^2 - \frac{C \varepsilon_{\omega_0} \sqrt{\omega_0}}{\gamma_B} || \partial_x z ||^2 . \]
Let us take $B$ large enough (depending on $\omega_0$) such that $\gamma_B \geqslant \frac{1}{2} \int_{\R} \frac{P_\infty}{\varepsilon_{\omega_0}} \geqslant 10 \max \left ( \frac{C_2}{C} \, , C \right ) \varepsilon_{\omega_0} \sqrt{\omega_0}$. This comes from $(H_2)$. Here, recall that $C_2$ is the constant involved in the control of $K_2$. We obtain
\[ \int_{\R} P_B |z|^2 \geqslant 10 C_2 \varepsilon_{\omega_0}^2 \omega_0 \int_{\R} \rho |z|^2 - \frac{1}{10} || \partial_x z ||^2. \]
First, let us take $\omega_0$ small enough such that
\[ |K_2| \leqslant \frac{1}{100} || \partial_x z ||^2 + C_2 \omega_0 \varepsilon_{\omega_0}^2 \int_{\R} \rho |z|^2. \]
Note that the control on $K_2$ does not imply $A$, $B$, $\alpha$ or $\epsilon$: it only depends on $\omega_0$. The fact that we have the quantity $\varepsilon_{\omega_0}^2 \omega_0$ in front of $\int_{\R} \rho |z|^2$ is crucial. It matches the analogous term in the inequality above given by Lemma 7. \\
\\ Now, we take $B$ large enough so that the previous assumption about $\gamma_B$ holds, and that 
\[ |J_1| \leqslant \frac{\varepsilon_{\omega_0}^2 \omega_0}{100} \int_{\R} \rho |z|^2 \, , \, \, \, |K_1| \leqslant \frac{1}{100} \left [ || \partial_x z ||^2 + C_2 \varepsilon_{\omega_0}^2 \omega_0 \int_{\R} \rho |z|^2 + \frac{1}{A^3 \omega_0^{3/2}} \left ( \alpha^{-4} || \eta_A \partial_x u ||^2 + \omega_0^4 || \eta_A u ||^2 \right ) \right ]. \]
From now on, $B$ is considered as a constant. Now, let us fix $\alpha$ small enough (depending on $\omega_0$ and $B$) such that
\[ |J_3|, |K_3| \leqslant \frac{1}{100} \left ( || \partial_x z ||^2 + C_2 \omega_0 \varepsilon_{\omega_0}^2 \int_{\R} \rho |z|^2 \right ) + \frac{C}{A^3 \omega_0^{3/2}} \left ( \alpha^{-4} || \eta_A \partial_x u ||^2 + \omega_0^4 || \eta_A u ||^2 \right ). \]
From now on, $\alpha$ is considered as a constant. We get
\[ |J_2| \leqslant \frac{C}{A} \left ( || \eta_A \partial_x u ||^2 + \frac{\omega_0^4}{A^2} || \eta_A u ||^2 + \omega_0^5 || \rho^2 u ||^2 \right ). \]
Now, $A$ remains to be fixed but the way we choose $A$ will be given a little bit later. We choose $\epsilon$ small enough (depending on $\omega_0$ and $A$) such that
\[ |J_4| , |K_4| \leqslant \frac{1}{100A} \left ( || \eta_A \partial_x u ||^2 + \frac{\omega_0^4}{A^2} || \eta_A u ||^2 \right ) \, \, \, \, \, \text{and} \, \, \, \, \, |J_5| , |K_5| \leqslant \frac{C_2 \varepsilon_{\omega_0}^2 \omega_0}{100} \int_{\R} \rho |z|^2. \]
All of this lead to
\[ \left | \sum_{j=1}^5 (J_j + K_j) \right | \leqslant 2C_2 \varepsilon_{\omega_0}^2 \omega_0 \int_{\R} \rho |z|^2 + \frac{1}{10} || \partial_x z ||^2 + C \left ( \frac{1}{A^3 \omega_0^{3/2}} + \frac{1}{A} \right ) || \eta_A \partial_x u ||^2 + \frac{C \omega_0^{5/2}}{A^3} || \eta_A u ||^2 + \frac{C \omega_0^5}{A} || \rho^2 u ||^2. \]
Now, we get
\[ \begin{array}{rcl} \dot{\mathcal{J}} + \dot{\mathcal{K}} & \geqslant & \displaystyle{\left ( 2 - \frac{1}{10} - \frac{1}{10} \right ) || \partial_x z ||^2 + C_2 \varepsilon_{\omega_0}^2 \omega_0 \left ( 10-2 \right ) \int_{\R} \rho |z|^2 - C \left ( \frac{1}{A^3 \omega_0^{3/2}} + \frac{1}{A} \right ) || \eta_A \partial_x u ||^2} \\ \\ & & \displaystyle{ \, \, \, \, \, \, \, \, - \, \frac{C \omega_0^{5/2}}{A^3} || \eta_A u ||^2 - \frac{C \omega_0^5}{A} || \rho^2 u ||^2} \\
\\ & \geqslant & \displaystyle{|| \partial_x z ||^2 + C_2 \varepsilon_{\omega_0}^2 \omega_0 \int_{\R} \rho |z|^2 - \frac{C}{A \sqrt{\omega_0}} || \eta_A \partial_x u ||^2 - \frac{C \omega_0^{5/2}}{A^3} || \eta_A u ||^2 - \frac{C \omega_0^5}{A} || \rho^2 u ||^2,} \end{array} \]
where we have noticed that $\frac{1}{A} + \frac{1}{A^3 \omega_0^{3/2}} \leqslant \frac{C}{A \sqrt{\omega_0}}$. Lemma 13 then gives
\[ \begin{array}{rcl} \displaystyle{|| \partial_x z ||^2 + C_2 \varepsilon_{\omega_0}^2 \omega_0 \int_{\R} \rho |z|^2} & \geqslant & \displaystyle{C \varepsilon_{\omega_0}^2 \omega_0 \left ( || \partial_x z ||^2 + \int_{\R} \rho |z|^2 \right )} \\
\\ & \geqslant & \displaystyle{C \varepsilon_{\omega_0}^2 \omega_0 || \rho v ||^2 - \frac{C \varepsilon_{\omega_0}^2}{A^3 \omega_0^{1/2}} || \eta_A \partial_x u ||^2 - \frac{C \varepsilon_{\omega_0}^2 \omega_0^{7/2}}{A^3} || \eta_A u ||^2} \end{array}. \]
Finally we obtain
\[ \dot{\mathcal{J}} + \dot{\mathcal{K}} \geqslant C \varepsilon_{\omega_0}^2 \omega_0 || \rho v ||^2 - \frac{C}{A \sqrt{\omega_0}} || \eta_A \partial_x u ||^2 - \frac{C \omega_0^{5/2}}{A^3} || \eta_A u ||^2 - \frac{C \omega_0^5}{A} || \rho^2 u ||^2. \]
By the definition of $\mathcal{J}$ and the upper bounds $| \Psi_{A,B} | \leqslant C \eta_A^2$ and $| \Psi_{A,B} ' | \leqslant C \eta_A^2$ (recall that $B$ is now a constant), we have, for any $T>0$,
\[ \begin{array}{rcl} | \mathcal{J} (T) | &=& \displaystyle{\left | \int_{\R} v_1 ( 2 \Psi_{A,B} \partial_x v_2 + \Psi_{A,B} ' v_2 ) \right | \, \, \leqslant \, \, C \left ( || \eta_A v(T)||^2 + || \eta_A \partial_x v(T)||^2 \right )} \\
\\ & \leqslant & \displaystyle{C \left ( || \eta_A u(T)||^2 + || \eta_A \partial_x u(T)||^2 \right ) \, \, \leqslant \, \, ||u(T)||_{H^1}^2 \, \, \leqslant \, \, C \epsilon^2.} \end{array} \]
Writing that $|z_k| \leqslant |v_k|$ and $|R_B| \leqslant C \rho^2 \leqslant C \eta_A^2$, the same argument gives $| \mathcal{K} (T) | \leqslant C \epsilon^2$ too. Therefore,
\[ \int_0^T ( \dot{\mathcal{J}} + \dot{\mathcal{K}} ) \, \text{d}t \leqslant | \mathcal{J} (T) | + | \mathcal{K} (T) | + | \mathcal{J} (0) | + | \mathcal{K} (0) | \leqslant C \epsilon^2. \]
Using the inequality on $\dot{\mathcal{J}} + \dot{\mathcal{K}}$ and integrating it on $[0 \, , T]$, we finally obtain:
\[ \varepsilon_{\omega_0}^2 \omega_0 \int_0^T || \rho v ||^2 \, \text{d}t \leqslant C \epsilon^2 + C \int_0^T \left ( \frac{1}{A \sqrt{\omega_0}} || \eta_A \partial_x u ||^2 + \frac{\omega_0^{5/2}}{A^3} || \eta_A u ||^2 + \frac{\omega_0^5}{A} || \rho^2 u ||^2 \right ) \, \text{d}t. \]
This is the result announced. \hfill \qedsymbol

\subsection{Coercivity property and conclusion}
\noindent Now we will need the following coercivity property.

\begin{leftbar}
\noindent \textbf{Proposition 5.} Assume $(H_1)$ and $(H_2)$. We have
\[ \omega_0^2 || \rho^2 u || \leqslant C || \rho v ||. \]
\end{leftbar}

\noindent \textit{Proof.} We follow the exact same proof as in \cite{Ma1}. We need two lemmas to obtain the desired result. First, if $q \in L^2 ( \R )$ satisfies $\langle q \, , \phi_\omega \rangle = \langle q \, , x \phi_\omega \rangle = 0$, then $|| \rho^2 q || \leqslant C \omega_0^{-2} || \rho (X_\alpha^2 S^2 L_+ q) ||$. We follow the proof in \cite{Ma1}. We recall that we know that $| \langle \phi_\omega \, , \Lambda_\omega \rangle | \geqslant C \sqrt{\omega}$. We only have to check that we can write
\[ \begin{array}{rl} & \displaystyle{\frac{q''}{\phi_\omega} = \left ( \frac{q}{\phi_\omega} \right ) '' + (f_3 q)' + f_2 q} \\ \\ \text{and} & \displaystyle{\frac{q''''}{\phi_\omega} = \left ( \frac{q}{\phi_\omega} \right ) '''' + (f_3 q)''' + (f_2 q)'' + (f_1 q)' + f_0 q} \end{array} \]
where $f_j$ are $\mathscr{C}^\infty$ functions (whose expression change from line to line) which satisfy $|f_j(x)| \leqslant C \omega^{-1/2} e^{\sqrt{\omega} |x|}$. This is easily checked thanks to the lower bound $\phi_\omega (x) \geqslant c \sqrt{\omega} \, e^{- \sqrt{\omega} |x|}$. For example, in the first line, $f_2= -2 \frac{\omega}{\phi_\omega} + \phi_\omega - 2 \frac{G(\phi_\omega^2)}{\phi_\omega^3}$. The rest of the proof is entirely identical to the proof of Lemma 11 in \cite{Ma1}. Note that we use the expression and the properties of $I_+$ here. \\
\\ The second lemma we need is the following one: if $q \in L^2 ( \R )$ satisfies $\langle q \, , \Lambda_\omega \rangle = \langle q \, , \phi_\omega ' \rangle = 0$, then $|| \rho^2 q || \leqslant C \omega_0^{-2} || \rho (X_\alpha^2 M_- S^2)q ||$. Here the proof is entirely identical to the proof of Lemma 12 in \cite{Ma1}. There is only an identity at the end of the proof which is different: in our case we have $\phi_\omega '' \phi_\omega - 2 ( \phi_\omega ')^2 = - \omega \phi_\omega^2 + \phi_\omega^2 g(\phi_\omega^2) - 2 G(\phi_\omega^2)$. The rest of the argument is unchanged. Note that we use the expression and the properties of $J_-$ here; that is why hypothesis $(H_2)$ is needed. \hfill \qedsymbol

\noindent \textcolor{white}{a} \\ Now we can conclude the proof of Theorem 2. Using propositions 3, 4 and 5, we obtain
\[ \begin{array}{rcl} \displaystyle{\int_0^T || \rho^2 u ||^2 \, \text{d}t} & \leqslant & \displaystyle{C \omega_0^{-4} \int_0^T || \rho v ||^2 \, \text{d}t} \\
\\ & \leqslant & \displaystyle{C \omega_0^{-5} \varepsilon_{\omega_0}^{-2} \epsilon^2 + C \int_0^T \left ( \frac{\omega_0^{-11/2} \varepsilon_{\omega_0}^{-2}}{A} || \eta_A \partial_x u ||^2 + \frac{\omega_0^{-5/2} \varepsilon_{\omega_0}^{-2}}{A^3} || \eta_A u ||^2 + \frac{\varepsilon_{\omega_0}^{-2}}{A} || \rho^2 u ||^2 \right ) \, \text{d}t} \\
\\ & \leqslant & \displaystyle{C \omega_0^{-5} \varepsilon_{\omega_0}^{-2} \epsilon^2 + \frac{C \omega_0^{-1/2} \varepsilon_{\omega_0}^{-2}}{A} \int_0^T \left ( || \eta_A \partial_x u ||^2 + \frac{1}{A^2} || \eta_A u ||^2 \right ) \, \text{d}t + \frac{C \varepsilon_{\omega_0}^{-2}}{A} \int_0^T || \rho^2 u ||^2 \, \text{d}t} \\
\\ & \leqslant & \displaystyle{C \omega_0^{-5} \varepsilon_{\omega_0}^{-2} \epsilon^2 + \frac{C \omega_0^{-11/2} \varepsilon_{\omega_0}^{-2}}{A} \left ( C \epsilon + C \omega_0 \int_0^T || \rho^2 u ||^2 \, \text{d}t \right ) + \frac{C \varepsilon_{\omega_0}^{-2}}{A} \int_0^T || \rho^2 u ||^2 \, \text{d}t.} \end{array} \]
Since $\omega_0^{-11/2} / A \leqslant \omega_0^{-5}$, we have
\[ \int_0^T || \rho^2 u ||^2 \, \text{d}t \leqslant C \omega_0^{-5} \varepsilon_{\omega_0}^{-2} \epsilon^2 + \frac{C \omega_0^{-9/2} \varepsilon_{\omega_0}^{-2}}{A} \int_0^T || \rho^2 u ||^2 \, \text{d}t. \]
Now we fix $A$. We choose $A$ (depending on $\omega_0$, $B$ and $\alpha$) such that $A>B> \omega_0^{-1/2}$ and $\frac{C \omega_0^{-9/2} \varepsilon_{\omega_0}^{-2}}{A} \leqslant \frac{1}{100}$. This gives
\[ \int_0^T || \rho^2 u ||^2 \, \text{d}t \leqslant C \omega_0^{-5} \varepsilon_{\omega_0}^{-2} \epsilon^2. \]
Using the first virial property, letting $T \to + \infty$ and recalling that $A$ is now a constant, we obtain
\[ \int_0^{+ \infty} \left ( || \eta_A \partial_x u||^2 + || \eta_A u ||^2 + \omega_0 || \rho^2 u ||^2 \right ) \leqslant C \epsilon + C \omega_0^{-4} \varepsilon_{\omega_0}^{-2} \epsilon^2 \leqslant C \omega_0^{-4} \varepsilon_{\omega_0}^{-2} \epsilon^2. \]
Now, we recall the system \eqref{Su} verified by $u$ and we integrate by parts, noticing that $u_2 \partial_x^2 u_1 - u_1 \partial_x^2 u_2 = \partial_x (u_2 \partial_x u_1 - u_1 \partial_x u_2)$:
\[ \begin{array}{rcl} \displaystyle{\frac{\text{d}}{\text{d}t} \left ( \frac{|| \rho^2 u ||^2}{2} \right )} & = & \displaystyle{\int_{\R} \rho^4 ( u_1 \partial_t u_1 + u_2 \partial_t u_2 )} \\
\\ & = & \displaystyle{\int_{\R} ( \rho^4 )' (u_1 \partial_x u_2 - u_2 \partial_x u_1 ) + \int_{\R} 2 \rho^4 u_1 u_2 \phi_\omega^2 ( 1 - g'(\phi_\omega^2))} \\
\\ & & \displaystyle{\, \, + \, \int_{\R} \rho^4 \left ( ( \theta_2 + m_2 - q_2) u_1 - (\theta_1 + m_1 - q_1) u_2 \right ).} \end{array} \]
We write that $| \rho ' | \leqslant C \rho$, so $| (\rho^4)' | \leqslant C \rho^4$. Hence,
\[ \left | \int_{\R} ( \rho^4 )' (u_1 \partial_x u_2 - u_2 \partial_x u_1 ) \right | \leqslant C \int_{\R} \rho^4 \left ( | \partial_x u |^2 + |u|^2 \right ). \]
Another easy bound is the following one (using $| \phi_\omega^2 - \phi_\omega^2 g'(\phi_\omega^2)| \leqslant C$):
\[ \left | \int_{\R} 2 \rho^4 u_1 u_2 \phi_\omega^2 ( 1 - g'(\phi_\omega^2)) \right | \leqslant C || \rho^2 u ||^2. \]
Recalling that $|q_1|,|q_2| \leqslant C \epsilon |u| \leqslant C |u|$, we have
\[ \left | \int_{\R} \rho^4 ( -q_2 u_1 + q_1u_2 ) \right | \leqslant C || \rho^2 u ||^2. \]
Now, using \eqref{orth} and $|x \phi_\omega |, |\phi_\omega|, |\Lambda_\omega| , |\phi_\omega ' | \leqslant C$, we find
\[ | \theta_1 | , | \theta_2 | \leqslant C || \rho^2 u ||^2. \]
On the other hand,
\[ |m_1| \leqslant | \dot{\beta} | \, |x u_1| + | \dot{\gamma} - \omega - \beta^2 | \, |u_1| + | \dot{\sigma} - 2 \beta | \, | \partial_x u_2 | + |\beta| \, | \dot{\sigma} - 2 \beta | \, |u_1| \leqslant C || \rho^2 u ||^2 ( 1 + |x| ) \]
and the same estimate holds for $m_2$. Since $\int_{\R} |x| \rho^4 < + \infty$, we finally obtain that:
\[ \left | \frac{\text{d}}{\text{d}t} || \rho^2 u ||^2 \right | \leqslant C \left ( || \rho^2 \partial_x u ||^2 + || \rho^2 u ||^2 \right ). \]
We recall that $\int_0^{+ \infty} || \rho^2 u ||^2 \, \text{d}t \leqslant C \omega_0^{-5} \varepsilon_{\omega_0}^{-2} \epsilon^2 < \infty$; therefore there exists a sequence $t_n \to + \infty$ such that
\[ || \rho^2 u(t_n) || \, \underset{n \to + \infty}{\longrightarrow} \, 0. \]
Now let us consider $t>0$ and $n$ such that $t_n > t$. We integrate the previous inequality on $[t \, , t_n]$, which gives
\[ || \rho^2 u(t) ||^2 \leqslant || \rho^2 u(t_n) ||^2 + C \int_t^{t_n} \left ( || \rho^2 \partial_x u ||^2 + || \rho^2 u ||^2 \right ) \, \text{d} \tau. \]
Passing to the limit $n \to + \infty$, we get
\[ || \rho^2 u (t)||^2 \leqslant C \int_t^{+ \infty} \left ( || \rho^2 \partial_x u ||^2 + || \rho^2 u ||^2 \right ) \, \text{d} \tau \, \underset{t \to + \infty}{\longrightarrow} \, 0. \]
The previous integral term exists (and converges to $0$ as $t \to + \infty$) because
\[ \int_0^{+ \infty} \left ( || \rho^2 \partial_x u ||^2 + || \rho^2 u ||^2 \right ) \leqslant \int_0^{+ \infty} \left ( || \eta_A \partial_x u ||^2 + || \eta_A u ||^2 \right ) < \infty. \]
Hence we have shown that
\[ || \rho^2 u(t) || \, \underset{t \to + \infty}{\longrightarrow} \, 0. \]
Now, let us take $x,y \in \R$. Using the Cauchy-Schwarz inequality and the basic inequality $| (\rho^2)'| \leqslant C \rho^2$, we write that
\[ \begin{array}{rcl} \rho^2 (x) | u(t \, , x) |^2 & = & \displaystyle{\rho^2 (y) | u(t \, , y) |^2 + \int_x^y \left ( 2 \, \text{Re} \left ( \overline{u(t)} \, \partial_x u(t) \right ) \rho^2 + |u(t)|^2 ( \rho^2 )' \right )} \\
\\ & \leqslant & \displaystyle{\rho^2 (y) | u(t \, , y)|^2 + C ||u(t)||_{H^1 ( \R )} || \rho^2 u(t) ||.} \end{array} \]
We integrate for $y \in [0 \, , 1]$ and use the Cauchy-Schwarz inequality again, as well as \eqref{orbstab}:
\[ \rho^2 (x) | u(t \, , x) |^2 \leqslant \int_{\R} \rho^2 |u(t)|^2 + C ||u(t)||_{H^1 ( \R )} || \rho^2 u (t) || \leqslant C || u(t)||_{H^1 ( \R )} || \rho^2 u(t) || \leqslant C \epsilon || \rho^2 u(t) ||. \]
Henceforth,
\[ \sup_{x \in \R} \rho^2 (x) | u(t \, , x) | \leqslant C \epsilon || \rho^2 u(t) || \, \underset{t \to + \infty}{\longrightarrow} \, 0. \]
This assures that, for any compact $I \subset \R$,
\[ \sup_{x \in I} | u(t \, , x)| \leqslant \frac{1}{\displaystyle{\min_I} (\rho^2)} \, \sup_{x \in \R} \rho^2 (x) |u(t \, , x)|^2 \, \underset{t \to + \infty}{\longrightarrow} \, 0. \]
Now, we recall from \eqref{orth} that $| \dot{\beta} | + | \dot{\omega} | \leqslant C || \rho^2 u ||^2$ thus
\[ \int_0^{+ \infty} | \dot{\beta} | \, \text{d}t + \int_0^{+ \infty} | \dot{\omega} | \, \text{d}t \leqslant C \int_0^{+ \infty} || \rho^2 u ||^2 \, \text{d}t < \infty, \]
which shows that $\omega (t)$ and $\beta (t)$ have finite limits when $t \to + \infty$ (namely respectively $\omega_+$ and $\beta_+$). Letting $t \to + \infty$ in \eqref{orbstab} we find that $| \beta_+ | + | \omega_+ - \omega_0 | \leqslant \epsilon$. Finally, to conclude we write that
\[ | \psi (t \, , x + \sigma (t)) - e^{i \gamma (t)} e^{i \beta_+ x} \phi_{\omega_+} (x) | \leqslant | e^{i \beta (t) x} \phi_{\omega (t)} (x) - e^{i \beta_+ x} \phi_{\omega_+} (x) | + | u(t \, , x) |. \]
First,
\[ \left | e^{i \beta_+ x} \phi_{\omega (t)} (x) - e^{i \beta_+ x} \phi_{\omega_+} (x) \right | = | \phi_{\omega (t)} (x) - \phi_{\omega_+} (x) |  = \left | \int_{\omega_+}^{\omega (t)} \partial_{\widetilde{\omega}} \phi_{\widetilde{\omega}} (x) \, \text{d} \widetilde{\omega} \right | \leqslant \frac{C | \omega (t) - \omega_+ |}{\sqrt{\omega_0}}. \]
This shows that
\[ \sup_{x \in \R} \left | e^{i \beta_+ x} \phi_{\omega (t)} (x) - e^{i \beta_+ x} \phi_{\omega_+} (x) \right | \, \underset{t \to + \infty}{\longrightarrow} \, 0. \]
And on the other hand,
\[ \left | e^{i \beta_+ x} \phi_{\omega (t)} (x) - e^{i \beta(t) x} \phi_{\omega (t)} (x) \right | \leqslant \left | e^{i \beta_+ x} - e^{i \beta (t) x} \right | = 2 \left | \sin \left ( \frac{\beta_+ - \beta (t)}{2} \, x \right ) \right | \]
which shows that, for any compact $I \subset \R$,
\[ \sup_{x \in I} \left | e^{i \beta_+ x} \phi_{\omega (t)} (x) - e^{i \beta(t) x} \phi_{\omega (t)} (x) \right | \leqslant \sup_{x \in I} 2 \left | \sin \left ( \frac{\beta_+ - \beta (t)}{2} \, x \right ) \right | \, \underset{t \to + \infty}{\longrightarrow} \, 0. \]
Gathering those two estimates and the fact that $\displaystyle{\sup_{x \in \R} |u(t \, , x) | \, \underset{t \to + \infty}{\longrightarrow} \, 0}$, we finally obtain that
\[ \sup_{x \in \R} | \psi (t \, , x + \sigma (t)) - e^{i \gamma (t)} e^{i \beta_+ x} \phi_{\omega_+} (x) | \, \underset{t \to + \infty}{\longrightarrow} \, 0, \]
which is the theorem we sought to establish. \hfill \qedsymbol

\end{document}